%% file: gibo.tex
\documentclass[english,11pt]{smfart}
\usepackage{etex}

\usepackage{amsbsy}
\usepackage{amsmath,amsfonts,amssymb,amsthm,mathrsfs,mathtools}
\usepackage{bm}
\usepackage[a4paper,vmargin={3cm,3cm},hmargin={2.5cm,2.5cm}]{geometry}
\usepackage[font= sf, labelfont={sf,bf}, margin=1cm]{caption}
\usepackage{graphicx}
\usepackage{epsfig}%pour les eps
\usepackage{latexsym}%encore des symboles
\usepackage{xcolor}
\usepackage{array,booktabs}
\usepackage{ae,aecompl}
\usepackage{soul,framed}
\usepackage{comment}
\usepackage[pdfpagemode=UseNone,bookmarksopen=false,colorlinks=true,urlcolor=blue,citecolor=blue,citebordercolor=blue,linkcolor=blue]{hyperref}
\usepackage{smfhyperref}
\usepackage[capitalize]{cleveref}
\usepackage{pgfplots}
\pgfplotsset{compat=1.18}

\usepackage{pstricks}
\usepackage{enumerate}
\usepackage{tikz,animate,media9}						%%%%%%%%%%
\usepackage{todonotes}
\usepackage{pifont}
\usepackage{bm,marvosym}
\usepackage{algorithm}
\usepackage{algorithmic}

\usepackage{bbm}

\usepackage{tcolorbox}

\definecolor{aleacolor}{rgb}{0.16,0.59,0.78}

\hypersetup{
breaklinks,
colorlinks=true,
linkcolor= aleacolor,
urlcolor= aleacolor,
citecolor= aleacolor}

\newtheorem{theorem}{Theorem}[section]

\newtheorem{corollary}[theorem]{Corollary}
\newtheorem{proposition}[theorem]{Proposition}
\newtheorem{lemma}[theorem]{Lemma}

\newtheorem{remark}[theorem]{Remark}
\newtheorem{example}[theorem]{Example}

\title{The Critical Cauchy Boundary for Gibbs Partitions}

\author{Fabien Baeriswyl}
\address[Fabien Baeriswyl]{Vienna University of Technology}
\email{fabien.baeriswyl at tuwien.ac.at}
\author{Benedikt Stufler}
\address[Benedikt Stufler]{Vienna University of Technology}
\email{benedikt.stufler at tuwien.ac.at}

\newcommand{\ndN}{\mathbb{N}}
\newcommand{\ndZ}{\mathbb{Z}}

\newcommand{\ndR}{\mathbb{R}}

\newcommand{\diff}{\mathrm{d}}

\newcommand{\Prb}[1]{\mathbb{P} \left(#1\right)}

\newcommand{\Ex}[1]{\mathbb{E}[#1]}

\newcommand{\E}[1]{\mathbb{E} \left[#1\right]}

\newcommand{\Va}[1]{\mathbb{V}[#1]}

\newcommand{\one}{{\mathbbm{1}}}

\newcommand{\convdis}{\,{\buildrel \mathrm{d} \over \longrightarrow}\,}
\newcommand{\convd}{\,{\buildrel \mathrm{d} \over \longrightarrow}\,}

\newcommand{\convp}{\,{\buildrel \mathrm{p} \over \longrightarrow}\,}

\newcommand{\eqdist}{\,{\buildrel \mathrm{d} \over =}\,}

\newcommand{\cA}{\mathcal{A}}

\newcommand{\cC}{\mathcal{C}}
\newcommand{\cD}{\mathcal{D}}

\newcommand{\cF}{\mathcal{F}}
\newcommand{\cG}{\mathcal{G}}

\newcommand{\cM}{\mathcal{M}}
\newcommand{\cN}{\mathcal{N}}

\newcommand{\cS}{\mathcal{S}}
\newcommand{\cT}{\mathcal{T}}

\newcommand{\N}{\mathbb{N}}
\newcommand{\Z}{\mathbb{Z}}

\newcommand{\R}{\mathbb{R}}

\newcommand{\I}[1]{\one_{\left\{#1\right\}}}

\begin{document}
\maketitle

\vspace{-0.5cm}

\begin{abstract}
	We study critical Gibbs partitions whose component law is in the domain of attraction of a spectrally positive stable law of index one. In the dense phase, we prove a local limit theorem for the number of components and functional limits for their sizes. We also describe convergent and divergent condensation regimes, and a phase transition at a boundary in which the dense and condensation contributions compete. 
\end{abstract}

\section{Introduction}
\label{sec:intro}

Gibbs partitions are a general model of random partitions generated according to two weight sequences that bias the number and size of components. Alongside other general models~\cite{arratiaBarbourTavare2000,arratiaBarbourTavare2003,MR2121024} they enable a unified study of a wide range of combinatorial structures that decompose into smaller components, such as connected components of graphs, cycles of permutations and trees in forests. See \cite{MR2245368,berestyckiPitman2007,erlihsonGranovsky2008,MR4780503,arXiv260114166}. In combinatorics, the term composition scheme is often used synonymously \cite{gourdon1998,flajoletSedgewick2009,banderierEtAl2001,arXiv241103930}.

We study critical Gibbs partitions for which the component law belongs to the domain of attraction of a spectrally positive stable law of index one. Spectrally positive refers here to the fact that the associated Lévy process has no negative jumps. Our main interest concerns the number and the size of the components. We describe a dense phase, in which the number of components concentrates around a deterministic sequence that is regularly varying with index $1$, and the largest component has sublinear size. We then uncover a condensation phase, in which the largest component is asymptotically equivalent to the total mass $n$. We distinguish two subregimes of condensation. In the first the number of components converges in distribution to a finite random variable. We refer to this as convergent condensation since the behaviour is identical to the convergent case for random partitions satisfying a conditioning relation~\cite{MR2121024}. In the other subregime the number of components tends to infinity. We refer to this as divergent condensation.

Our main results encompass a local limit theorem for the number of components in the dense phase. In the divergent condensation case, we describe the components (outside of the largest one) with functional limit theorems and point process limits. We describe a phase transition, where we observe an asymptotic mixture behaviour between the dense phase and the condensation phase.  We apply our framework to critical and subcritical Bienaymé--Galton--Watson trees for which the offspring law belongs to the domain of attraction of a Cauchy law. Here we complement the comprehensive results by Kortchemski and Richier~\cite{zbMATH07057468} and Kortchemski and Vetter~\cite{zbMATH08114289} by local limit theorems for the maximal degree.

\subsection{Gibbs partitions and composition schemes}\label{sec:probsetup}

For a detailed introduction to Gibbs partitions, we refer to \cite{MR2245368}. Let $(v_k)_{k \ge 1}$ be a nonnegative sequence weighting the number of components and $(w_j)_{j \ge 0}$ be a nonnegative sequence weighting their individual sizes. We call these the outer, respectively the inner, weight sequences, and we write their generating functions as
\[
	V(z):= \sum_{k \ge 1}v_kz^k, \qquad W(z):= \sum_{j \ge 0}w_jz^j, \qquad U(z):= V(W(z)) = \sum_{n \ge0}u_nz^n.
\]
Their radii of convergence are denoted by $\rho_v$ and $\rho_w$. For sizes $n$ with $0 < u_n < \infty$, we define the ordered component-size vector $P_n := (K_1,\ldots,K_{N_n})$ (of random size) by its distribution
\[
	u_n= \sum_{k \ge 1}v_k \sum_{\substack{n_1+\cdots+n_k= n \\n_1,\ldots,n_k \ge 0}} \prod_{i=1}^kw_{n_i}, \qquad \Prb{P_n=(n_1,\ldots,n_k)} := \frac{v_k}{u_n}\prod_{i=1}^kw_{n_i}.
\]
We call $u_n$ the partition function and $N_n$ the number of components. We write $K_{(1)} \ge K_{(2)} \ge \cdots$ for the ranked component sizes. Unlabelled analogues are discussed in \cite{stufler2020unlabelled,zbMATH07831642,arXiv:2203.15543}.

We mainly consider critical composition schemes, for which
\begin{equation}\label{eq:critical}
	W(\rho_w) = \rho_v.
\end{equation}
For a summary of known subcritical and supercritical regimes see~\cite{MR4780503} and the references therein. Define a random nonnegative integer $X$ and its independent copies $(X_i)_{i \ge 1}$ by its distribution
\begin{equation}\label{eq:lawX}
	\Prb{X = j} := \frac{w_j \rho_w^j}{W(\rho_w)}, \qquad j \ge 0,
\end{equation}
and write $S_0 := 0$ and $S_k := X_1 + \cdots + X_k$ for $k \ge 1$. Multiplying the partition function by $\rho_w^n$, and referring to this multiplicative quantity as the tilted partition function, the definition of $X$ and its partial sums give
\begin{equation}\label{eq:intro-random-sum}
	u_n \rho_w^n = \sum_{k \ge 1} v_k W(\rho_w)^k \Prb{S_k = n},
\end{equation}
while the law of the number of components satisfies
\begin{equation}\label{eq:components}
	\Prb{N_n = k} = \frac{v_k W(\rho_w)^k \Prb{S_k = n}}{u_n \rho_w^n}.
\end{equation}
For every $k \ge 1$ with $\Prb{N_n = k}>0$, the Kolchin representation~\cite{kolchin1986,MR2908619} reads
\begin{equation}\label{eq:KasX}
	\left((K_1,\ldots,K_{N_n})\mid N_n = k\right) \eqdist \left((X_1,\ldots,X_k)\mid S_k = n \right).
\end{equation}
In studying the different regimes, we estimate the contributions to~\eqref{eq:intro-random-sum} from different ranges of $k$. Equation~\eqref{eq:components} tells us where the distribution of $N_n$ concentrates. The Kolchin representation~\eqref{eq:KasX} hence connects the study of the component sizes to that of independent copies of $X$ conditioned on their sum to equal $n$, for those values of $k.$

\subsection{The Cauchy boundary regime}
We assume throughout that the inner weights satisfy
\begin{equation}\label{eq:wasymptotic}
	w_j \sim L_w(j)j^{-2} \rho_w^{-j},
\end{equation}
with $L_w$ a slowly varying function. Hence, $(w_j \rho_w^j)$ is a regularly varying sequence with index $-2$. This asymptotic relation justifies the Cauchy domain terminology for the component law. Note that setting $L := L_w/W(\rho_w)$, Karamata's theorem, recalled in Theorem~\ref{thm:karamata} below, gives
\begin{equation}\label{eq:xtail}
	\Prb{X>x} \sim L(x)x^{-1}.
\end{equation}
By the choice in~\eqref{eq:wasymptotic} and the specification in~\eqref{eq:lawX}, $X$ has a regularly varying tail with index $-1$; it may (or may not) have a finite first moment, depending on the slowly varying term $L$.

Further, let $(a_k)$ be the smooth increasing normalisation constructed in Section~\ref{sec:background}, which satisfies $k \Prb{X > a_k} \rightarrow 1$, and set
\[
	\mu(x) := \E{X;X \le x}, \qquad b_k=k\mu(a_k), \qquad m_n := \min \{k \ge 1:b_k\ge n \}, \qquad c_n := \frac{a_{m_n}}{\mu(a_{m_n})}.
\]
Then, the partial sums $S_k$ of $(X_i)$ satisfy, by the classical stable domain of attraction theorem (see~\cite[Section~1.1]{zbMATH07109856}), \[(S_k-b_k)/a_k \convdis Y, \quad k \rightarrow \infty,\] where $Y$ has the asymmetric Cauchy law specified in~\eqref{eq:chara}. We denote its density by $f$. As shown in Section~\ref{sec:scales}, it holds that $b_{m_n} \sim n$, $c_n = o(m_n)$ and $a_{m_n} = o(n)$. If $\E{X} = \mu < \infty$, then $m_n \sim n/\mu$, while otherwise it holds that $m_n = o(n)$. For example, $L(x) \rightarrow c > 0$ gives $m_n \sim n/(c \log n)$.

Under~\eqref{eq:wasymptotic}, Theorem~\ref{thm:local-big-jump} gives a uniform local limit for the largest summand conditioned on $S_k = s_k$ when $(s_k-b_k)/a_k \rightarrow \infty$.

\subsection{Main results}

For the outer weights, suppose that $v_k \sim L_v(k)k^{-b}W(\rho_w)^{-k}$, where $L_v$ is a slowly varying function and $b > 1$ is the outer weights exponent. Write
\[
	H(x) := \sum_{1 \le k \le x}k v_k W(\rho_w)^k, \qquad H(\infty) := \lim_{x \rightarrow \infty}H(x).
\]
\subsubsection*{Cauchy dense phase: $1 < b < 2$}

Theorem~\ref{thm:numco} gives the following asymptotic for the tilted partition function:
\[
	u_n \rho_w^n \sim \frac{H(m_n)-H(m_n-1)}{b_{m_n}} \sim \frac{L_v(m_n)m_n^{-b}}{\mu(a_{m_n})}.
\]
Furthermore, the number of components satisfies a uniform local limit theorem: for every $\delta > 0,$
\begin{equation*}
	\sup_{k \ge \delta m_n} \left|c_n \Prb{N_n = k} -f\left(-\frac{k-m_n}{c_n}\right)\right|\longrightarrow 0.
\end{equation*}
Hence $(N_n-m_n)/c_n \convdis-Y$ and $N_n/m_n \convp 1$. The ``dense'' terminology refers to concentration around $m_n$, which is not necessarily linear in $n$. 

Lemma~\ref{lem:tv} compares an initial segment of the component-size vector with a corresponding vector of independent copies of $X$. In turn, this gives convergence of the centered and rescaled partial-sum process in the Skorokhod $J_1$ topology, as stated in Theorem~\ref{thm:functional}. A direct consequence is Corollary~\ref{cor:largest}, which gives a point-process convergence of the component sizes
\[
	\sum_{\substack{1 \le i \le N_n \\ K_i > 0}} \delta_{K_i/a_{m_n}} \convdis \cN := \sum_{j \ge 1}\delta_{J_j}, \qquad \nu(\diff x) := x^{-2} \diff x,
\]
where $\cN$ is a Poisson point process with intensity measure $\nu$. Its points represent the limiting rescaled component sizes. Listed in decreasing order, these satisfy $(J_j)_{j \ge 1} \eqdist (\Gamma_j^{-1})_{j \ge 1}$, where $(\Gamma_j)$ are the arrival times of a unit-rate Poisson process. In particular, $K_{(1)} = O_p(a_{m_n}) = o_p(n)$. It also covers convergence of vectors of the largest component sizes. The same Poisson order statistics arise for the largest remaining
outdegrees after deletion of the maximum in critical Cauchy BGW
trees by Kortchemski and Richier~\cite[Theorem~3]{zbMATH07057468}.

\subsubsection*{Convergent condensation: $b > 2$}

In this case, $H(\infty) < \infty$, and known limits~\cite{stufler2018,MR4780503} (recalled in Proposition~\ref{thm:condensation}) apply. In particular, the tilted partition function satisfies
\[
    u_n \rho_w^n \sim \Prb{X =n} H(\infty) \sim L(n)n^{-2} H(\infty).
\] 
A typically unique largest component contains asymptotically all the mass $n$, i.e. $K_{(1)}/n \convp 1.$ We call it the giant. The remaining mass after removing it, $n-K_{(1)},$ is called the residual mass. In the convergent case, both the residual mass and the number of components have finite distributional limits, i.e.
\[
	N_n \convdis \widehat{N}, \qquad n-K_{(1)} \convdis S_{\widehat{N}-1}, \qquad \Prb{\widehat{N} =k} = \frac{k v_k W(\rho_w)^k}{H(\infty)}.
\]
The limiting number of components $\widehat{N}$ is independent of the sequence $(X_i)_{i \ge 1}$. After deleting the largest component, the remaining sizes, in their original order, are approximated in total variation by $(X_1, \ldots, X_{\widehat{N}-1})$.

\subsubsection*{Phase transition: $b = 2$}
We set
\[
	\cC_n^{\cD}  := \frac{H(m_n)-H(m_n-1)}{b_{m_n}} = \frac{v_{m_n} W(\rho_w)^{m_n}}{\mu(a_{m_n})}, \quad \cC_n^{\cG}  := \Prb{X = n}H(m_n), \quad \Theta_n := \cC_n^{\cG} /\cC_n^{\cD} . 
\]
Theorem~\ref{thm:b2} establishes the competition
\[
	u_n \rho_w^n \sim \cC_n^{\cD} +\cC_n^{\cG} ,
\]
between the two main regimes
\[
    \qquad \Prb{\cD_n} := \Prb{N_n \ge m_n/2} \sim\frac{1}{1+\Theta_n}, \qquad \Prb{\cG_n} := \Prb{N_n < m_n/2} \sim\frac{\Theta_n}{1+\Theta_n}. 
\]
The term $\cC_n^{\cD} $ approximates the contribution to~\eqref{eq:intro-random-sum} from numbers of components $k = m_n + O(c_n)$. The term $\cC_n^{\cG} $ approximates the contribution from $k = o(m_n)$, and where one component accounts for asymptotically almost all the total mass $n$. This is the so-called single-big-jump principle. The remaining components are approximated, when removing the largest one from consideration, by independent copies of $X$. Conditionally on $\cD_n$, the dense-phase limits hold true. Conditionally on $\cG_n$, a unique giant emerges, and the remaining component sizes, in their original order, are approximated in total variation by $(X_1, \ldots, X_{\widehat{N}_n-1})$. Here, $\widehat{N}_n$ is independent of $(X_i)_{i \ge 1}$ and has distribution
\[
 \Prb{\widehat{N}_n = k} = \frac{k v_k W(\rho_w)^k}{H(m_n)}, \qquad 1 \le k \le m_n,
\]
as in the construction used in~\eqref{eq:defNhatboundary}.

If $H(\infty) < \infty$, then, conditionally on $\cG_n$, the number of components $N_n$ converges in distribution to $\widehat{N}$ as defined in the convergent condensation case. If $H(\infty) = \infty$, then, conditionally on $\cG_n$, $N_n \convp \infty$, $N_n/m_n \convp 0$ and $K_{(1)}/n \convp 1$. We call this last regime divergent condensation. In this case, Theorem~\ref{thm:divergent condensationtv} and the conditional approximation in Theorem~\ref{thm:b2} give
\[
	\left(\frac{H(N_n \wedge m_n)}{H(m_n)}, \frac{n-K_{(1)}-b_{N_n}}{a_{N_n}}, \sum_{\substack{2 \le i \le N_n \\ K_{(i)} > 0}} \delta_{K_{(i)}/a_{N_n}} \right) \convdis (U, Y, \cN).
\]
In the above, $U$ is uniform on $[0,1]$ and independent of the point process $\cN$. The variable $Y$ is determined by the same Poisson points: we sum the points larger than some threshold $\varepsilon$, subtract a compensator $\log(1/\varepsilon)$, and let $\varepsilon \rightarrow 0$. This subtraction defines the compensated sum constructed in Lemma~\ref{lem:meso-poisson-cauchy}. Therefore, the limiting residual-mass fluctuation and secondary sizes are constructed from the same point process. Note that these limits hold conditionally on $\cG_n$. When $\Theta_n \rightarrow \infty$, $\Prb{\cG_n} \rightarrow 1$, so the same limits hold without having to condition.

We call all components other than the giant the secondary components. Their ranked sizes are $K_{(2)},K_{(3)}, \ldots$. Deleting the giant and the following $r$ largest components leaves the mass $n - \sum_{i = 1}^{r+1}K_{(i)}$. This operation is referred to as trimming. For each fixed $r \ge 1$, Corollary~\ref{cor:meso-ranked} gives the joint limit of $(K_{(2)}/a_{N_n}, \ldots, K_{(r+1)}/a_{N_n})$ and
\[
 \frac{n - \sum_{i = 1}^{r+1}K_{(i)}-b_{N_n}}{a_{N_n}}.
\]
The residual-mass limit also holds for $r = 0$, in which case only the giant is removed. Corollary~\ref{cor:meso-functionals} gives convergence of the entire ranked sequence of rescaled, secondary sizes in $\ell^q$ for $q > 1$, and consequently of the sums of their $q$th powers. Proposition~\ref{pro:meso-functional} gives joint convergence of the centered and rescaled partial sums of the remaining component sizes and the point process encoding their rescaled sizes and positions. Corollary~\ref{cor:meso-conditional-llt} gives a local limit for the residual mass conditional on $N_n$. Finally, Example~\ref{ex:logs} illustrates these regimes for logarithmic weights.

A summary of the regimes and their main features is given in Table~\ref{tab:intro-regimes}. 

\begin{table}[htbp]
	\centering
	\small
	\renewcommand{\arraystretch}{1.3}
    \setlength{\tabcolsep}{5pt}
	\begin{tabular}{@{}
		>{\raggedright\arraybackslash}p{.27\linewidth}
		>{\raggedright\arraybackslash}p{.25\linewidth}
		>{\raggedright\arraybackslash}p{.24\linewidth}
		>{\raggedright\arraybackslash}p{.16\linewidth}@{}}
		\toprule
		Range and law & Number of components & Largest component & Source\\
		\midrule
		\multicolumn{4}{@{}l}{\textit{Dense phase}}\\[.2ex]
		$1<b<2$
			& $N_n = m_n + O_p(c_n)$
			& $K_{(1)} = o_p(n)$
			& This paper\\
		$b=2$, given $\cD_n$
			& $N_n = m_n + O_p(c_n)$
			& $K_{(1)} = o_p(n)$
			& This paper\\

		\addlinespace[.7ex]
		\multicolumn{4}{@{}l}{\textit{Convergent condensation}}\\[.2ex]
		$b=2$, given $\cG_n$, $H(\infty) < \infty$
			& $N_n \convdis \widehat{N}$
			& $n-K_{(1)} = O_p(1)$
			& This paper\\
		$b > 2$
			& $N_n \convdis \widehat{N}$
			& $n-K_{(1)} = O_p(1)$
			& Known~\cite{MR4780503}\\

		\addlinespace[.7ex]
		\multicolumn{4}{@{}l}{\textit{Divergent condensation}}\\[.2ex]
		$b=2$, given $\cG_n$, $H(\infty) = \infty$
			& $N_n \convp \infty$, $N_n/m_n \convp 0$
			& $K_{(1)}/n \convp 1$
			& This paper\\
		\bottomrule
	\end{tabular}
	\caption{Regimes at inner weight index $-2$.
	The dense limits follow from Theorem~\ref{thm:numco} and Corollary~\ref{cor:largest}. The conditional $b = 2$ limits follow from Theorem~\ref{thm:b2}, with $\Prb{\cD_n}$ and $\Prb{\cG_n}$ determined by $\Theta_n$. The known $b > 2$ case is recalled in Proposition~\ref{thm:condensation}.}
	\label{tab:intro-regimes}
\end{table}

\subsection{Applications}

Kortchemski and Richier~\cite{zbMATH07057468} studied critical Bienaymé--Galton--Watson (BGW) trees conditioned to be large with an offspring distribution belonging to the Cauchy domain of attraction. Their main result describes a condensation phenomenon, with a typically unique vertex exhibiting a large degree that is regularly varying with index $1$. This is in sharp contrast with other regimes for critical trees, since condensation has otherwise only been observed in the subcritical regime~\cite{MR2764126,MR2908619,MR3335012,MR4144886}. Kortchemski and Vetter~\cite{zbMATH08114289} described a condensation in the case of a Cauchy domain of attraction in the subcritical regime, completing the phase diagram.

For both critical and subcritical BGW trees with an offspring law $\xi$ satisfying $\Prb{\xi = k} \sim L_\xi(k) k^{-2}$ with $L_\xi$ a slowly varying function, Theorem~\ref{thm:bgwmax} below strengthens the distributional limits for the maximal degree of \cite{zbMATH07057468,zbMATH08114289} to a local limit theorem. This answers a question specifically raised in~\cite[Sec. 4.1]{zbMATH08114289}. A motivation specifically for providing local limit theorems for the maximal degree of random trees is the recently discovered common-core transfer method~\cite{stunonseparxiv}, which exploits them for proving scaling limits of random planar maps.

The Cauchy dense regime of Gibbs partitions also appears in various combinatorial settings. We are going to provide applications in a companion paper.

\subsection{Organisation of the paper}
Section~\ref{sec:background} collects the probabilistic tools and Cauchy normalisations. Section~\ref{sec:scales} studies the centering and auxiliary sequences and proves a local limit theorem for the largest summand under the assumption that the sum is itself large. Section~\ref{sec:dense} treats the dense phase, Section~\ref{sec:condensation} concerns condensation, and Section~\ref{sec:phase} studies the phase transition at the outer weight exponent $b = 2$. Section~\ref{sec:bgw} applies the results and tools from previous sections to conditioned BGW trees and forests.

\subsection{Notation}\label{sec:not}

We  let $\ndN = \{1, 2, \ldots\}$ and $\ndN_0 = \{0, 1, 2, \ldots\}$ denote the sets of positive and nonnegative integers. For each $n \in \ndN_0$, we set $[n] := \{1, \ldots, n \}$, so that~$[0]= \emptyset$. Random variables are defined on a common probability space, with its associated probability measure denoted by $\mathbb{P}$. Unspecified limits are taken as $n$ tends to infinity. We use $\convdis$ and $\convp$ to denote convergence in distribution and probability. Equality in distribution is denoted by $\eqdist$.

For an event $A$ and a random variable $X$ we write $\Ex{X; A} = \Ex{X \one_A}$, with $\one_A$ denoting the indicator variable of $A$. We let $O_p(1)$ denote an unspecified random variable $X_n$ of a stochastically bounded sequence $(X_n)_{n \ge 1}$, and write $o_p(1)$ for an unspecified random variable $Y_n$ with $Y_n \convp 0$. For probability measures $\lambda,\lambda'$ defined on the same measurable space, we write
\[
	d_{\mathrm{TV}}(\lambda,\lambda') = \sup_B| \lambda(B)-\lambda'(B)|,
\]
where the supremum is over measurable sets. By a slight abuse of notation, we also use $d_{\mathrm{TV}}(A,B)$ for the distance between the laws of random elements $A,B$. We sometimes write $\lambda \approx_{\mathrm{TV}} \lambda'$ to mean that the total-variation distance between $\lambda$ and $\lambda'$ goes to 0. We refer to e.g. \cite[Chapter~3]{thorisson2000} for background on total variation.

For a nonempty finite vector $\mathbf{x} = (x_1,\ldots,x_k) \in \ndN_0^k$, let $i$ be the position of its first maximal coordinate. Define
\begin{equation}\label{eq:deletion-map}
    \tau(\mathbf{x}) := (x_1, \ldots, x_{i-1}, x_{i+1}, \ldots, x_k).
\end{equation}
$\tau$ deletes the first maximal coordinate and preserves the original order of all remaining coordinates. In particular, $\tau((x_1)) = \emptyset$. We also set $\tau(\emptyset) = \emptyset$.

We write $\mathbb{D}([0,1], \ndR)$ for the space of real-valued càdlàg functions (right-continuous with left limits), and we endow it with Skorokhod $J_1$ topology. For a textbook treatment, see~\cite[Chapter~12]{whitt2002}. For a random element $Z$ of $\mathbb{D}([0,1], \ndR)$, we write $\Delta Z(t) := Z(t) - Z(t-)$ for its jump at time $t \in ]0,1]$. For a locally compact, second countable space $E$, we let $\cM_p(E)$ denote the space of Radon point measures (i.e., measures that are finite on compact sets), and we endow it with the vague topology. We say that $\xi_n \rightarrow\xi$ vaguely if $\int f \diff \xi_n \rightarrow\int f \diff \xi$ for every continuous, compactly supported function $f$ on $E$; see~\cite[Section~3.4]{resnick1987}. $\delta_x$ denotes the Dirac point mass at $x$. On $]0,\infty]$, we use the topology induced by $x \mapsto 1/x$, with $1/\infty = 0$. Whenever necessary, finite ranked sequences are extended by zeros. For $q \ge 1$, $\ell^q$ is endowed with its norm $\lvert \mathbf{x} \rvert_q=(\sum_{j \ge 1}|x_j|^q)^{1/q}$.

For a power series $f(z) \in \ndR[[z]]$, we write $[z^n]f(z)$ for the coefficient of $z^n$. We use the classical, deterministic asymptotic notation $O$, $o$ and $\Theta$. For positive sequences, $u_n \asymp v_n$ means that, for $n$ large enough, $u_n/v_n$ is bounded above and bounded away from zero. We write $x \wedge y := \min\{x,y\}$, $x\vee y := \max\{x,y\}$ and $x_+ := \max\{x,0\}$. 

\section{Probabilistic framework, background and regular variation}\label{sec:background}

\subsection*{Total variation}

On a countable space, total variation is half the $\ell^1$ distance between the probability mass functions. We refer to~\cite[Section~4.1]{levinPeres2017} for a comprehensive background. A Markov kernel $Q$ assigns to each input $x$ a probability measure $Q(x,\cdot)$, with $x \mapsto Q(x,B)$ measurable for every measurable set $B$. It describes the same possibly random operation applied to either input law. For probability measures $\lambda$ and $\lambda'$, applying this operation gives
\begin{equation}\label{eq:tv-kernel-contract}
 d_{\mathrm{TV}}(\lambda Q,\lambda'Q) \le d_{\mathrm{TV}}(\lambda,\lambda'),
\end{equation}
where $(\lambda Q)(B) = \int Q(x,B)\,\lambda(\diff x)$. To see this, note that $0 \le Q(x,B) \le 1$, so the variational characterisation of total variation gives
\[
 |(\lambda Q)(B)-(\lambda'Q)(B)| = \left| \int Q(x,B)(\lambda-\lambda')(\diff x)\right| \le d_{\mathrm{TV}}(\lambda,\lambda').
\]
Taking the supremum over $B$ proves~\eqref{eq:tv-kernel-contract}. Note that the kernel may depend on $n$. In particular, a measurable map $T_n$ corresponds to $Q_n(x,B) = \one_{\{T_n(x)\in B\}}$, and hence
\[
 d_{\mathrm{TV}}\big(T_n(A_n),T_n(B_n)\big) \le d_{\mathrm{TV}}(A_n,B_n).
\]
Therefore, if $d_{\mathrm{TV}}(A_n,B_n) \rightarrow 0$, every distributional limit obtained for $T_n(B_n)$ also holds for $T_n(A_n)$.

A coupling of two probability laws is a pair of random elements defined jointly with those marginal laws. For any coupling of $A$ and $B$ on a countable space, $d_{\mathrm{TV}}(A,B) \le \Prb{A\ne B}$, and a coupling attaining equality exists; see~\cite[Proposition~4.7]{levinPeres2017}.

\subsection*{Cauchy normalisations: the essentials}

We use the distribution $X$ from~\eqref{eq:lawX} and its tail asymptotic~\eqref{eq:xtail}. We now specify a smooth increasing normalisation $(a_{\cdot})$ satisfying~\eqref{eq:defan}.

The function $x \mapsto 1/L(x)$ is asymptotically equivalent to a smooth, positive function $A^*(x)$ by \cite[Proposition 1.2.5]{Mikosch1999regular}. The function $A(x) = \int_{1}^{x} A^*(t) \diff t$ is strictly increasing and, by Karamata's theorem~\ref{thm:karamata}, it is asymptotically equivalent to $x A^*(x)$, as $x \rightarrow\infty$. Choose $A^*$ on a bounded interval so that $A(j_0) < 1$ for some $j_0 \ge 1$ with $\Prb{X = j_0} > 0$. We let $a(x)$ denote the inverse function of $A(x)$ and, for integer values, we set $a_n = a(n)$. Then $a(1) > j_0$, so $\mu(a_n) > 0$ for every $n \ge 1$. Hence $a(x)$ is regularly varying at infinity with index $1$ (see \cite[Theorem 1.5.12]{bgt1989}) and
\begin{align}\label{eq:defan}
	n \Prb{X>a_n} \rightarrow 1.
\end{align}
Equivalently,
\begin{align}\label{eq:defan2}
	a_n \sim n L(a_n).
\end{align}
We therefore construct $a(\cdot)$ as the inverse of the strictly increasing function $A(\cdot)$. Hence $a(\cdot)$ increases strictly as well. Moreover, since $A'(x) =A^*(x)$, the (ordinary) inverse function theorem gives
\begin{equation}\label{eq:aderivative}
	a'(x) = \frac{1}{A'(a(x))} = \frac{1}{A^*(a(x))} \sim L(a(x)), \qquad x \rightarrow \infty.
\end{equation}
The truncated first moment of $X$ is denoted by
\[
	\mu(x) := \Ex{X; X \le x}
\]
and
\begin{equation}\label{eq:defbn}
	b_n := n \mu(a_n).
\end{equation}
Clearly, the centering $(b_{\cdot})$ is nondecreasing and even eventually increasing. The stable limit theorem reads
\begin{align}\label{eq:snbnclt}
	\frac{S_n - b_n}{a_n} \convdis Y,
\end{align}
where $Y$ is an asymmetric Cauchy random variable with characteristic function given, for $t > 0,$ by
\begin{align}\label{eq:chara}
	\E{\exp(itY)} = \exp(-t (\pi / 2 - i(1-\gamma_E)) - it \log t), 
\end{align}
with $\gamma_E \approx 0.577\ldots$ denoting Euler's constant. Its density $f$ is analytic on $\ndR$, see \cite[Theorem 4.5]{janson2022stabledistributions}, and admits the representation
\begin{align*}
	f(x) = \frac{1}{\pi} \sum_{n = 0}^\infty (-1)^n  \frac{(x - (1-\gamma_E))^n}{n!} \int_{0}^\infty \exp(- t \log t) \sin(\pi t)t^n \,\diff t, \qquad x \in \ndR.
\end{align*}
This follows from~\cite[Theorem~4.4]{janson2022stabledistributions} by simply shifting the $\delta$-parameter therein to $1-\gamma_E$. The asymptotic given in Equation~\eqref{eq:wasymptotic} implies that, eventually, $w_n > 0$, so that the support of $X$ contains consecutive integers and has maximal lattice span $1$. This justifies stating the following local limit theorem~\cite[Theorem~4.2.1]{MR0322926}
\begin{align}\label{eq:llt}
	\sup_{k \in \ndZ} \left\vert a_n \Prb{S_n = k} - f\left(\frac{k-b_n}{a_n}\right)  \right\vert \rightarrow  0.
\end{align}
The characteristic function in~\eqref{eq:chara} is integrable. Hence, $f$ is bounded and uniformly continuous, and, as $M \rightarrow \infty,$
\begin{equation}\label{eq:ftail}
	\sup_{|x| > M} f(x) \longrightarrow 0.
\end{equation}
In particular,~\eqref{eq:llt} gives $\sup_z \Prb{S_h=z} \le C/a_h$ for all sufficiently large $h$. Since $(a_{\cdot})$ is increasing and regularly varying with index $1$, for every fixed $\theta > 0$ there is $C_\theta < \infty$ such that, for all sufficiently large $k,$
\begin{equation}\label{eq:common-concentration}
	\sup_{\substack{h \in \ndN \\ h \ge \theta k}} \sup_{z \in \ndZ} \Prb{S_h=z} \le C_\theta a_k^{-1}.
\end{equation}
By Berger's local large deviation bounds~\cite[Theorem~2.4]{zbMATH07109856}, there exists a constant $C > 0$ such that, for all $n \ge 1$ and all $k \in \ndZ$ with $|k| \ge a_n$,
\begin{align}\label{eq:bedev}
	\Prb{S_n - \lfloor b_n \rfloor =  k} \le C n L(|k|) |k|^{-2}.
\end{align}
For negative deviations, apply the same theorem to $-X$. The required local tail condition therein is automatic in our case because $X$ is nonnegative. 

The asymptotic part of the same theorem~\cite[Theorem~2.4, Equation~(2.10)]{zbMATH07109856} applies because the local mass satisfies $\Prb{X = j} \sim L(j)j^{-2}$ by~\eqref{eq:lawX} and~\eqref{eq:wasymptotic}. Therefore, if $k = k_n > 0$ and $k/a_n \rightarrow \infty$ as $n \rightarrow \infty$, then it holds that
\begin{align}\label{eq:beloc}
	\Prb{S_n - \lfloor b_n \rfloor =  k} \sim n L(k) k^{-2}.
\end{align}

\subsection*{Regular variation.}
We collect the following regular-variation results, which we use extensively. See~\cite{bgt1989} for a textbook treatment.
\begin{theorem}[Uniform convergence theorem, Theorem 1.5.2 in \cite{bgt1989}]\label{thm:uct}
	Let $f: ]0,\infty[ \rightarrow ]0,\infty[$ be regularly varying with index $\alpha \in \ndR$. Then, for every $0 < a \le b < \infty$,
	\begin{equation*}
		\sup_{\lambda \in [a,b]} \left| \frac{f(\lambda x)}{f(x)}-\lambda^\alpha \right| \longrightarrow 0, \quad x \rightarrow \infty.
	\end{equation*}
\end{theorem}

\begin{theorem}[Potter bounds, Theorem 1.5.6 in \cite{bgt1989}]\label{thm:pb}
	Let $f:]0, \infty[ \rightarrow ]0, \infty[$ be regularly varying with index $\alpha \in \R$. For every $C > 1$ and $\delta > 0$, there exists $x_0 = x_0(f,C,\delta) > 0$ such that, for all $x,y \ge x_0$,
	\begin{equation}\label{eq:pb}
		\frac{f(y)}{f(x)} \le C \max \left\{\left(\frac{y}{x}\right)^{\alpha+\delta}, \left(\frac{y}{x}\right)^{\alpha-\delta} \right\}.
	\end{equation}
\end{theorem}

\begin{theorem}[Karamata's theorem, Theorem 1.5.11 in \cite{bgt1989}]\label{thm:karamata}
	Let $f:]0, \infty[ \rightarrow ]0, \infty[$ be regularly varying with index $\alpha \in \R$ and locally bounded on some interval $[x_0, \infty[$.
	\begin{itemize}
		\item If $\beta>-\alpha-1$, then
			\[
				\int_{x_0}^{x} t^\beta f(t)\,dt \sim \frac{x^{\beta+1}f(x)}{\beta+\alpha+1}, \quad x \rightarrow \infty.
			\]
		\item If $\beta<-\alpha-1$, then
			\[
				\int_x^\infty t^\beta f(t)\,dt \sim -\frac{x^{\beta+1}f(x)}{\beta+\alpha+1}, \quad x \rightarrow \infty.
			\]
	\end{itemize}
\end{theorem}
\begin{remark}
	The discrete statements of Karamata's theorem~\ref{thm:karamata} follow by applying the integral version to the function $x \mapsto u_{\lfloor x \rfloor}$. To see this, suppose that $(u_{\cdot})$ is a regularly varying sequence with index $\alpha$, and let
	\[
		f(x) := u_{\lfloor x \rfloor}, \text{ for } x \ge 1.
	\]
	Then, for every fixed $\lambda > 0$, and because $\lfloor \lambda x \rfloor / \lfloor x \rfloor \rightarrow \lambda \in ]0, \infty[$, the (sequential form of the) uniform convergence theorem~\ref{thm:uct} gives
	\[
		\frac{f(\lambda x)}{f(x)} = \frac{u_{\lfloor \lambda x \rfloor}}{u_{\lfloor x \rfloor}} \rightarrow \lambda^{\alpha}, \text{ as } x \rightarrow \infty.
	\]
	Now, for $\alpha > -1,$ by the first statement of Theorem~\ref{thm:karamata},
	\[
		\sum_{k = 1}^{n} u_k = \int_{1}^{n+1} f(x) dx  \sim \frac{(n+1) f(n+1)}{\alpha+1} = \frac{(n+1) u_{n+1}}{\alpha+1} \sim \frac{nu_n}{\alpha+1}.
	\]
	A similar conclusion holds for the second statement where $\alpha < -1.$
\end{remark}

\begin{remark}
	We use Theorem~\ref{thm:uct} in the following form: for a regularly varying function $f$ of index $\alpha$ and fixed $0 < \lambda < \Lambda < \infty$,
	\begin{equation}\label{eq:common-rv-window}
		\sup_{\lambda x \le y \le \Lambda x} \left| \frac{f(y)}{f(x)}-\left(\frac{y}{x}\right)^\alpha\right| \longrightarrow 0.
	\end{equation}
	Therefore, if $\varepsilon_x \rightarrow 0$, then
	\begin{equation}\label{eq:common-rv-shrinking}
		\sup_{|y/x-1| \le \varepsilon_x} \left| \frac{f(y)}{f(x)}-1\right| \longrightarrow 0.
	\end{equation}
	When the ratio of arguments tends to zero or infinity, we use the Potter bounds~\eqref{eq:pb} instead.
\end{remark}

\section{Properties of the Cauchy normalisations and local estimates}
\label{sec:scales}
We first study the growth of the centering sequence $b_k$ and its generalised inverse $m(q)$. We then prove a local limit theorem for the largest summand conditioned on a large value of the sum.

\subsection{Growth rate of the centering sequence}\label{sec:centering} The centering $b_n = n \mu(a_n)$ depends on whether $X$ has finite mean. In particular,
\[
	\E{X} < \infty \quad \text{if and only if} \quad \sum_{n \ge 1} n^{-1}L(n) < \infty.
\]
We next determine the growth of $b_k$ in the finite and infinite-mean cases, and then study its generalised inverse $m_n = \min\{k \ge 1: b_k \ge n\}$, which gives the typical number of components in the dense phase.

\subsubsection*{The infinite mean case}

Suppose that~$\E{X} = \infty$. We set
\begin{align*}
	\ell(x) = \int_{1}^x t^{-1} L(t) \diff t.
\end{align*}
A by-product of Karamata's theorem~\ref{thm:karamata}, as formulated in~\cite[Rem. 1.2.7]{Mikosch1999regular}, is that $\ell(x)$ is slowly varying, $\ell(x) \rightarrow \infty$ as $x \rightarrow \infty$, and
\begin{align}\label{eq:qll}
	\ell(x) / L(x) \rightarrow \infty, \qquad x \rightarrow \infty.
\end{align}
We have
\begin{align}\label{eq:mul}
	\mu(x) \sim \ell(x), \qquad x \rightarrow \infty,
\end{align}
and from~\eqref{eq:defan}
\begin{align}\label{eq:defbninfinite}
	b_n \sim  n \ell(a_n).
\end{align}
This implies that $(b_{\cdot})$ is regularly varying with index $1$. Moreover, using~\eqref{eq:defan},~\eqref{eq:xtail} and~\eqref{eq:qll} it follows that
\begin{align}\label{eq:infinitemean}
	a_n \sim n L(a_n) = o(b_n).
\end{align}

\subsubsection*{The finite mean case}

Suppose that~$\mu := \E{X}<\infty$. We set
\begin{align}\label{eq:ellstar}
	\ell^*(x) = \int_{x}^\infty t^{-1} L(t)\,\diff t.
\end{align}
Another by-product of Karamata's theorem~\ref{thm:karamata}, again given in~\cite[Rem. 1.2.7]{Mikosch1999regular}, is that $\ell^*(x)$ is slowly varying,  $\ell^*(x) \rightarrow 0$ and $L(x) \rightarrow 0$ as $x \rightarrow \infty$, and
\begin{align}\label{eq:ellstarl}
	\ell^*(x) / L(x) \rightarrow \infty, \qquad x \rightarrow \infty.
\end{align}
Moreover,
\begin{align}\label{eq:ellstarasymp}
	\mu - \mu(n) \sim \ell^*(n)
\end{align}
and
\begin{align*}
	b_n \sim n \mu .
\end{align*}
Similarly to the infinite mean case, the normalisation is negligible relative to the centering, i.e.
\begin{align}
	\label{eq:finitemean}
	a_n \sim n L(a_n) =  o(n) = o(b_n).
\end{align}

\subsection{Auxiliary sequences}
\label{sec:auxiliary}

To estimate the terms in~\eqref{eq:intro-random-sum} with $k = m_n + O(c_n)$, we need precise estimates for the inverse of the centering sequence $b_k$. Since $(b_k)$ is regularly varying with index $1$ and tends to infinity, define, for every real $q > 0$,
\begin{equation}\label{eq:defmn}
	m(q) := \min \{k \ge 1: b_k \ge q\}.
\end{equation}
The function $m(\cdot)$ is nondecreasing and tends to infinity. Define further
\begin{equation}\label{eq:defcn}
	c(q) := \frac{a_{m(q)}}{\mu(a_{m(q)})}.
\end{equation}
For integer targets $n \in \N$, we use the abbreviations
\[
	m_n :=m(n), \quad c_n := c(n).
\]

Note that $\mu$ is slowly varying: in the finite-mean case, it converges to $\mu$, whereas in the infinite-mean case, $\mu(x) \sim \ell(x)$ by~\eqref{eq:mul}, and $\ell$ is slowly varying as argued above.

\begin{proposition}\label{pro:scale}
	As $q \rightarrow \infty$,
	\[
		c(q) \longrightarrow \infty, \qquad c(q) = o(m(q)).
	\]
	In particular,
	\[
		c_n \longrightarrow \infty, \qquad c_n= o(m_n).
	\]
\end{proposition}

\begin{proof}
	Since $\mu$ is slowly varying, the function $x/\mu(x)$ is regularly varying with index $1$ and therefore tends to infinity. As $a_n \rightarrow \infty$, it follows that
	\begin{equation}\label{eq:muan}
		\frac{a_n}{\mu(a_n)} \longrightarrow \infty.
	\end{equation}
	Since $m(q) \rightarrow \infty$, the above implies, by composition, that $c(q) \rightarrow \infty$.

	Furthermore,
	\begin{equation}\label{eq:omun}
		L(x) = o(\mu(x)), \qquad x \rightarrow \infty.
	\end{equation}
	In the finite-mean case, $\mu(x) \rightarrow \mu$ and $L(x) \rightarrow 0$ (as argued around Equation~\eqref{eq:ellstarl}). In the infinite-mean case, it follows from $L(x) = o(\ell(x))$ (Equation~\eqref{eq:qll}) and $\mu(x) \sim \ell(x)$. Consequently, by~\eqref{eq:defan2},
	\[
		\frac{c(q)}{m(q)} = \frac{a_{m(q)}}{m(q)\mu(a_{m(q)})} \sim \frac{L(a_{m(q)})}{\mu(a_{m(q)})} \longrightarrow 0.
	\]
\end{proof}

We state these estimates for real $q > 0$, because we later use them for residual sums of the form $q = n - \sum_i x_i$.
\begin{lemma}
	\label{lem:seq}
	Let $m(q)$ and $c(q)$ be defined by~\eqref{eq:defmn} and~\eqref{eq:defcn}. As $q \rightarrow \infty,$ the following statements hold.

	\begin{itemize}
		\item[(a)] Let $u(q) \ge 1$ satisfy $u(q) = o(m(q))$. Uniformly over integers $h$ with $|h| \le u(q)$,
			\begin{equation}\label{eq:adiff}
				|a_{m(q)+h}-a_{m(q)}| = O\left(|h|L(a_{m(q)})\right),
			\end{equation}
			\begin{equation}\label{eq:mudiff}
				|\mu(a_{m(q)+h})-\mu(a_{m(q)})| = \left(|h|L(a_{m(q)})+1\right) O(1/m(q)),
			\end{equation}
			and
			\begin{equation}\label{eq:bdiff}
				b_{m(q)+h}-b_{m(q)} = h\mu(a_{m(q)+h}) + O\left(|h|L(a_{m(q)}) +1\right).
			\end{equation}

		\item[(b)] The inverse relation satisfies
			\begin{equation}\label{eq:bdiff2}
				0 \le b_{m(q)}-q \le b_{m(q)}-b_{m(q)-1} = o(a_{m(q)}).
			\end{equation}
			Consequently,
			\begin{equation}\label{eq:bdiff2co}
				b_{m(q)} \sim q
			\end{equation}
			and
			\begin{equation}\label{eq:amo}
				a_{m(q)} = o(q).
			\end{equation}

		\item[(c)] It holds that
			\begin{equation}\label{eq:la}
				L(q) = o(\mu(a_{m(q)})).
			\end{equation}

		\item[(d)] For every $k \in \ndN$,
			\begin{equation}\label{eq:nice1}
				|b_k-b_{m(q)}| \ge |k-m(q)| \mu(a_{m(q)}).
			\end{equation}

		\item[(e)] The functions $m(q)$ and $c(q)$ are regularly varying with index $1$. Moreover, if $(q_j)$ is a sequence of positive real numbers and $(\ell_j)$ is a sequence of positive integers such that
			\[
				q_j \rightarrow \infty, \qquad \ell_j \rightarrow \infty, \qquad q_j \sim b_{\ell_j},
			\]
			then
			\begin{equation}\label{eq:inversestability}
				m(q_j) \sim \ell_j.
			\end{equation}

		\item[(f)] Fix $M>0$ and, for $|x| \le M$, set
			\[
				k_q(x) := m(q) + \lfloor xc(q)\rfloor.
			\]
			Then, uniformly for $|x| \le M$,
			\[
				\frac{a_{k_q(x)}}{a_{m(q)}} \longrightarrow 1, \qquad \frac{\mu(a_{k_q(x)})}{\mu(a_{m(q)})} \longrightarrow 1,
			\]
			and
			\begin{equation}\label{eq:argument}
				\frac{q-b_{k_q(x)}}{a_{k_q(x)}} \longrightarrow -x.
			\end{equation}
	\end{itemize}
\end{lemma}

\begin{proof}
	\textit{Proof of (a).} Uniformly for $|h| \le u(q)$, we have $m(q)+h \sim m(q)$. By~\eqref{eq:aderivative} and~\eqref{eq:common-rv-shrinking}, the derivative $a'(x)$ is $O(L(a_{m(q)}))$ uniformly between $m(q)$ and $m(q)+h$. The mean value inequality proves~\eqref{eq:adiff}. Moreover,~\eqref{eq:defan2} gives
	\begin{equation}\label{eq:cla}
		|h|L(a_{m(q)}) = O(u(q)L(a_{m(q)})) = o(a_{m(q)}).
	\end{equation}
	Hence, $a_{m(q)+h}/a_{m(q)} \rightarrow 1$ uniformly. For every integer $i$ between $\min(a_{m(q)},a_{m(q)+h})$ and $\max(a_{m(q)},a_{m(q)+h})$, the asymptotic of~\eqref{eq:lawX} and~\eqref{eq:common-rv-shrinking} give
	\[
		i\Prb{X =i} = O\left(\frac{L(a_{m(q)})}{a_{m(q)}}\right).
	\]
	By~\eqref{eq:adiff} and~\eqref{eq:defan2},
	\begin{align*}
		|\mu(a_{m(q)+h})-\mu(a_{m(q)})| \le C\left(|a_{m(q)+h}-a_{m(q)}|+1\right)\frac{L(a_{m(q)})}{a_{m(q)}} = \left(|h|L(a_{m(q)})+1\right)O(1/m(q)).
	\end{align*}
	Finally,
	\[
		b_{m(q)+h}-b_{m(q)} = h\mu(a_{m(q)+h}) + m(q)\left(\mu(a_{m(q)+h})-\mu(a_{m(q)})\right),
	\]
	which proves~\eqref{eq:bdiff}.

	\textit{Proof of (b).} By definition,
	\[
		b_{m(q)} \ge q>b_{m(q)-1}.
	\]
	Using~\eqref{eq:bdiff} with $h=-1$ gives
	\begin{equation}\label{eq:ineqbq}
		0 \le b_{m(q)}-q \le b_{m(q)}-b_{m(q)-1} = \mu(a_{m(q)-1}) + O(L(a_{m(q)})+1).
	\end{equation}
	The proof of Proposition~\ref{pro:scale} now yields
	\[
		\mu(a_{m(q)-1}) = o(a_{m(q)}), \qquad L(a_{m(q)}) = o(a_{m(q)}),
	\]
	which then proves~\eqref{eq:bdiff2}. Since $a_n = o(b_n)$ by~\eqref{eq:finitemean} and~\eqref{eq:infinitemean}, it follows, by dividing both sides of~\eqref{eq:ineqbq} by $b_{m(q)}$, and using $a_{m(q)} = o(b_{m(q)}),$ that $b_{m(q)} \sim q$ and $a_{m(q)} = o(q)$.

	\textit{Proof of (c).} By~\eqref{eq:defan2}, ~\eqref{eq:bdiff2co}, the Potter bounds~\eqref{eq:pb}, and~\eqref{eq:amo},
	\begin{align*}
		\frac{L(q)}{\mu(a_{m(q)})} &= \frac{L(q)}{L(a_{m(q)})} \frac{L(a_{m(q)})}{\mu(a_{m(q)})} \\
		&\sim \frac{L(q)}{L(a_{m(q)})} \frac{a_{m(q)}}{m(q)\mu(a_{m(q)})} \sim \frac{L(q)}{L(a_{m(q)})} \frac{a_{m(q)}}{q} \\
		&\le O\left( \max \left\{ \left(\frac{q}{a_{m(q)}}\right)^{1/2}, \left(\frac{a_{m(q)}}{q}\right)^{1/2} \right\} \frac{a_{m(q)}}{q}\right) = o(1).
	\end{align*}

	\textit{Proof of (d).} If $k > m(q)$, the definition of $b_k$ and monotonicity of $\mu(a_k)$ give
	\begin{align*}
		b_k-b_{m(q)} &= (k-m(q))\mu(a_k) + m(q)\left(\mu(a_k)-\mu(a_{m(q)})\right) \ge (k-m(q))\mu(a_{m(q)}).
	\end{align*}
	If $k < m(q)$, similar computations yield
	\begin{align*}
		b_{m(q)}-b_k = (m(q)-k)\mu(a_{m(q)}) + k \left(\mu(a_{m(q)}) -\mu(a_k)\right) \ge (m(q)-k)\mu(a_{m(q)}),
	\end{align*}
	which together prove~\eqref{eq:nice1}.

	\textit{Proof of (e).} Since $(b_{\cdot})$ is eventually increasing and regularly varying with index $1$, the inversion theorem for regularly varying functions (see \cite[Theorem 1.5.12]{bgt1989}) implies that its generalised inverse $m(q)$ is regularly varying with index $1$. Hence $a_{m(q)}$ is regularly varying with index $1$, whereas $\mu(a_{m(q)})$ is slowly varying. It follows that $c(q)$ is regularly varying with index $1$.

	Now suppose that $q_j \sim b_{\ell_j}$. By~\eqref{eq:bdiff2co} and this assumption,
	\[
		\frac{b_{m(q_j)}}{b_{\ell_j}} = \frac{b_{m(q_j)}}{q_j} \frac{q_j}{b_{\ell_j}} \longrightarrow 1.
	\]
	If $m(q_j)/\ell_j$ failed to converge to $1$, then, along a subsequence, for some $\varepsilon > 0$,
	\[
		\text{ either } m(q_j)\ge(1+\varepsilon)\ell_j \quad \text{or} \quad  m(q_j)\le(1-\varepsilon) \ell_j.
	\]
	$(b_{\cdot})$ being monotone and regularly varying with index $1$, this would imply
	\[
		\liminf_{j \rightarrow \infty} \frac{b_{m(q_j)}}{b_{\ell_j}} \ge 1 +\varepsilon \quad \text{or} \quad   \limsup_{j \rightarrow \infty} \frac{b_{m(q_j)}}{b_{\ell_j}} \le 1 -\varepsilon,
	\]
	which would be a contradiction, hence~\eqref{eq:inversestability} follows.

	\textit{Proof of (f).} Let $h_q(x) = \lfloor xc(q) \rfloor$. From~\eqref{eq:adiff} and~\eqref{eq:omun},
	\[
		\sup_{|x| \le M} \left| \frac{a_{m(q)+h_q(x)}}{a_{m(q)}}-1 \right| \longrightarrow 0,
	\]
	while slow variation of $\mu$ and the uniform convergence theorem~\ref{thm:uct} give
	\[
		\sup_{|x| \le M} \left| \frac{\mu(a_{m(q)+h_q(x)})}{\mu(a_{m(q)})}-1 \right| \longrightarrow 0.
	\]
	Moreover,~\eqref{eq:bdiff} and~\eqref{eq:cla} give, uniformly for $|x| \le M$,
	\begin{align*}
		\frac{b_{m(q)+h_q(x)}-b_{m(q)} }{a_{m(q)} } &= \frac{h_q(x)}{c(q)} \frac{\mu(a_{m(q)+h_q(x)})}{\mu(a_{m(q)})} + o(1) \\
		&= x + o(1)
	\end{align*}
	Combining this with~\eqref{eq:bdiff2} and $a_{m(q)+h_q(x)} \sim a_{m(q)}$ proves~\eqref{eq:argument}.
\end{proof}
For later use, note that adjacent normalisations satisfy
\begin{equation}\label{eq:mesocenterincrement}
	a_{k-1} \sim a_k, \qquad b_k-b_{k-1} = o(a_k).
\end{equation}
The first assertion follows from regular variation of $a_{\cdot}$. For the second, use Lemma \ref{lem:seq}(b) with $q = b_k$ therein; as argued above, $b$ is eventually strictly increasing, which gives $m(b_k) = k$.

Also for later use, the discrete version of Karamata~\ref{thm:karamata} gives, for $p > 1,$
\[
\E{X^p; X \le x}\sim L(x)x^{p-1}/(p-1).
\]
Hence, for every fixed $p > 1$ and $\varepsilon > 0$,
\begin{equation}\label{eq:common-truncated-power}
	\frac{k}{a_k^p} \E{X^p; X \le \varepsilon a_k} \longrightarrow \frac{\varepsilon^{p-1}}{p-1}.
\end{equation}
Using~\eqref{eq:common-rv-window} to compare $L(j)$ with $L(a_k)$ in the sum defining $\mu(\lambda a_k)-\mu(a_k)$ gives, for every fixed $\lambda > 0$,
\begin{equation}\label{eq:common-truncated-first}
	\frac{k}{a_k} \left(\mu(\lambda a_k) - \mu(a_k)\right) \longrightarrow \log\lambda.
\end{equation}

The following example gives explicit scales when $L(n) \sim c \log^r n$.
\begin{example}
	\label{ex:addme}
	Suppose that, for constants $c>0$ and $r \in \ndR$,
	\[
		L(n) \sim c \log^{r} n.
	\]
	Then
	\begin{align}\label{eq:example-an}
		a_n \sim cn \log^{r} n.
	\end{align}
	The remaining auxiliary sequences admit the asymptotic relations
	\[
		b_n \sim 
        \begin{cases} 
            \displaystyle \frac{c}{r+1}n \log^{r+1}n, & r > -1, \\[1.2ex]
			\displaystyle c n \log \log n, & r = -1, \\[1.2ex]
			\displaystyle n \mu , & r < -1, 
        \end{cases}
	\]
	where $\mu= \E{X} < \infty$ in the case $r < -1$. Consequently,
	\[
		m_n \sim 
        \begin{cases} 
            \displaystyle \frac{r+1}{c}\frac{n}{\log^{r+1}n}, & r > -1, \\[1.2ex]
			\displaystyle \frac{n}{c \log \log n}, & r = -1, \\[1.2ex]
			\displaystyle \frac{n}{\mu}, & r < -1, 
        \end{cases}
	\]
	and
	\[
		c_n \sim 
        \begin{cases} 
            \displaystyle \frac{(r+1)^2}{c}\frac{n}{\log^{r+2}n}, & r > -1,\\[1.4ex]
			\displaystyle \frac{n}{c \log n(\log \log n)^2}, & r = -1,\\[1.4ex]
			\displaystyle \frac{c}{\mu^2}n \log^{r}n, & r < -1. 
        \end{cases}
	\]
\end{example}
\begin{proof}
	Regular variation of $(a_n)$ gives $\log a_n \sim \log n$, so~\eqref{eq:defan2} yields
	\[
		a_n \sim n L(a_n) \sim cn \log^{r}n,
	\]
	which proves~\eqref{eq:example-an}. Integration of $L(x)/x$, together with~\eqref{eq:mul} and~\eqref{eq:ellstarasymp}, gives
	\[
		\mu(a_n) \sim
		\begin{cases}
			\dfrac{c}{r+1}\log^{r+1}n, & r > -1,\\
			c \log \log n, & r = -1,\\
			\mu, & r < -1.
		\end{cases}
	\]
	In particular, for $r > -1$,
	\begin{equation}\label{eq:muanrgm1}
		\mu(a_n) \sim \frac{c}{r+1} \log^{r+1}n.
	\end{equation}
	For $r < -1$, the tail integral also gives
	\begin{equation}\label{eq:mudiffn}
		\mu - \mu(a_n) \sim \frac{c}{-r-1} \log^{r+1}n,
	\end{equation}
	and, hence,
	\[
		b_n = n \mu-\frac{c}{-r-1} n \log^{r+1}n (1+o(1)).
	\]
	Multiplication by $n$ proves the claimed equivalents for $b_n$. Let $\widetilde{m}_n$ be the proposed equivalent for $m_n$ in the corresponding case. In each case, $\log \widetilde{m}_n \sim \log n$ and $b_{\lfloor\widetilde{m}_n \rfloor} \sim n$. Lemma~\ref{lem:seq}(e) gives $m_n \sim \widetilde{m}_n$. Finally, substitute these equivalents into $c_n = a_{m_n}/\mu(a_{m_n})$ to obtain the three formulas for $c_n$.
\end{proof}

\subsection{A local limit for the largest summand}
\label{sec:local-big-jump}

We consider $S_k$ conditioned on a value $s_k$ such that $(s_k-b_k)/a_k \rightarrow \infty$. Berger, Birkner and Yuan~\cite[Corollary~3.5]{bby24} give an independent-sample approximation after deleting the largest summand, together with a distributional limit for that summand. Also, for  The theorem below strengthens the latter conclusion to a local limit uniform over all integer values of the maximum. It only requires the local asymptotic relation $\Prb{X = j} \sim L(j)j^{-2}$. Recall from~\eqref{eq:infinitemean} and~\eqref{eq:finitemean} that $a_n = o(b_n)$. 

\begin{theorem}[Local limit for the largest summand]\label{thm:local-big-jump}

	Suppose that~\eqref{eq:wasymptotic} holds, so that $X$ satisfies $\Prb{X = j} \sim L(j)j^{-2}$ as $j \rightarrow \infty$. Keep the normalisations $a_k, b_k$ constructed earlier. Let $(s_k)$ be a sequence of integers such that
	\[
		d_k:= s_k - b_k, \qquad \frac{d_k}{a_k} \longrightarrow \infty.
	\]
	Denote by $M_k := \max_{1 \le i \le k} X_i$. Then
	\begin{equation}\label{eq:local-big-jump-global}
		\sup_{j \in \ndZ} \left| a_k \Prb{M_k = j \mid S_k = s_k} -f\left(\frac{d_k-j}{a_k}\right) \right| \longrightarrow 0.
	\end{equation}
\end{theorem}

On $\{S_k = s_k\}$, the identity $M_k = s_k-(S_k-M_k)$ shows that the fluctuations of the largest summand are the negative of those of the remaining sum. This explains the reflected argument $(d_k-j)/a_k$ in~\eqref{eq:local-big-jump-global}.

\begin{proof}
	The local large-deviation estimate~\eqref{eq:beloc} gives
	\begin{equation}\label{eq:local-bj-normaliser}
		\Prb{S_k = s_k} = \Prb{S_k -\lfloor b_k \rfloor = \lceil d_k \rceil} \sim k \Prb{X = \lfloor d_k \rfloor}.
	\end{equation}
	The rounding difference in~\eqref{eq:local-bj-normaliser} is negligible by~\eqref{eq:common-rv-shrinking}. We use the common concentration estimate~\eqref{eq:common-concentration}, for all sufficiently large $k,$ in the form
	\begin{equation}\label{eq:local-bj-concentration}
		\sup_{k/3\le h\le k}\sup_{z\in \ndZ}\Prb{S_h=z} \le Ca_k^{-1}.
	\end{equation}
	Let $|j - d_k| \le D a_k$, for some fixed $D$. First, note that
	\begin{align*}
		k \Prb{X = j} \Prb{S_{k-1} = s_k-j, M_{k-1} < j} &\le \Prb{S_k = s_k, M_k = j} \\
		&\le k \Prb{X = j} \Prb{S_{k-1} = s_k - j}.
	\end{align*}
	The lower bound holds by using independence of the $(X_i)$ and requiring one specific coordinate (we have $k$ many) to equal $j$ while the others are smaller. The upper bound holds by requiring at least one coordinate to equal $j$ and applying a union bound. Now, subtracting the middle term from the upper bound yields
	\begin{align*}
		0 &\le k \Prb{X = j} \Prb{S_{k-1} = s_k-j} - \Prb{S_k = s_k, M_k = j} \\
		&\le k \Prb{X = j} \Prb{S_{k-1} = s_k-j, M_{k-1} \ge j}.
	\end{align*}
	Another union bound and \eqref{eq:local-bj-concentration} give, uniformly for $\lvert j-d_k \rvert < D a_k$,
	\begin{align*}
		\Prb{S_{k-1} = s_k-j, M_{k-1} \ge j} &\le \sum_{i=1}^{k-1} \Prb{S_{k-1} = s_k -j, X_i \ge j} \\
		& \le (k-1) \Prb{X \ge j} \sup_{z} \Prb{S_{k-2} = z} \\
		& \le k \Prb{X\ge j} Ca_k^{-1} = o(a_k^{-1}).
	\end{align*}
	For the uniformity, recalling that by~\eqref{eq:xtail} and~\eqref{eq:defan2}, it holds that $k \Prb{X \ge a_k} \rightarrow 1$, for any $0 < \eta < 1$, the Potter bounds~\eqref{eq:pb} give
	\[
		k \Prb{X \ge d_k/2} \le k \Prb{X \ge a_k} \frac{\Prb{X \ge d_k/2}}{\Prb{X \ge a_k}} \le C_{\eta} (d_k/a_k)^{-1+\eta} \longrightarrow 0.
	\]
	Moreover, the uniform convergence theorem~\ref{thm:uct} gives, uniformly for $|j-d_k| \le D a_k$,
	\begin{equation}\label{eq:Xratjdk}
		\Prb{X = j}/\Prb{X = \lfloor d_k \rfloor} \longrightarrow 1.
	\end{equation}
	Using $a_{k-1} \sim a_k$ and $b_k - b_{k-1} = o(a_k)$ from~\eqref{eq:mesocenterincrement}, the local limit theorem~\eqref{eq:llt} gives, uniformly for $|j-d_k| \le D a_k$,
	\[
		a_k \Prb{S_{k-1} = s_k-j} = \frac{a_k}{a_{k-1}} f\left(\frac{s_k-j-b_{k-1}}{a_{k-1}}\right)+o(1) = f\left(\frac{d_k-j}{a_k}\right) + o(1).
	\]
	By~\eqref{eq:local-bj-normaliser} and the uniform ratio~\eqref{eq:Xratjdk}, $k \Prb{X = j}/\Prb{S_k = s_k} \rightarrow 1$ uniformly in the same window. Collecting all of the above, and uniformly for $|j-d_k| \le D a_k,$
    \begin{align*}
        a_k \Prb{M_k = j \mid S_k = s_k} &= \frac{a_k \Prb{S_k = s_k, M_k = j}}{\Prb{S_k = s_k}} \\
			&= \frac{a_k \left(k \Prb{X = j} \left(\Prb{S_{k-1} = s_k - j} + o(a_k^{-1})\right) \right)}{\Prb{S_k = s_k}} \\
            &= \frac{k \Prb{X = j}}{\Prb{S_k = s_k}} \left(a_k \Prb{S_{k-1} = s_k-j} + o(1)\right)\\
			&=f\left(\frac{d_k-j}{a_k}\right) + o(1).
    \end{align*}
	This proves~\eqref{eq:local-big-jump-global} uniformly for $|j-d_k| \le D a_k$, for every fixed $D > 0$.

	We turn to values of $j$ outside of the above range. A union bound over the coordinates attaining the maximum, independence, and~\eqref{eq:local-bj-normaliser} give
	\begin{align}\label{eq:local-bj-marking}
		\Prb{M_k = j \mid S_k = s_k} &\le \frac{k \Prb{X = j}}{\Prb{S_k = s_k}} \Prb{S_{k-1} = s_k-j, M_{k-1} \le j} \\
		&\le C\frac{\Prb{X = j}}{\Prb{X = \lfloor d_k\rfloor}} \Prb{S_{k-1} = s_k-j, M_{k-1} \le j}\notag.
	\end{align}
	For $j \ge d_k/2$, the Potter bounds~\eqref{eq:pb}, applied to the sequence $(\Prb{X = \cdot})$ of index $-2$, bound the ratio uniformly by a constant. Moreover, by~\eqref{eq:mesocenterincrement}, for each fixed $D > 0$ and all sufficiently large $k$, and for $|j-d_k| > D a_k,$
	\[
		\left|\frac{s_k-j-b_{k-1}}{a_{k-1}}\right| \ge \frac{D a_k-|b_k-b_{k-1}|}{a_{k-1}} \ge D/2.
	\]
	Bounding from above by dropping the maximum restriction in~\eqref{eq:local-bj-marking}, and appealing to the uniform local limit theorem~\eqref{eq:llt}, yields
	\[
		\limsup_{k \rightarrow \infty} \sup_{\substack{j \in \ndN_0, j \ge d_k/2\\ |j-d_k| > D a_k}} \left|a_k \Prb{M_k = j \mid S_k = s_k} - f\left(\frac{d_k-j}{a_k}\right) \right| \le C \sup_{|x| \ge D/2}f(x)\underset{D \rightarrow \infty}{\longrightarrow} 0.
	\]
	Next fix $E > 1$ and consider $Ea_k \le j < d_k/2$. Let $r_{k,j} := s_k-j - \lfloor b_{k-1} \rfloor$.
	By~\eqref{eq:mesocenterincrement}, uniformly in this range, it eventually holds that $d_k/2 \le r_{k,j} \le 2d_k$. Equation~\eqref{eq:bedev} and~\eqref{eq:common-rv-window} yield $\Prb{S_{k-1} = s_k-j} \le Ck L(d_k) d_k^{-2}$.
	Together with~\eqref{eq:local-bj-marking},~\eqref{eq:defan2} and the Potter bounds~\eqref{eq:pb}, this gives, for any fixed $0 < \eta < 1$,
	\[
		a_k \Prb{M_k = j \mid S_k = s_k} \le C a_k k \Prb{X = j} \le C(j/a_k)^{-2+\eta}.
	\]
	The supremum over this range is therefore at most $C E^{-2+\eta}$.

    It remains to control values of the maximum with $j \le E a_k$. The restriction on $M_{k-1}$ in~\eqref{eq:local-bj-marking} must be kept in this range. For $0 \le u \le 1$, the inequality $e^u \le 1+u+(e-2)u^2$ gives, applied to $u = X/j$,
    \[
		\E{e^{X/j}; X \le j} \le 1 - \Prb{X > j} + \frac{{\mu(j)}}{j} + (e-2)\frac{\E{X^2; X \le j}}{j^2}.
	\]
	By the discrete analogue of Karamata's theorem~\ref{thm:karamata} and~\eqref{eq:xtail},
	\[
		\frac{\E{X^2; X \le j}}{j^2} \sim \Prb{X>j} \quad \text{and} \quad e - 2 < 1.
	\]
    Since $\E{X^2; X \le j}/j^2 \sim \Prb{X>j}$ and $e-2 < 1$, choose $j_0 \ge 1$ such that, for $j \ge j_0,$
    \[
        (e-2)\frac{\E{X^2; X \le j}}{j^2} \le \Prb{X>j}.
    \]
    Substituting this into the above bound, and using the basic inequality $1+x \le e^x$ gives, for $j \ge j_0$,
    \begin{equation}\label{eq:local-bj-truncated-mgf}
        \E{e^{X/j}; X \le j} \le 1+\frac{\mu(j)}{j} \le \exp\left(\frac{\mu(j)}{j}\right).
    \end{equation}
	By~\eqref{eq:common-truncated-first}, for fixed $E > 1$,
	\[
		\frac{k}{a_k} \left(\mu(Ea_k)-\mu(a_k)\right) \longrightarrow \log E.
	\]
	Since $\mu$ is nondecreasing, $(\mu(j)-\mu(a_k))_+ = 0$ for $j \le a_k$. For $a_k < j \le Ea_k$, monotonicity of $\mu$ and the bound $1/j \le 1/a_k$ give
	\[
		\sup_{j_0 \le j \le E a_k} \frac{k}{j} \left(\mu(j)-\mu(a_k) \right)_+ \le \frac{k}{a_k} \left(\mu(E a_k)-\mu(a_k) \right) = \log E  + o(1).
	\]
	Applying the exponential Markov inequality with \eqref{eq:local-bj-truncated-mgf} yields, using the independence of the $(X_i)$, for $ h \le k$, and uniformly $j_0 \le j \le Ea_k$, for some constant $C_E > 0,$
	\begin{equation}\label{eq:local-bj-truncated-tail}
		\Prb{S_h \ge h \mu(a_k) + d_k/4, M_h \le j} \le \exp\left( -\frac{d_k}{4j} + \frac{k}{j} \left(\mu(j) - \mu(a_{k})\right)_+ \right) \le C_E \exp\big(-d_k/(4j)\big).
	\end{equation}

	Split the $k-1$ summands into two independent blocks (and denote their respective sums by $S_{h_1}$ and $S_{h_2}$) of lengths
	\[
		h_1 = \lfloor (k-1)/2 \rfloor, \qquad h_2 = k-1-h_1.
	\]
	For all sufficiently large $k$, both $h_1$ and $h_2$ belong to $[k/3, k]$. For $j \le E a_k = o(d_k)$, eventually
	\[
		s_k-j-(k-1) \mu(a_k) = d_k+ \mu(a_k) -j \ge d_k - j \ge d_k/2.
	\]
	On the event $\{S_{k-1} = s_k-j\}$, the two block sums satisfy
	\[
		(S_{h_1}-h_1\mu(a_k)) + (S_{h_2} - h_2\mu(a_k)) = d_k+\mu(a_k)-j\ge d_k/2. 
	\]
	At least one block satisfies $S_{h_i} \ge h_i\mu(a_k)+d_k/4$. Apply~\eqref{eq:local-bj-truncated-tail} to that block and~\eqref{eq:local-bj-concentration} to the other. Convolution and independence give, for $j_0 \le j \le E a_k$,
	\begin{align}\label{eq:local-bj-truncated-local}
		\Prb{S_{k-1} = s_k-j,\ M_{k-1} \le j} &\le \sum_{i=1}^2 \Prb{S_{h_i} \ge h_i \mu(a_k)+d_k/4,\ M_{h_i} \le j} \sup_{z \in \ndZ} \Prb{S_{h_{3-i}} = z} \notag\\
		&\quad\le a_k^{-1} C_E \exp\left(-d_k/(4j)\right).
	\end{align}
	Using \eqref{eq:local-bj-marking}, \eqref{eq:local-bj-truncated-local}, and the Potter bounds~\eqref{eq:pb} for the ratio involving $X$ (choosing a higher $j_0$ if needed), we obtain
	\[
		\sup_{j_0 \le j\le E a_k} a_k \Prb{M_k = j \mid S_k = s_k} \le C_E \sup_{y \ge d_k/(Ea_k)} y^{2+\eta} e^{-y/4} \longrightarrow 0.
	\]
	For $0 \le j < j_0$, choose $t > 0$ so small that $\E{e^{tX}; X \le j_0} \le 1$. Such a choice is possible because $\Prb{X \le j_0} < 1$ and the truncated exponential moment is continuous at $t = 0$. The same argument by splitting into two blocks, using an exponential Markov inequality with this fixed $t$, gives
	\[
		\sup_{0 \le j < j_0} a_k \Prb{M_k = j \mid S_k = s_k} \le \frac{C e^{-t d_k/4}}{\Prb{X = \lfloor d_k \rfloor}} \longrightarrow 0.
	\]
	The last limit follows by regular variation of $X$, and in particular that of the mass $\Prb{X = j}$ with index $-2$.

	The density term tends uniformly to zero for $j < d_k/2$ by~\eqref{eq:ftail}. Letting first $k \rightarrow \infty$, then $E \rightarrow \infty$ and finally $D \rightarrow \infty$, this proves the claim for $j \ge 0$. For $j<0$, the conditional probability is zero and $(d_k-j)/a_k \ge d_k/a_k \rightarrow \infty$, so~\eqref{eq:ftail} extends the estimate to all $j \in \ndZ$.
\end{proof}

\section{The Cauchy dense phase}
\label{sec:dense}

In this section, we consider the Gibbs partition model with outer weight sequence satisfying
\begin{align}\label{eq:vasymptotic-dense}
	v_n \sim L_v(n) n^{-b} W(\rho_w)^{-n},
\end{align}
for a slowly varying function $L_v$ and an exponent $b$ such that $1 < b < 2$.

We first bound $\Prb{S_k = n}$ when the number of summands satisfies $k > 2m_n$. This estimate depends only on the inner weights and might (and in fact will) also be used when the outer exponent is $b = 2$.

\begin{lemma}[Exponential bound when the number of components is large]
	\label{lem:large-count-exponential}
	Under~\eqref{eq:wasymptotic}, there exist constants $C, \kappa > 0$ such that, for all sufficiently large $n$ and every integer $k > 2m_n$,
	\[
		\Prb{S_k = n} \le \frac{C}{a_{m_n}} \exp\left(-\kappa\frac{k}{c_n}\right).
	\]
\end{lemma}

\begin{proof}
	Set $\eta := 1-e^{-1/2} > 1/3$. By concavity, $1-e^{-u/2} \ge \eta u$ for $0 \le u \le 1$. Because $X$ is nonnegative, it holds that
	\[
		\E{\exp\left(-\frac{X}{2a_{m_n}}\right)} \le 1-\eta\frac{\mu(a_{m_n})}{a_{m_n}} \le \exp\left(-\frac{\eta}{c_n}\right).
	\]
	Let $h = \lfloor 3k/4 \rfloor$. The exponential Markov inequality and $n \le b_{m_n} = m_n \mu(a_{m_n})$ give
	\[
		\Prb{S_h \le n} \le \exp\left(\frac{n}{2a_{m_n}}-\frac{\eta h}{c_n}\right) \le \exp\left(\frac{m_n/2-\eta h}{c_n}\right).
	\]
	For $k > 2m_n$, and all sufficiently large $n$,
	\[
		\frac{m_n}{2} - \eta h \le -\frac{3\eta-1}{4}k + \eta \le -\kappa k,
	\]
	where we can take $\kappa = (3\eta-1)/8$. Since $k-h \ge k/4 > m_n/2$, the estimate~\eqref{eq:common-concentration}, with reference index $m_n$ and $\theta = 1/2$, gives
	\[
		\sup_{z \in \ndZ} \Prb{S_{k-h} = z} \le C a_{m_n}^{-1}.
	\]
	By independence and nonnegativity of the summands, combining with the above bounds,
	\[
		\Prb{S_k = n} \le \Prb{S_h\le n}\sup_{z\in \ndZ}\Prb{S_{k-h} =z} \le C' a_{m_n}^{-1}\exp\left(-\kappa \frac{k}{c_n}\right).
	\]
\end{proof}
The following lemma estimates sums of $\Prb{S_h = q}$ over numbers of summands $h$ close to $m(q)$. In the infinite-mean case, the central-window calculation below appears in Berger's proof of the Cauchy renewal theorem~\cite[Section~8.1.2, Equations~(8.4)--(8.7)]{zbMATH07109856}. Here, we state the estimates for both finite and infinite mean.
\begin{lemma}[Sums of local probabilities]\label{lem:renewalwindow}
	Suppose that~\eqref{eq:wasymptotic} holds. For every fixed $A > 0$,
	\begin{equation}\label{eq:renewalcentral}
		\mu(a_{m(q)}) \sum_{h\,:\, |h-m(q)| \le Ac(q)}\Prb{S_h=q} \longrightarrow \int_{-A}^{A}f(-x)\diff x, \qquad q \rightarrow \infty, \qquad q \in \ndN.
	\end{equation}
	For any $0 < \lambda < 1 < \Lambda < \infty$ and $0 < \eta < 1$, there is $C < \infty$ such that, for all sufficiently large integers $q$ and all integers $h$ with $\lambda m(q)\le h\le\Lambda m(q)$,
	\begin{equation}\label{eq:common-window-envelope}
		\mu(a_{m(q)})\Prb{S_h=q} \le \frac{C}{c(q)} \left(1+\frac{|h-m(q)|}{c(q)}\right)^{-2+\eta}.
	\end{equation}
	It follows that, for all $A \ge 1$,
	\begin{equation}\label{eq:renewaltail}
		\limsup_{\substack{q\to\infty\\q\in \ndN}} \mu(a_{m(q)})\sum_{\substack{h\,:\,\lambda m(q)\le h\le\Lambda m(q)\\|h-m(q)|>Ac(q)}} \Prb{S_h=q} \le CA^{-1+\eta}.
	\end{equation}
	Note that the various constants $C$ may depend on $\lambda, \Lambda, \eta$, but not on $A$.
\end{lemma}

\begin{proof}
	In order to simplify notation, write $m = m(q)$ and $c = c(q)$.  Lemma~\ref{lem:seq}(f), the local limit theorem~\eqref{eq:llt}, and uniform continuity of $f$ argued earlier give, uniformly for $|h-m| \le Ac$,
	\[
		\mu(a_m)\Prb{S_h=q} = \frac{1}{c} \left(f\left(-\frac{h-m}{c}\right) + o(1)\right).
	\]
	A Riemann summation argument proves~\eqref{eq:renewalcentral}.

	By Lemma~\ref{lem:seq}(b),(d), for $\lambda m \le h \le \Lambda m$,
	\[
		r_{q,h} := |q-\lfloor b_h\rfloor| \ge |h-m| \mu(a_m)-o(a_m).
	\]
	Also $a_h\asymp a_m$ uniformly by~\eqref{eq:common-rv-window}. Choose a fixed $A_0$ large enough. For $|h-m| > A_0 c$ and all sufficiently large $q$, we then have $r_{q,h} \ge \frac{1}{2} |h-m| \mu(a_m) \ge a_h$ and $r_{q,h} \ge a_m$. Equation~\eqref{eq:bedev} and the Potter bounds~\ref{thm:pb} give, by using $mL(a_m) \sim a_m$,
	\begin{align*}
		\Prb{S_h=q} \le ChL(r_{q,h})r_{q,h}^{-2} \le C m L(a_m)a_m^{-\eta} \left(|h-m|\mu(a_m)\right)^{-2+\eta} \le Ca_m^{-1} \left(\frac{|h-m|}{c}\right)^{-2+\eta}.
	\end{align*}
	Since $a_m = c \mu(a_m)$, these bounds prove~\eqref{eq:common-window-envelope}. Summing this estimate over $\lambda m \le h \le \Lambda m$ with $|h-m| > Ac$, and using $\sum_{j > Ac}j^{-2+\eta} = O((Ac)^{-1+\eta})$, proves~\eqref{eq:renewaltail}.
\end{proof}

Next, we provide a bound when $k < m_n/2$. Monotonicity and regular variation of $b_k$, together with $b_{m_n} \sim n$, imply that, uniformly in this range and for $n$ large enough, $b_k \le 2n/3$. Hence $n- \lfloor b_k\rfloor \in [n/3,2n]$ and $n-\lfloor b_k \rfloor \ge a_k$. Equations~\eqref{eq:bedev} and~\eqref{eq:common-rv-window} give

\begin{equation}\label{eq:i2b}
	\Prb{S_k = n} \le Ck L(n)n^{-2} \le C'k \Prb{X = n}.
\end{equation}
This bound and Lemma~\ref{lem:renewalwindow} depend only on the inner weight asymptotic.

\subsection{Number of components}

For tilted inner weights with index strictly below $-2$, a dense-phase local limit theorem for the number of components is given in~\cite[Theorem~3.1(i)]{MR4780503}. The following theorem gives the Cauchy counterpart, with centering $m_n$ and scale $c_n$.

\begin{theorem}[Partition function and number of components]\label{thm:numco}
	Suppose that~\eqref{eq:wasymptotic},~\eqref{eq:critical} and~\eqref{eq:vasymptotic-dense} hold. Then
	\[
		u_n \sim \frac{L_v(m_n)}{\mu(a_{m_n})}  m_n^{-b} \rho_w^{-n}.
	\]
	For each $\delta > 0$, the number of components satisfies the uniform local limit theorem
	\[
		\sup_{k \ge \delta m_n}	\left| c_n \Prb{N_n = k} - f\left(-\frac{k - m_n}{c_n}\right) \right| \longrightarrow 0.
	\]
	Hence, the number of components satisfies
	\[
		\frac{N_n - m_n}{c_n} \convdis -Y.
	\]
\end{theorem}
\begin{proof}
	Fix $M \ge 1$ and $0 < \eta < 1$. By Proposition~\ref{pro:scale}, $c_n \rightarrow \infty$ and $c_n = o(m_n)$. Hence, uniformly for $|k-m_n| \le Mc_n$, it holds that $k/m_n \rightarrow 1$. Lemma~\ref{lem:seq}(f), applied with $q = n$, gives, uniformly in this range,
	\[
		\frac{a_k}{a_{m_n}} \longrightarrow 1, \qquad \frac{n-b_k}{a_k} = -\frac{k-m_n}{c_n}+o(1).
	\]
	Moreover,~\eqref{eq:vasymptotic-dense} and \eqref{eq:common-rv-shrinking}, applied to the outer weights, give $v_k W(\rho_w)^k \sim v_{m_n} W(\rho_w)^{m_n}$ uniformly on the same window. Combining these estimates with the local limit theorem~\eqref{eq:llt} and uniform continuity of $f$ yields
	\begin{equation}\label{eq:almostllt}
		v_k W(\rho_w)^k \Prb{S_k = n} = \frac{v_{m_n} W(\rho_w)^{m_n}}{c_n \mu(a_{m_n})} \left(f\left(-\frac{k-m_n}{c_n}\right) + o(1)\right),
	\end{equation}
	where we used $a_{m_n} =c_n \mu(a_{m_n})$ from~\eqref{eq:defcn}. The error is uniform over the window $|k-m_n| \le Mc_n$.

	Split the positive integers into four sets
	\[
		\begin{split}
			I_1&= \{k \ge 1:|k-m_n| \le Mc_n \},\\
			I_2&= \{k \ge 1:k<m_n/2\},\\
			I_3&= \{k \ge 1:k > 2m_n \},\\
			I_4&= \{k \ge 1:m_n/2 \le k \le 2m_n,\ |k-m_n|>Mc_n \}.
		\end{split}
	\]
	Since $c_n= o(m_n)$ by Proposition~\ref{pro:scale}, the above sets are disjoint and, for each fixed $M$ and $n$ large enough, they cover all positive integers.

	Summing the uniform error in~\eqref{eq:almostllt} over $I_1$ contributes $o(v_{m_n}W(\rho_w)^{m_n}/\mu(a_{m_n}))$. The spacing between consecutive points $(k-m_n)/c_n$ is $1/c_n$, which tends to zero by Proposition~\ref{pro:scale}. A Riemann sum argument gives
	\begin{equation} \label{eq:i1}
		\sum_{k \in I_1} v_k W(\rho_w)^k \Prb{S_k = n} \sim \frac{v_{m_n} W(\rho_w)^{m_n}}{\mu(a_{m_n})} \int_{-M}^{M}f(-x) \diff x.
	\end{equation}

	For $I_2$,~\eqref{eq:i2b} gives the upper bound
	\[
		\sum_{k \in I_2}v_k W(\rho_w)^k \Prb{S_k = n} \le CL(n)n^{-2} \sum_{k<m_n/2}k v_k W(\rho_w)^k.
	\]
	By~\eqref{eq:vasymptotic-dense}, $k v_k W(\rho_w)^k \sim L_v(k)k^{1-b}$. Since $1-b > -1$, the discrete version of Karamata's theorem~\ref{thm:karamata} yields
	\[
		\sum_{k < m_n/2}k v_k W(\rho_w)^k = O\left(L_v(m_n)m_n^{2-b}\right).
	\]
	Also,~\eqref{eq:defbn} and~\eqref{eq:bdiff2co} give $n \sim b_{m_n} = m_n \mu(a_{m_n})$. It follows that
	\begin{align}\label{eq:tu}
		\sum_{k\in I_2}v_k W(\rho_w)^k \Prb{S_k = n} = O\left(L(n)n^{-2}L_v(m_n)m_n^{2-b}\right) = O\left(L_v(m_n)m_n^{-b} \frac{L(n)}{\mu(a_{m_n})^2}\right).
	\end{align}
	Since $L(n) = o(\mu(a_{m_n}))$ by~\eqref{eq:la}, and $v_{m_n} W(\rho_w)^{m_n} \sim L_v(m_n)m_n^{-b}$
	by~\eqref{eq:vasymptotic-dense}, we obtain
	\begin{equation}\label{eq:i2}
		\sum_{k \in I_2} v_k W(\rho_w)^k \Prb{S_k = n} = o\left(\frac{v_{m_n} W(\rho_w)^{m_n}}{\mu(a_{m_n})}\right).
	\end{equation}

	For $I_3$, use the Potter bounds~\eqref{eq:pb} on the tilted outer weights, who are regularly varying with index $-b$ by~\eqref{eq:vasymptotic-dense}. Since $0 < \eta < 1 < b$, for $k > 2m_n,$ it holds that
	\[
		\frac{v_k W(\rho_w)^k}{v_{m_n} W(\rho_w)^{m_n}} \le C \left(\frac{k}{m_n}\right)^{-b+\eta}.
	\]
	In particular, this yields for the supremum
	\begin{equation}\label{eq:outer-weight-tail}
		\sup_{k > 2m_n} \frac{v_k W(\rho_w)^k}{v_{m_n} W(\rho_w)^{m_n}} \le C.
	\end{equation}
	Note that this argument also applies when $b=2$. Combining~\eqref{eq:outer-weight-tail} with Lemma~\ref{lem:large-count-exponential} gives
	\[
		\sum_{k \in I_3}v_k W(\rho_w)^k \Prb{S_k = n} \le \frac{Cv_{m_n} W(\rho_w)^{m_n}}{a_{m_n}} \sum_{k > 2m_n} e^{-\kappa k/c_n}.
	\]
	Since $c_n \rightarrow \infty$,
	\[
		\sum_{k \in I_3} e^{-\kappa k/c_n} = \frac{e^{-\kappa(2m_n+1)/c_n}}{1-e^{-\kappa/c_n}} = O\left(c_n e^{-2\kappa m_n/c_n}\right).
	\]
	Using $c_n/a_{m_n} = 1/\mu(a_{m_n})$ from~\eqref{eq:defcn}, and $m_n/c_n \rightarrow \infty$ from Proposition~\ref{pro:scale}, we conclude that
	\begin{align}\label{eq:i3}
		\sum_{k \in I_3} v_k W(\rho_w)^k \Prb{S_k = n} = O\left(\frac{v_{m_n} W(\rho_w)^{m_n}}{\mu(a_{m_n})} e^{-2\kappa m_n/c_n}\right) = o\left(\frac{v_{m_n} W(\rho_w)^{m_n}}{\mu(a_{m_n})}\right).
	\end{align}

	For $I_4$,~\eqref{eq:vasymptotic-dense} and \eqref{eq:common-rv-window}, applied on $[1/2,2]$, give $v_k W(\rho_w)^k \le Cv_{m_n} W(\rho_w)^{m_n}$ uniformly for $m_n/2 \le k \le 2m_n$. Hence
	\[
		\frac{\mu(a_{m_n})}{v_{m_n} W(\rho_w)^{m_n}} \sum_{k \in I_4} v_k W(\rho_w)^k \Prb{S_k = n} \le C\mu(a_{m_n})\sum_{k \in I_4}\Prb{S_k = n}.
	\]
	Use~\eqref{eq:renewaltail} with $q = n$, $\lambda=1/2$, $\Lambda=2$, and $A=M$. This gives
	\begin{equation}\label{eq:i4}
		\limsup_{n \rightarrow \infty} \frac{\mu(a_{m_n})}{v_{m_n} W(\rho_w)^{m_n}} \sum_{k\in I_4}v_k W(\rho_w)^k \Prb{S_k = n} \le CM^{-1+\eta},
	\end{equation}
	where $C$ is independent of $M$.

	By~\eqref{eq:intro-random-sum}, the partition function is the sum of the contributions from $I_1,\ldots,I_4$. The terms are nonnegative, so~\eqref{eq:i1}--\eqref{eq:i4} imply
	\begin{align*}
		\liminf_{n \rightarrow \infty} \frac{u_n \rho_w^n \mu(a_{m_n})}{v_{m_n} W(\rho_w)^{m_n}} &\ge \int_{-M}^{M}f(-x) \diff x,\\
		\limsup_{n \rightarrow \infty} \frac{u_n \rho_w^n \mu(a_{m_n})}{v_{m_n} W(\rho_w)^{m_n}} &\le \int_{-M}^{M}f(-x)\diff x+CM^{-1+\eta}.
	\end{align*}
	Letting $M \rightarrow \infty$ and using $\int_{\ndR} f(x)\diff x = 1$, followed by~\eqref{eq:vasymptotic-dense}, proves
	\begin{equation}\label{eq:upart}
		u_n \rho_w^n \sim\frac{v_{m_n} W(\rho_w)^{m_n}}{\mu(a_{m_n})} \sim\frac{L_v(m_n)m_n^{-b}}{\mu(a_{m_n})}.
	\end{equation}

	We now prove a local limit for the number of components. It is sufficient to consider the range $0 < \delta < 1/2$, since increasing the lower bound on $k$ can only decrease the supremum. Equations~\eqref{eq:components} and~\eqref{eq:upart}, together with $c_n\mu(a_{m_n}) = a_{m_n}$, give 
	\begin{equation}\label{eq:numcoupper}
		c_n \Prb{N_n = k} \le Ca_{m_n} \frac{v_k W(\rho_w)^k}{v_{m_n} W(\rho_w)^{m_n}} \Prb{S_k = n},
		\qquad k \ge 1.
	\end{equation}
    For every fixed $M > 0$, combining~\eqref{eq:almostllt},~\eqref{eq:components} and~\eqref{eq:upart} gives, uniformly over $|k-m_n| \le M c_n$,
	\begin{equation} \label{eq:lltcompact}
		\sup_{|k-m_n| \le Mc_n} \left|c_n \Prb{N_n = k} - f\left(-\frac{k-m_n}{c_n}\right)\right| \longrightarrow 0.
	\end{equation}

	On $\delta m_n \le k < m_n/2$,~\eqref{eq:vasymptotic-dense} and~\eqref{eq:common-rv-window}, used on $[\delta,1/2]$, bound the weight ratio in~\eqref{eq:numcoupper} by a constant depending only on $\delta$. Using~\eqref{eq:i2b}, we obtain, using $a_{m_n}\sim m_nL(a_{m_n})$ from~\eqref{eq:defan2} and $n \sim m_n \mu(a_{m_n})$ from~\eqref{eq:defbn} and~\eqref{eq:bdiff2co},
	\begin{align}\label{eq:ng1}
		\sup_{\delta m_n \le k < m_n/2} c_n \Prb{N_n = k} \le C_\delta a_{m_n} m_n L(n)n^{-2} = O\left(\frac{L(a_{m_n})L(n)}{\mu(a_{m_n})^2}\right) \longrightarrow 0.
	\end{align}
	The final convergence follows from~\eqref{eq:omun} and~\eqref{eq:la}.

	For $k > 2m_n$, combine~\eqref{eq:numcoupper}, \eqref{eq:outer-weight-tail}, and Lemma~\ref{lem:large-count-exponential}. The factors $a_{m_n}$ cancel, giving, for $k > 2m_n,$
	\begin{equation}\label{eq:count-exponential-tail}
		\Prb{N_n = k} \le \frac{C}{c_n} e^{-\kappa k/c_n}.
	\end{equation}
	Consequently, and because $m_n/c_n \rightarrow \infty$ by Proposition~\ref{pro:scale},
	\begin{equation}\label{eq:ng2}
		\sup_{k > 2m_n} c_n \Prb{N_n = k} \le Ce^{-2\kappa m_n/c_n} \longrightarrow 0.
	\end{equation}

	For $m_n/2 \le k \le 2m_n$, the weight ratio is uniformly bounded by~\eqref{eq:vasymptotic-dense} and~\eqref{eq:common-rv-window}. Using~\eqref{eq:common-window-envelope} with $q= n$, $\lambda=1/2$ and $\Lambda=2$. Together with \eqref{eq:numcoupper} and~\eqref{eq:defcn}, this gives
	\begin{align*}
		c_n \Prb{N_n = k} \le Ca_{m_n}\Prb{S_k = n} =Cc_n \mu(a_{m_n})\Prb{S_k = n} \le C\left(1+\frac{|k-m_n|}{c_n}\right)^{-2+\eta}.
	\end{align*}
	For $|k-m_n| > M c_n$, this is at most $C M^{-2+\eta}$. The ranges $\delta m_n \le k < m_n/2$ and $k > 2m_n$ contribute $o(1)$ to the supremum by~\eqref{eq:ng1} and~\eqref{eq:ng2}, while the density term is at most $\sup_{|x| > M} f(x)$. This results in
	\begin{equation}\label{eq:llttail}
		\limsup_{n \rightarrow \infty} \sup_{\substack{k \ge \delta m_n \\ |k-m_n| > Mc_n}} \left|c_n \Prb{N_n = k} - f\left(-\frac{k-m_n}{c_n}\right)\right| \le CM^{-2+\eta}+\sup_{|x| > M}f(x).
	\end{equation}
	The above constant is independent of $M$. Letting $M \rightarrow \infty$ and using~\eqref{eq:ftail}, combining it with \eqref{eq:lltcompact}, proves the uniform local limit.

	Finally, for fixed real $x < y$,~\eqref{eq:lltcompact} gives
	\[
		\Prb{x < \frac{N_n-m_n}{c_n} \le y} = \frac{1}{c_n} \sum_{\substack{k \ge 1 \\ x < (k-m_n)/c_n \le y}} f\left(-\frac{k-m_n}{c_n}\right) + o(1).
	\]
    For fixed $x < y$, the points $(k-m_n)/c_n$ with $m_n + xc_n < k \le m_n + yc_n$ lie in $]x,y]$ and have spacing $1/c_n \rightarrow 0$. Since $c_n = o(m_n)$, for $n$ large enough, every integer in this range is positive. The Riemann sum therefore converges to $\int_{x}^{y} f(-t) \diff t$. Since $f(-\cdot)$ is a probability density, this proves $(N_n-m_n)/c_n \convdis -Y$.\end{proof}

Even when $\E{X} = \mu <\infty$, the first-order equivalent $m_n \sim n/\mu$ cannot replace $m_n$ as the centering in the fluctuation limit. This is because Lemma~\ref{lem:seq}(b),~\eqref{eq:ellstarasymp} and~\eqref{eq:ellstarl} give together
\[
	\frac{m_n-n/\mu}{c_n} \sim \frac{\mu-\mu(a_{m_n})}{L(a_{m_n})} \longrightarrow \infty.
\]
This shows why it is crucial to keep $m_n$ in the local and distributional limits. 

\subsection{Component size asymptotics}

The next lemma approximates $(K_1, \ldots, K_{\min(N_n,m_n-s_n)})$ in total variation by $(X_1, \ldots, X_{m_n-s_n})$, where $c_n = o(s_n)$ and $s_n = o(m_n)$. It gives the Cauchy analogue of Stufler's total-variation approximation~\cite[Lemma~3.5]{MR4780503}.

\begin{lemma}\label{lem:tv}
    Suppose that~\eqref{eq:wasymptotic},~\eqref{eq:critical} and~\eqref{eq:vasymptotic-dense} hold.
	For any sequence $(s_n)_{n \ge 1}$ of positive integers satisfying $c_n = o(s_n)$ and $s_n = o(m_n)$, the following total-variation convergence holds:
	\begin{align*}
		d_{\mathrm{TV}} \left( (K_1, \ldots, K_{\min(N_n, m_n-s_n)}), (X_1, \ldots, X_{m_n-s_n}) \right) \longrightarrow 0.
	\end{align*}
	In particular, the assumptions hold (and the result is valid) for $s_n = \max(\lfloor \sqrt{c_n m_n} \rfloor, 1)$.
\end{lemma}
\begin{proof}
	Since $s_n = o(m_n)$, we can restrict ourselves to $n$ sufficiently large so that
	\[
		1 \le s_n < m_n.
	\]
	Set
	\[
		r_n := m_n-s_n, \qquad t_n := \left\lfloor\sqrt{c_ns_n}\right\rfloor.
	\]
	The assumptions on $s_n$ imply
	\begin{equation}\label{eq:orderseq}
		c_n= o(t_n), \qquad t_n= o(s_n), \qquad s_n= o(r_n), \qquad r_n \sim m_n.
	\end{equation}
	In particular, for all sufficiently large $n$,
	\[
		0 < c_n < t_n < s_n < r_n.
	\]

	Fixing $\varepsilon > 0$, and since $r_n \rightarrow \infty$, the stable convergence in~\eqref{eq:snbnclt} allows us to choose $M > 0$ such that
	\begin{align}\label{eq:sumclose}
		\Prb{|S_{r_n}-b_{r_n}| \le M a_{r_n}} >1-\varepsilon
	\end{align}
	for all sufficiently large $n$. Accordingly, define 
	\[
		B_n := \left\{\mathbf{x} \in \ndN_0^{r_n}: \Prb{(X_1,\ldots,X_{r_n}) = \mathbf{x}}>0,\quad \left|\sum_{i=1}^{r_n} x_i-b_{r_n}\right|\le Ma_{r_n}\right\}
	\]
	and
	\[
		E_n := \left\{ |N_n-m_n| \le t_n \right\}.
	\]

	On $E_n$,
	\[
		N_n \ge m_n-t_n > m_n-s_n = r_n.
	\]
	Therefore, for every $\mathbf{k} =(k_1,\ldots,k_{r_n})\in B_n$, Kolchin representation~\eqref{eq:KasX} gives
	\begin{equation*}
		\Prb{(K_1,\ldots,K_{r_n}) = \mathbf{k},\,E_n} = \sum_{j: |j-m_n| \le t_n} \Prb{(X_1, \ldots, X_{r_n}) = \mathbf{k}\mid S_j= n} \Prb{N_n= j}.
	\end{equation*}

	Independence of the variables $X_i$ gives
	\[
		\Prb{(X_1, \ldots, X_{r_n}) = \mathbf{k} \mid S_j = n} = \frac{\Prb{(X_1, \ldots, X_{r_n}) = \mathbf{k}} \Prb{S_{j-r_n} = n- \sum_{i=1}^{r_n} k_i}}{\Prb{S_j = n}}.
	\]

	Using the law of $N_n$ in~\eqref{eq:components}, this further yields
	\[
		\frac{\Prb{(K_1,\ldots,K_{r_n}) = \mathbf{k},\,E_n}}{\Prb{(X_1, \ldots, X_{r_n}) = \mathbf{k}}} = \frac{1}{u_n \rho_w^n} \sum_{j: |j-m_n| \le t_n} v_jW(\rho_w)^j \Prb{S_{j-r_n} = n-\sum_{i=1}^{r_n}k_i }.
	\]

	Using the asymptotics of $u_n$ in Theorem~\ref{thm:numco}, as well as~\eqref{eq:vasymptotic-dense} and slow variation of $L_v$, which uniformly gives, for $|j-m_n| \le t_n$, that
	\[
		v_jW(\rho_w)^j \sim L_v(m_n)m_n^{-b},
	\]
	it follows, uniformly for $\mathbf{k} \in B_n$, that
	\begin{equation}\label{eq:likrat}
		\frac{\Prb{(K_1, \ldots, K_{r_n}) = \mathbf{k},\,E_n}}{\Prb{(X_1, \ldots, X_{r_n}) = \mathbf{k}}} \sim \mu(a_{m_n}) \sum_{j: |j-m_n| \le t_n} \Prb{S_{j-r_n} = n-\sum_{i=1}^{r_n}k_i}.
	\end{equation}

     We now estimate the sum of probabilities in~\eqref{eq:likrat}, uniformly over $B_n$. For $\mathbf{k} \in B_n$, the remaining summands need to contribute
    \[
     \tilde{n} := n - \sum_{i=1}^{r_n} k_i.
    \]	

	Since $r_n \sim m_n$,~\eqref{eq:common-rv-shrinking} gives
	\begin{equation}\label{eq:arnamn}
		a_{r_n} \sim a_{m_n}.
	\end{equation}
	Lemma~\ref{lem:seq}(a), taking $q= n$ and $h=-s_n$ therein, also gives
	\begin{equation}\label{eq:rn}
		r_n| \mu(a_{m_n})-\mu(a_{r_n})| = O(s_nL(a_{m_n})+1) = o(a_{m_n}).
	\end{equation}
	Now, by~\eqref{eq:wasymptotic} and the definition of $X$, for $i$ large enough,
	\[
		i\Prb{X =i} = O\left(i^{-1} L(i)\right).
	\]
	The Potter bounds~\eqref{eq:pb} applied to $L$ give, for some $C_1 > 1$ and $n$ sufficiently large,
	\[
		\frac{L(i)}{L(a_{m_n})} \le C_1 \left(\frac{a_{m_n}}{i}\right)^{1/4}, \qquad a_{s_n}<i \le a_{m_n},
	\]
	and the same Potter bounds~\eqref{eq:pb} applied to the regularly varying function $(a_{\cdot})$ give, for some $C_2 > 1$ and $n$ large enough,
	\[
		\frac{a_{m_n}}{a_{s_n}} \le C_2 \left(\frac{m_n}{s_n}\right)^2.
	\]
	Combining all of the above bounds yields
	\begin{align}\label{eq:snass}
		s_n \left| \mu(a_{m_n})-\mu(a_{s_n}) \right| &= s_n \sum_{a_{s_n}<i \le a_{m_n}} i\Prb{X =i}\nonumber \\
		&\le C s_n L(a_{m_n}) a_{m_n}^{1/4} \sum_{i > a_{s_n}} i^{-5/4}\nonumber \\
		&\le C s_n L(a_{m_n}) \left(\frac{a_{m_n}}{a_{s_n}}\right)^{1/4}\nonumber \\
		&\le C L(a_{m_n}) \sqrt{m_n s_n} = o(m_n L(a_{m_n})) = o(a_{m_n}),
	\end{align}
	where the last equality follows from Equation~\eqref{eq:defan2}. Since we can write
	\[
		b_{m_n}-b_{r_n}-b_{s_n} = r_n \left(\mu(a_{m_n})-\mu(a_{r_n}) \right) + s_n (\mu(a_{m_n})-\mu(a_{s_n})),
	\]
	combining~\eqref{eq:rn} and~\eqref{eq:snass}, we obtain
	\[
		b_{m_n}-b_{r_n} = b_{s_n}+o(a_{m_n}).
	\]
	Equation~\eqref{eq:bdiff2} gives
	\[
		n-b_{m_n} = o(a_{m_n}).
	\]
	Hence, for $\mathbf{k} \in B_n$, and together with~\eqref{eq:arnamn}, this yields
	\[
		\sum_{i=1}^{r_n}k_i = b_{r_n} + O(a_{r_n}) = b_{r_n} + O(a_{m_n}).
	\]
	Collecting these estimates finally gives, uniformly for $\mathbf{k} \in B_n$,
	\begin{equation}\label{eq:ntildeorder}
		\tilde{n} = n-\sum_{i=1}^{r_n} k_i = (n-b_{m_n}) + (b_{m_n} - b_{r_n}) + \left( b_{r_n} - \sum_{i=1}^{r_n} k_i \right) = b_{s_n}+O(a_{m_n}).
	\end{equation}
	Hence, there exists a constant $D > 0$, which depends only on $M$, such that, for every $\mathbf{k} \in B_n$ and for all sufficiently large $n$,
	\begin{equation}\label{eq:nbsnbigO}
		|\tilde{n}-b_{s_n}| \le D a_{m_n}.
	\end{equation}

	We use Lemma~\ref{lem:renewalwindow} with $q = \tilde{n}$ therein. Writing $h := j-r_n$, the summation range in~\eqref{eq:likrat} becomes $|h-s_n| \le t_n$. We compare the center $m_{\tilde{n}}$ with $s_n$ and the scale $c_{\tilde{n}}$ with $t_n$. First, estimate~\eqref{eq:snass} and Equation~\eqref{eq:defcn} give, when combined with the order between sequences in~\eqref{eq:orderseq},
	\[
		\frac{|\mu(a_{m_n})-\mu(a_{s_n})|}{\mu(a_{m_n})} = O\left(\frac{a_{m_n}}{s_n \mu(a_{m_n})}\right) = O\left(\frac{c_n}{s_n}\right) = o(1).
	\]
	Hence,
	\[
		\mu(a_{s_n}) \sim \mu(a_{m_n}).
	\]
	Since $b_{s_n} = s_n \mu(a_{s_n})$ by~\eqref{eq:defbn},
	\[
		\frac{a_{m_n}}{b_{s_n}} = \frac{c_n}{s_n} \frac{\mu(a_{m_n})}{\mu(a_{s_n})} \longrightarrow 0.
	\]
	It follows from~\eqref{eq:ntildeorder} that, uniformly for $\mathbf{k} \in B_n$,
	\[
		\tilde{n} \sim b_{s_n}.
	\]
	We use the centering and scaling quantities from~\eqref{eq:defmn} and~\eqref{eq:defcn} in abbreviated forms for $\tilde{n},$ i.e.
	\[
		m_{\tilde{n}} := m(\tilde{n}) = \min \{k \ge 1:b_k \ge \tilde{n}\}, \qquad c_{\tilde{n}} := c(\tilde{n}) = \frac{a_{m_{\tilde{n}}}}{\mu(a_{m_{\tilde{n}}})}.
	\]
	Lemma~\ref{lem:seq}(e), together with a subsequence argument, gives, uniformly for $\mathbf{k} \in B_n$,
	\[
		m_{\tilde{n}} \sim s_n.
	\]
	Regular variation of $(a_{\cdot})$ and the uniform convergence theorem~\ref{thm:uct} give
	\[
		a_{m_{\tilde{n}}} \sim a_{s_n},
	\]
	while Potter bounds~\eqref{eq:pb} and $s_n = o(m_n)$ give
	\[
		a_{s_n} = o(a_{m_n}).
	\]
	By estimate~\eqref{eq:snass}, and using the slow variation of $\mu$, this yields, uniformly for $\mathbf{k} \in B_n$,
	\[
		\mu(a_{m_{\tilde{n}}}) \sim \mu(a_{s_n}) \sim \mu(a_{m_n}).
	\]

	We also need a sharper uniform bound, in the form of
	\begin{equation}\label{eq:mntildecnorder}
		|m_{\tilde{n}}-s_n| = O(c_n).
	\end{equation}
	To establish it, note that Lemma~\ref{lem:seq}(d), used with $q = \tilde{n}$ and $k = s_n$ therein, together with Lemma~\ref{lem:seq}(b) and~\eqref{eq:nbsnbigO}, gives
	\begin{align*}
		|m_{\tilde{n}}-s_n| \mu(a_{m_{\tilde{n}}}) \le |b_{m_{\tilde{n}}}-b_{s_n}| \le |b_{m_{\tilde{n}}}-\tilde{n}| + |\tilde{n}-b_{s_n}| = o(a_{m_{\tilde{n}}}) + O(a_{m_n}) = O(a_{m_n}).
	\end{align*}
    In the last equality, we have used the above estimates as well as~\eqref{eq:orderseq}, which gives $m_{\tilde{n}} \sim s_n = o(m_n)$. Dividing by $\mu(a_{m_{\tilde{n}}})$ proves the claim. Furthermore, uniformly for $\mathbf{k} \in B_n$,
	\[
		\frac{c_{\tilde{n}}}{c_n} = \frac{a_{m_{\tilde{n}}}}{a_{m_n}} \frac{\mu(a_{m_n})}{\mu(a_{m_{\tilde{n}}})} \longrightarrow 0.
	\]
	Hence,
	\begin{equation}\label{eq:cntildeorder}
		c_{\tilde{n}} = o(c_n) = o(t_n).
	\end{equation}

	Fix $A > 0$. Equations~\eqref{eq:mntildecnorder} and~\eqref{eq:cntildeorder} imply, uniformly on $B_n$, that, eventually,
	\[
		\{|h-m_{\tilde{n}}| \le Ac_{\tilde{n}}\} \subseteq \{|h-s_n| \le t_n\}.
	\]
	Since $m_{\tilde{n}} \sim s_n$ and $t_n = o(s_n)$, the interval $|h-s_n| \le t_n$ lies in $[m_{\tilde{n}}/2, 2m_{\tilde{n}}]$. Also, note that, uniformly on $B_n$, $\tilde{n} \sim b_{s_n} \rightarrow \infty$. Apply~\eqref{eq:renewalcentral} on $|h-m_{\tilde{n}}| \le A c_{\tilde{n}}$ and~\eqref{eq:renewaltail} to the rest of this interval. Using $\mu(a_{m_{\tilde{n}}}) \sim \mu(a_{m_n})$, we obtain, for $0 < \eta < 1$,
	\[
		\limsup_{n \rightarrow \infty} \sup_{\mathbf{x} \in B_n} \left|\mu(a_{m_n})\sum_{h: |h-s_n|\le t_n}\Prb{S_h= \tilde{n}}-\int_{-A}^{A}f(-x)\diff x\right| \le CA^{-1+\eta}.
	\]
	Letting $A \rightarrow \infty$ yields, uniformly on $B_n,$
	\begin{equation}\label{eq:renewaluniform}
		\sum_{h: |h-s_n|\le t_n} \Prb{S_h= \tilde{n}} \sim \frac{1}{\mu(a_{m_n})}.
	\end{equation}
	With $h = j-r_n$, the range in~\eqref{eq:likrat} becomes $|h-s_n| \le t_n$, and all such $h$ are positive. Hence,
	\begin{equation}\label{eq:likrat2}
		\delta_n := \sup_{\mathbf{x} \in B_n} \left|\frac{\Prb{(K_1, \ldots, K_{r_n}) = \mathbf{x},E_n}}{\Prb{(X_1, \ldots, X_{r_n}) = \mathbf{x}}}-1\right| \longrightarrow 0.
	\end{equation}

    We have proved that the ratio of probability masses in~\eqref{eq:likrat2} tends uniformly to one on $B_n$. Now, we bound the probability outside $B_n$ and remove the restriction to $E_n$. Theorem~\ref{thm:numco}, together with $t_n/c_n \rightarrow \infty$ from~\eqref{eq:orderseq}, gives	
    \[
		\Prb{E_n^c} = \Prb{\left| \frac{N_n-m_n}{c_n} \right| > \frac{t_n}{c_n}} \longrightarrow 0.
	\]
	Furthermore, by the choice of $M$ in~\eqref{eq:sumclose}, and for $n$ sufficiently large,
	\begin{equation}\label{eq:Xvecbn}
		\Prb{(X_1, \ldots, X_{r_n}) \notin B_n} < \varepsilon.
	\end{equation}
	Equation~\eqref{eq:likrat2} implies that
	\[
		\delta_n := \sup_{\substack{\mathbf{k} \in B_n \\ \Prb{(X_1, \ldots, X_{r_n}) = \mathbf{k}} > 0}} \left\vert \frac{\Prb{(K_1, \ldots, K_{r_n}) = \mathbf{k}, E_n}}{\Prb{(X_1, \ldots, X_{r_n}) = \mathbf{k}}}  - 1 \right\vert \longrightarrow 0.
	\]
	Since $\Prb{E_n} \rightarrow 1$, it follows that
	\begin{equation}\label{eq:Kvecbn}
		\Prb{(K_1, \ldots, K_{r_n}) \notin B_n,\, E_n} \le \varepsilon + o(1).
	\end{equation}
	Collecting everything, let $\cS$ denote the space of finite sequences, i.e.
	\[
		\cS := \bigsqcup_{k > 0} \N_0^k.
	\]
	Then, for any measurable set $\cA \subseteq \cS$, since $(K_1, \ldots, K_{\min(N_n, r_n)})$ and $(X_1, \ldots, X_{r_n})$ are random elements of $\cS,$ and because on $E_n$ it holds that $N_n \ge r_n$, and hence that
	\[
		(K_1, \ldots, K_{\min(N_n,r_n)}) = (K_1,\ldots,K_{r_n}),
	\]
	this gives, splitting the first probability of interest,
	\begin{align*}
		\lvert \Prb{(K_1, \ldots, K_{\min(N_n,r_n)}) \in \cA} & - \Prb{(X_1, \ldots, X_{r_n}) \in \cA} \rvert \le \Prb{E_n^c} \\
		&\hspace{-0.7cm}\quad + \lvert \Prb{(K_1,\ldots,K_{r_n}) \in \cA, E_n} - \Prb{(X_1, \ldots, X_{r_n}) \in \cA} \rvert.
	\end{align*}
	Split the second term further, according to $B_n$, and let $\cA_n := \cA \cap \N_0^{r_n}.$ Then,
	\begin{align*}
		\lvert \Prb{(K_1,\ldots,K_{r_n}) \in \cA, E_n} & - \Prb{(X_1, \ldots, X_{r_n}) \in \cA} \rvert \\
		&\hspace{-2cm} \le \lvert \Prb{(K_1,\ldots,K_{r_n}) \in \cA_n \cap B_n, E_n} - \Prb{(X_1, \ldots, X_{r_n}) \in \cA_n \cap B_n} \rvert \\
		&\hspace{-2cm}\quad +\Prb{(K_1, \ldots, K_{r_n}) \notin B_n, E_n} + \Prb{(X_1, \ldots, X_{r_n}) \notin B_n}.
	\end{align*}
	Using now Equations~\eqref{eq:Xvecbn} and~\eqref{eq:Kvecbn} yields
	\[
		d_{\mathrm{TV}} \left((K_1, \ldots, K_{\min(N_n,r_n)}),(X_1, \ldots, X_{r_n})\right) \le \Prb{E_n^c} + \delta_n + 2\varepsilon + o(1).
	\]
	Taking the supremum over $\cA,$ letting $n \rightarrow \infty$ and finally $\varepsilon \rightarrow 0$, shows that
	\[
		d_{\mathrm{TV}} \left((K_1, \ldots, K_{\min(N_n,r_n)}),(X_1, \ldots, X_{r_n})\right) \longrightarrow 0.
	\]
	This proves the main claim of the lemma. It remains to verify the validity of the proposed choice for $s_n$
	\[
		s_n = \max \left( \left\lfloor\sqrt{c_n m_n}\right\rfloor, 1 \right).
	\]
	By Proposition~\ref{pro:scale}, $c_n \rightarrow \infty$ and $c_n = o(m_n)$. Hence,
	\[
		s_n \sim \sqrt{c_nm_n},
	\]
	and therefore, by~\eqref{eq:orderseq},
	\[
		\frac{c_n}{s_n} \sim \sqrt{\frac{c_n}{m_n}} \longrightarrow 0, \qquad \frac{s_n}{m_n} \sim \sqrt{\frac{c_n}{m_n}} \longrightarrow 0.
	\]
	This shows that the choice of $(s_n)$ satisfies all the assumptions, completing the proof.
\end{proof}

We use Lemma~\ref{lem:tv} to extend the dense-phase functional limit of~\cite[Theorem~3.3]{MR4780503} to the Cauchy boundary. Let $(Y_s)_{s \ge 0}$ be the spectrally positive one-stable Lévy process. Its value at time $t = 1$ has the law of $Y$ in~\eqref{eq:chara}. Its characteristic function is given by
\[
	\E{ e^{itY_s} } = \exp\left(s \int_0^\infty \frac{1}{x^2} \left(e^{itx}-1-itx \one_{\{x \le 1\}} \right)\,\diff x \right), \quad t > 0.
\]

\begin{theorem}[Functional limit theorem]\label{thm:functional} 
Suppose that~\eqref{eq:wasymptotic},~\eqref{eq:critical} and~\eqref{eq:vasymptotic-dense} hold. Then
	\[
		\left(\frac{\sum_{i=1}^{\lfloor tN_n \rfloor}K_i - t \mu(a_{m_n}) N_n }{a_{m_n}},\,\,0 \le t \le 1 \right) \convd (Y_t,\, 0 \le t \le 1)
	\]
	in \(\mathbb{D}([0,1], \ndR)\) endowed with the Skorokhod $J_1$-topology.
\end{theorem}

\begin{proof}
	First, note that, for all $x > 0,$ it holds, by regular variation of $X$ and the choice of $(a_{\cdot})$, that
	\begin{equation*}
		m_n \Prb{X > a_{m_n} x} = m_n \Prb{X > a_{m_n}} \frac{\Prb{X > a_{m_n}x}}{\Prb{X > a_{m_n}}} \longrightarrow x^{-1}.
	\end{equation*}
	Equation~\eqref{eq:common-truncated-power}, taking $p = 2$ and $k = m_n$ therein, gives
	\[
		\lim_{\varepsilon \rightarrow 0} \limsup_{n \rightarrow \infty} \frac{m_n}{a_{m_n}^2}\E{X^2; X \le \varepsilon a_{m_n}} =0.
	\]
    The terms with $X_i \le \varepsilon a_{m_n}$ are the jumps of size at most $\varepsilon$ in the rescaled partial-sum process. The above second-moment estimate verifies the so-called small-jump condition~\cite[Equation~(7.6)]{resnick07}. Hence~\cite[Theorem~7.1]{resnick07} gives, with respect to the $J_1$ topology,	    
    \[
		\left(\frac{\sum_{i=1}^{\lfloor t m_n \rfloor} X_i - t m_n \mu(a_{m_n})}{a_{m_n}}, \, 0 \le t \le 1\right) \convd (Y_t,\, 0 \le t \le 1).
	\]
	Let $r_n := m_n-s_n \sim m_n$, with $s_n$ as in Lemma~\ref{lem:tv}, and define the time change $\gamma_n(t) := tr_n/m_n$. Composing the above independent-sample process with $\gamma_n$ replaces $\lfloor tm_n \rfloor$ by $\lfloor tr_n \rfloor$ and the centering $tm_n \mu(a_{m_n})$ by $tr_n \mu(a_{m_n})$. The maps $\gamma_n$ converge uniformly to the identity map, and are both continuous and strictly increasing. Whitt's composition theorem~\cite[Theorem~3.1 and the following remark]{whitt1980functions} applies, since $Y$ has no fixed-time jump at $1$. Lemma~\ref{lem:tv} and~\eqref{eq:tv-kernel-contract}, applied to the partial-sum map, give the same limit for the first $r_n$ component sizes. On the event $\{N_n < r_n\}$, which tends to zero in probability, extend missing entries by zero if necessary. Therefore, with respect to the $J_1$ topology,
	\begin{equation}\label{eq:kconv}
		\left(\frac{\sum_{i=1}^{\lfloor tr_n \rfloor} K_i - t r_n \mu(a_{m_n})}{a_{m_n}}, \, 0 \le t \le 1\right) \convd (Y_t,\, 0 \le t \le 1).
	\end{equation}
    Fix $0 < \delta < 1$. We next replace $r_n$ by $N_n$ using the (random) time change $\lambda_n(t) := tN_n/r_n$. Since $N_n/r_n \convp 1$,
    \[
        \sup_{0\le t \le \delta} |\lambda_n(t)-t| = \delta\left|\frac{N_n}{r_n}-1\right| \convp 0.
    \]
    Moreover, $\Prb{\delta N_n \le r_n} \rightarrow 1$, so $\lambda_n$ sends $[0,\delta]$ into $[0,1]$ with probability tending to one. Joint convergence with this asymptotically deterministic map follows from Slutsky's theorem. Use again \cite[Theorem~3.1 and the following remark]{whitt1980functions}, now using that $Y$ has no fixed-time jump at $\delta$. Composing the process in~\eqref{eq:kconv} with $\lambda_n$ changes its centering exactly to $tN_n \mu(a_{m_n})$. Hence
	\begin{equation}\label{eq:Znconvdist}
		\left(Z_n(t), \, 0 \le t \le \delta \right) := \left(\frac{\sum_{i=1}^{\lfloor t N_n \rfloor} K_i - t N_n \mu(a_{m_n})}{a_{m_n}},\, 0 \le t \le \delta \right) \convd (Y_t,\, 0 \le t \le \delta).
	\end{equation}
	We have proved convergence on every $[0, \delta]$ with $\delta < 1$. It remains to control the terminal interval near $1$. We use a time-reversal argument as in the proof of Stufler's theorem~\cite[Theorem~3.3]{MR4780503}. Define
	\[
		Z_n^{\leftarrow}(t) := \frac{ \sum_{i=N_n-\lfloor tN_n \rfloor+1}^{N_n}K_i-tN_n \mu(a_{m_n})}{a_{m_n}}, \qquad 0 \le t \le 1.
	\]

	Conditionally on $N_n$, the elements of $(K_1,\ldots,K_{N_n})$ are exchangeable. Consequently,
	\[
		Z_n^{\leftarrow} \eqdist Z_n.
	\]
	It follows from the convergence in~\eqref{eq:Znconvdist} that, for every fixed $\eta \in ]0,1[$, with respect to the $J_1$-topology,
	\[
		\left(Z_n^{\leftarrow}(t),\,0 \le t \le \eta\right) \convd \left(Y_t,\,0 \le t \le \eta\right)
	\]
	The mapping $x \mapsto \sup_{0 \le t \le \eta} |x(t)|$ is continuous on $\mathbb{D}([0,\eta],\ndR)$ endowed with the $J_1$-topology. Hence,
	\[
		\sup_{0 \le t \le \eta} |Z_n^{\leftarrow}(t)| \convd \sup_{0 \le t \le \eta}|Y_t|.
	\]
	By the usual Portmanteau theorem, and because $(\mathbb{D}([0, \eta], \R), J_1)$ is a metric space (see~\cite[Chapter~12]{whitt2002}), it holds that, for every $\epsilon > 0$,
	\[
		\limsup_{n \rightarrow \infty}\Prb{ \sup_{0 \le t \le \eta}|Z_n^{\leftarrow}(t)|>\epsilon} \le \Prb{\sup_{0 \le t \le \eta}|Y_t| \ge\epsilon }.
	\]
	Since $Y_0 = 0$ and $Y$ is càdlàg, its paths are right-continuous at $0$, and therefore
	\[
		\sup_{0 \le t \le \eta}|Y_t| \longrightarrow 0\,\, \text{almost surely, as }\eta \rightarrow 0.
	\]
	It follows that
	\[
		\lim_{\eta \rightarrow 0}\limsup_{n \rightarrow \infty}\Prb{\sup_{0 \le t \le \eta}|Z_n^{\leftarrow}(t)|>\epsilon} = 0.
	\]

	For $s\in [1-\eta,1]$, we can write
	\[
		Z_n(1)-Z_n(s) =Z_n^{\leftarrow} \left( \frac{N_n-\lfloor sN_n \rfloor}{N_n} \right) + \frac{ \left( N_n-\lfloor sN_n \rfloor-(1-s)N_n \right)\mu(a_{m_n})}{a_{m_n}}.
	\]
	Since
	\[
		N_n-\lfloor sN_n \rfloor-(1-s)N_n = sN_n-\lfloor sN_n \rfloor \in [0,1[,
	\]
	and
	\[
		0 \le \frac{N_n-\lfloor sN_n \rfloor}{N_n} \le \eta+\frac{1}{N_n},
	\]
	we obtain
	\[
		\sup_{1-\eta \le s\le 1}|Z_n(1)-Z_n(s)| \le \sup_{0 \le t \le \eta+1/N_n}|Z_n^{\leftarrow}(t)| + \frac{\mu(a_{m_n})}{a_{m_n}}.
	\]
	Moreover, $\frac{\mu(a_{m_n})}{a_{m_n}} \rightarrow 0 $ as argued in~\eqref{eq:muan}. Then, for every fixed $\eta > 0$,
	\[
		\Prb{\sup_{0 \le t \le \eta+1/N_n}|Z_n^{\leftarrow}(t)|>\epsilon} \le \Prb{N_n<1/\eta} + \Prb{\sup_{0 \le t \le 2\eta}|Z_n^{\leftarrow}(t)| > \epsilon}.
	\]
	Consequently, for every $\epsilon > 0$,
	\[
		\lim_{\eta \rightarrow 0} \limsup_{n \rightarrow \infty} \Prb{\sup_{1-\eta \le s \le 1}|Z_n(s)-Z_n(1)|>\epsilon} = 0.
	\]

	For $0 < \eta < 1$, we freeze the path by setting $Z_n^{[\eta]}(t) = Z_n(t\wedge(1-\eta))$ and $Y^{[\eta]}_t = Y_{t\wedge(1-\eta)}$. Equation~\eqref{eq:Znconvdist} gives $Z_n^{[\eta]} \convdis Y^{[\eta]}$ with respect to the $J_1$ topology. Moreover,
    \[
        \|Z_n-Z_n^{[\eta]}\|_\infty \le 2\sup_{1-\eta \le s \le 1}|Z_n(s)-Z_n(1)|.
    \]
    The above estimate makes this difference negligible as first $n \rightarrow \infty$ and then $\eta \rightarrow 0$. Since $Y$ is almost surely continuous at $1$, $Y^{[\eta]}$ converges uniformly to $Y$ almost surely. The converging-together theorem~\cite[Theorem~3.5]{resnick07} finally yields
	\[
		\left(Z_n(t),\,0 \le t \le 1\right) \convd \left(Y_t,\,0 \le t \le 1\right)
	\]
	in $\mathbb{D}([0,1],\ndR)$ with respect to the $J_1$-topology.
\end{proof}

The positive jumps of the partial sum process are the rescaled component sizes. As in~\cite[Corollary~3.4]{MR4780503}, the ranked sizes in the limit are the ordered jumps of the limiting process. The endpoint of that process also determines the limiting $Y$ of the component count.
\begin{corollary}[Largest components]\label{cor:largest}
    Suppose that~\eqref{eq:wasymptotic},~\eqref{eq:critical} and~\eqref{eq:vasymptotic-dense} hold. Define the point process of the largest scaled sizes by
    \[
        \cN_n := \sum_{\substack{1\le i \le N_n \\ K_i > 0}} \delta_{K_i/a_{m_n}}.
    \]
    Then, 
    \begin{equation}\label{eq:dense-count-extremes-joint}
        \left(\frac{N_n-m_n}{c_n}, \cN_n\right) \convdis (-Y,\cN)
    \end{equation}
     in $\ndR \times \cM_p(]0,\infty])$, with respect to the vague topology on the second coordinate. Here, we note that $Y = Y_1$ and $\cN = \sum_{j \ge 1} \delta_{J_j}$ records the positive jump sizes of the same limiting Lévy process on $[0,1]$. The points $J_j$ are the atoms of $\cN$, listed in decreasing order. The point process is Poisson with intensity measure given by $\nu(\diff x) = x^{-2}\,\diff x$.

    Let $\Gamma_j = E_1 + \cdots + E_j$, where $(E_i)_{i \ge 1}$ are independent exponential random variables with unit mean. For every fixed $r \ge 1$, and jointly with~\eqref{eq:dense-count-extremes-joint}
    \[
        \left(\frac{K_{(1)}}{a_{m_n}},\ldots,\frac{K_{(r)}}{a_{m_n}}\right) \convdis (J_1,\ldots,J_r) \eqdist(\Gamma_1^{-1},\ldots,\Gamma_r^{-1}).
    \]

    In particular, for every $j \ge 1$ and $x > 0$,
    \[
        \Prb{\frac{K_{(j)}}{a_{m_n}} \le x} \longrightarrow e^{-1/x}\sum_{k=0}^{j-1}\frac{x^{-k}}{k!}.
    \]
\end{corollary}

\begin{proof}
    Let $Z_n$ be the process in Theorem~\ref{thm:functional}, and define
    \[
        \Phi(z) := \sum_{\substack{0 < t \le 1 \\ \Delta z(t) > 0}} \delta_{\Delta z(t)}.
    \]
    Since the centering of $Z_n$ is continuous, its positive jumps are the nonzero terms in $(K_i/a_{m_n})_{1 \le i \le N_n}$, so $\cN_n \equiv \Phi(Z_n)$.
    
    To apply the continuous mapping theorem, we need to check continuity of the map $\Phi$. If $z_m \rightarrow z$ in $J_1$, and $z$ is continuous at the endpoints, then, for every $\varepsilon > 0$ outside the jump sizes of $z$, the positive jumps larger than $\varepsilon$ of $z_m$ are eventually in one-to-one correspondence with those of $z$, with converging times and sizes (we refer to e.g.~\cite[Propositions~VI.2.1 and~VI.2.7]{JacodShiryaev2003}). Every continuous, compactly supported test function on $]0,\infty]$ tends to 0 near zero, so its integrals against $\Phi(z_m)$ and $\Phi(z)$ are finite sums of corresponding terms. Hence $\Phi(z_m) \rightarrow \Phi(z)$ vaguely.
    
    The limiting process $Y_{\cdot}$ is almost surely continuous at the endpoints, and its positive jump sizes form a Poisson point process $\cN := \Phi(Y_{\cdot})$ with intensity measure given by $\nu(\diff x) = x^{-2}\diff x$, see~\cite[Section~5.5]{resnick07}. Since evaluation at $t = 1$ is also continuous with respect to the $J_1$ topology, Theorem~\ref{thm:functional} and the continuous mapping theorem~\cite[Theorem~3.1]{resnick07} give
    \[
        (Z_n(1), \cN_n) \convdis (Y_1, \cN).
    \]
    By Lemma~\ref{lem:seq}(b),
    \[
        \frac{N_n-m_n}{c_n} + Z_n(1) = \frac{n-b_{m_n}}{a_{m_n}} \longrightarrow 0.
    \]
    Slutsky's theorem therefore yields
    \[
        \left(\frac{N_n-m_n}{c_n}, \cN_n\right) \convdis (-Y_1, \cN),
    \]
    which proves~\eqref{eq:dense-count-extremes-joint}. Lemma~\ref{lem:meso-poisson-cauchy} gives the compensated-sum representation of $Y_1$ from the same Poisson points.
    
    Since the intensity measure $\nu$ is diffuse, $\cN$ is almost surely simple by~\cite[Proposition~6.9]{lastPenrose2017}. Furthermore,
    \[
        \nu(]\varepsilon,\infty]) = \varepsilon^{-1}, \qquad \nu(]0,\infty]) = \infty, \qquad \nu(\{\infty\}) = 0.
    \]
    The Poisson count distributions~\cite[Definition~3.1]{lastPenrose2017} imply that $\cN$ has infinitely many positive points, finitely many above each positive threshold, and none at $\infty$. Write its points as $J_1 > J_2 > \cdots$.
    
    The ranking map is continuous at any point measure $\xi$ which satisfy the above properties. This can be seen by writing its points as $x_1 > x_2 > \cdots$, fixing $r \ge 1$, and choosing $\varepsilon \in (x_{r+1}, x_r)$. If $\xi_m \rightarrow \xi$ vaguely, then~\cite[Proposition~3.12]{resnick1987} implies that $\xi_m$ has (eventually) exactly $r$ points above $\varepsilon$, and exactly one in each of arbitrarily small disjoint intervals around $x_1, \ldots, x_r$. Hence, its first $r$ ranked points converge to $(x_1, \ldots, x_r)$. Applying the continuous mapping theorem~\cite[Theorem~3.1]{resnick07} once more gives, jointly with~\eqref{eq:dense-count-extremes-joint},
    \[
        \left(\frac{K_{(1)}}{a_{m_n}}, \ldots, \frac{K_{(r)}}{a_{m_n}}\right) \convdis (J_1, \ldots, J_r).
    \]
    
    By the Poisson mapping theorem~\cite[Theorem~5.1]{lastPenrose2017}, $x \mapsto 1/x$ changes $\cN$ into a unit-rate Poisson point process. Hence, we have the distributional identity $(J_j)_{j \ge 1} \eqdist (\Gamma_j^{-1})_{j \ge 1}$. Finally, for $x > 0$,
    \[
        \Prb{J_j \le x} = \Prb{\cN(]x,\infty]) \le j-1} = e^{-1/x} \sum_{k=0}^{j-1} \frac{x^{-k}}{k!}.
    \]
\end{proof}

\begin{remark}\label{rem:atom-matching}
    The continuity argument above also applies jointly with other convergent coordinates. We use it below for sums over points above a fixed positive threshold and for ranked component sizes.
\end{remark}

\begin{example}
	\label{ex:running}
	Suppose that the assumptions of Theorem~\ref{thm:numco} hold and that, as in Example~\ref{ex:addme}, for constants $c>0$ and $r \in \ndR$,
	\[
		L(n)\sim c \log^r n.
	\]
	The asymptotic estimate for the partition function in Theorem~\ref{thm:numco} becomes
	\[
		u_n \sim \begin{cases} \displaystyle \left(\frac{c}{r+1}\right)^{b-1} L_v\left( \frac{r+1}{c}\frac{n}{\log^{r+1}n} \right) n^{-b}\log^{(b-1)(r+1)}(n)\,\rho_w^{-n}, & r > -1,\\[2.2ex]
			\displaystyle c^{b-1} L_v\left( \frac{n}{c \log \log n} \right) n^{-b}(\log \log n)^{b-1}\rho_w^{-n}, & r = -1,\\[2.2ex]
			\displaystyle \mu^{b-1} L_v\left(\frac{n}{\mu}\right) n^{-b}\rho_w^{-n}, & r < -1. \end{cases}
	\]
	For the total variation approximation in Lemma~\ref{lem:tv}, we may take
	\[
		s_n= \max \left( \left\lfloor\sqrt{c_nm_n}\right\rfloor,1 \right).
	\]
	In the three regimes this sequence satisfies
	\[
		s_n \sim \begin{cases} \displaystyle \frac{(r+1)^{3/2}}{c}\frac{n}{\log^{r+3/2}n}, & r > -1,\\[1.4ex]
			\displaystyle \frac{n}{c\sqrt{\log n}(\log \log n)^{3/2}}, & r = -1,\\[1.4ex]
			\displaystyle \frac{\sqrt{c}}{\mu^{3/2}} n \log^{r/2}n, & r < -1. \end{cases}
	\]
	The largest components scale at order
	\begin{align*}
		a_{m_n} \sim \begin{cases} \displaystyle (r+1)\frac{n}{\log n}, & r > -1,\\[1.2ex]
			\displaystyle \frac{n}{\log n \log \log n}, & r = -1,\\[1.2ex]
			\displaystyle \frac{c}{\mu}n \log^r n, & r < -1. \end{cases}
	\end{align*}
\end{example}
\begin{proof}
	The asymptotic relations for $m_n$ and $c_n$ were established in Example~\ref{ex:addme}. In all three regimes, it holds that
	\begin{align}\label{eq:logmnlogn}
		\log m_n \sim\log n.
	\end{align}

	We only show the above identities in the case $r > -1$, but the others are obtained straightforwardly using the asymptotics from the last sections. When $r > -1$, and using Example~\ref{ex:addme} as well as Equations~\eqref{eq:muanrgm1} and~\eqref{eq:logmnlogn} yields
	\[
		\frac{m_n^{-b}}{\mu(a_{m_n})} \sim \frac{\left(\frac{r+1}{c} \cdot \frac{n}{\log^{r+1}n}\right)^{-b}}{\frac{c}{r+1}\log^{r+1}m_n} \sim \left(\frac{c}{r+1}\right)^{b-1} n^{-b} \log^{(b-1)(r+1)}(n),
	\]
	which, when combining with Theorem~\ref{thm:numco}, gives
	\[
		u_n \sim \left(\frac{c}{r+1}\right)^{b-1} L_v\left(\frac{r+1}{c} \cdot \frac{n}{\log^{r+1}n} \right) n^{-b} \log^{(b-1)(r+1)}(n) \rho_w^{-n},
	\]
	where the slow variation of $L_v$ allows to substitute $m_n$ as its argument by the equivalent relationship given in Example~\ref{ex:addme}, establishing the first identity.

	Lemma~\ref{lem:tv} allows $s_n = \max(\lfloor \sqrt{c_n m_n}\rfloor, 1)$ and, as $n \rightarrow \infty,$ since $c_n m_n \rightarrow \infty,$
	\[
		s_n \sim \sqrt{c_n m_n}.
	\]
	Using the identities for $c_n, m_n$ in Example~\ref{ex:addme} with~\eqref{eq:logmnlogn} gives precisely
	\[
		\sqrt{c_n m_n} \sim \sqrt{\frac{(r+1)^2}{c} \frac{n}{\log^{r+2}n} \frac{r+1}{c} \frac{n}{\log^{r+1}n}} \sim \frac{(r+1)^{3/2}}{c} \frac{n}{\log^{r+3/2}n}.
	\]
	Using the asymptotic relation for $a_{m_n}$ in Example~\ref{ex:addme}, as well as the relation for $m_n$ and~\eqref{eq:logmnlogn} directly gives
	\[
		a_{m_n} \sim c m_n \log^{r}m_n \sim (r+1) \frac{n}{\log n}.
	\]
	All remaining claims follow in the same fashion.
\end{proof}

\section{The condensation regime}
\label{sec:condensation}

\subsection{Convergent case}

Recall the definition of the unnormalised truncated first moment, for $x \ge 1,$
\begin{equation}\label{eq:defHboundary}
	H(x) := \sum_{1 \le k \le x} k v_k W(\rho_w)^k,
\end{equation}
and write $H(\infty) := \lim_{x \rightarrow \infty}H(x)\in ]0, \infty]$. This way, $H(\infty) < \infty$ if and only if $V'(W(\rho_w)) < \infty$.

If $V(W(\rho_w)) < \infty$, we can define a random variable $N$, independent of $(X_i)_{i \ge 1}$, and with distribution given by
\begin{equation}\label{eq:Nconv}
	\Prb{N = k} = \frac{v_k W(\rho_w)^k}{V(W(\rho_w))}.
\end{equation}
The above definition gives $N_n \eqdist (N\mid S_N = n)$. When $H(\infty) < \infty$, the variable $N$ has finite mean. Let $\widehat{N}$ be independent of $(X_i)_{i \ge 1}$ and have the size-biased distribution of $N$, where each probability $\Prb{N = k}$ is multiplied by $k$, and is then normalised. Rigorously, 
\begin{align*}
 \Prb{\widehat{N} = k} = \frac{k\Prb{N = k}}{\E{N}} = \frac{k v_k W(\rho_w)^k}{H(\infty)} = \frac{k v_k W(\rho_w)^{k-1}}{V'(W(\rho_w))},\qquad k \ge 1.
\end{align*}
The factor $k$ corresponds to the possible positions of the (distinguished) large component.

Conditionally on $\widehat{N}$, we let $J$ be uniformly distributed on $[\widehat{N}]$ and independent of $(X_i)_{i \ge 1}$. For each $n \ge 1$ we form the random vector
\begin{align}\label{eq:Xi}
	\boldsymbol{\Xi}_n := (X_1,\ldots,X_{J-1},n-S_{\widehat{N}-1},X_J,\ldots,X_{\widehat{N}-1}).
\end{align}
If $S_{\widehat{N}-1} > n$, we instead set $\boldsymbol{\Xi}_n = (n,0,\ldots,0)$ with exactly $\widehat{N}$ coordinates.

The total-variation approximation below is of the same form as in Stufler's convergent-condensation approximation in~\cite[Theorem~3.11]{MR4780503}. We give here a proof under the stated assumptions, using the local random-sum criteria of Bloznelis~\cite{MR4038060} for the coefficient asymptotics.

\begin{proposition}[Convergent condensation]
	\label{thm:condensation}
	Suppose that~\eqref{eq:wasymptotic} holds and $H(\infty) < \infty$. Each of the following conditions is sufficient for
	\[
		u_n \sim V'(W(\rho_w)) w_n
	\]
	and
	\[
		d_{\mathrm{TV}}\left((K_1,\ldots,K_{N_n}),\boldsymbol{\Xi}_n \right) \longrightarrow 0.
	\]
	In particular,
	\[
		N_n \convdis \widehat{N}, \qquad n-K_{(1)} \convdis S_{\widehat{N}-1}.
	\]
	\begin{itemize}
		\item[1.] Subcritical composition scheme:
			\begin{equation}\label{cond:subcritical}
				W(\rho_w) < \rho_v.
			\end{equation}
		\item[2.] Critical scheme with outer exponent $b > 2$:
			\begin{equation}\label{cond:critWdom}
				W(\rho_w) = \rho_v, \qquad v_n \sim L_v(n)n^{-b}W(\rho_w)^{-n}, \quad b > 2.
			\end{equation}
		\item[3.] For some $\delta > 0$,
			\begin{equation}\label{cond:genmomN}
				\E{N^{1+\delta}} < \infty.
			\end{equation}
		\item[4.] If $L_w(n) \rightarrow c\in ]0,\infty[$, it is enough that,
			for some $\delta > 0$,
			\begin{equation}\label{cond:logmomN}
				\E{N(\log N)^{2+\delta}} < \infty.
			\end{equation}
	\end{itemize}
\end{proposition}
\begin{proof}
	By \cite[Theorem~3]{MR4038060}, with $\alpha=2$, condition~\eqref{cond:genmomN} gives
	\begin{equation}\label{eq:sbj}
		\Prb{S_N= n} \sim \E{N}\Prb{X = n}.
	\end{equation}
	Under condition~\eqref{cond:logmomN}, the same conclusion follows from \cite[Theorem~2(ii)]{MR4038060}. Condition~\eqref{cond:subcritical} implies an exponential moment, since
	\[
		\E{s^N} = \frac{V(sW(\rho_w))}{V(W(\rho_w))}<\infty, \text{ for some } s > 1.
	\]
	Under~\eqref{cond:critWdom}, the Potter bounds~\eqref{eq:pb} give $\E{N^{1+\delta}}<\infty$ for every $0 < \delta < b-2$. Then~\eqref{eq:sbj} holds under each stated condition. By~\eqref{eq:intro-random-sum} and~\eqref{eq:Nconv},
	\[
		\Prb{S_N = n} = \frac{u_n \rho_w^n}{V(W(\rho_w))}, \qquad \E{N} = \frac{H(\infty)}{V(W(\rho_w))}.
	\]
	Consequently,
	\[
		u_n \rho_w^n \sim H(\infty)\Prb{X = n},
	\]
	which is the claimed coefficient asymptotic.

	We now prove the total-variation approximation. Let $I_n$ be the position of the first maximal component. The triples below take values in the countable space
	\[
		\bigsqcup_{k \ge 1} \left(\{k\} \times \ndN_0^{k-1} \times [k]\right).
	\]

	Fix $k \ge 1$, $\mathbf{x} =(x_1,\ldots,x_{k-1})$ and $j \in [k]$. For sufficiently large $n$,
	$n-\sum_{i}x_i > \max(0,x_1,\ldots,x_{k-1})$. Equations~\eqref{eq:components} and~\eqref{eq:KasX} then give
	\begin{align*}
		\Prb{N_n = k, \tau(K_1,\ldots,K_{N_n}) = \mathbf{x},I_n= j} &= \frac{v_k W(\rho_w)^k \Prb{X = n-\sum_{i}x_i}}{u_n \rho_w^n}\prod_{i=1}^{k-1}\Prb{X =x_i}\\
		&\quad\longrightarrow \frac{\Prb{\widehat{N} =k}}{k} \prod_{i=1}^{k-1}\Prb{X =x_i}\\
		&\quad= \Prb{\widehat{N} =k,\,(X_1,\ldots,X_{k-1}) = \mathbf{x}, J= j}.
	\end{align*}
	The limit uses the coefficient asymptotic just proved and $\Prb{X = n-\sum_{i}x_i}/\Prb{X = n} \rightarrow 1$, which follows from Theorem~\ref{thm:uct}. Since the limiting masses form a probability distribution, Scheffé's lemma gives the total-variation approximation
	\[
		d_{\mathrm{TV}}\left((N_n,\tau(K_1,\ldots,K_{N_n}),I_n),(\widehat{N},(X_1,\ldots,X_{\widehat{N}-1}),J)\right) \longrightarrow 0.
	\]

	For each $n$, apply inequality~\eqref{eq:tv-kernel-contract} to the measurable maps
	\[
		\begin{aligned}
			(k,\mathbf{x},j) &\longmapsto (x_1,\ldots,x_{j-1},n-\sum_{i}x_i,x_j,\ldots,x_{k-1}),\\
			(k,\mathbf{x},j) &\longmapsto (k,\sum_{i}x_i).
		\end{aligned}
	\]
	On $\{\sum_{i} x_i > n\}$, extend the first map using the convention in~\eqref{eq:Xi}. The first map sends the two triples to $(K_1,\ldots,K_{N_n})$, respectively $\boldsymbol{\Xi}_n$, proving the claimed full-vector approximation. The second gives
	\[
		d_{\mathrm{TV}}\left((N_n,n-K_{(1)}), (\widehat{N},S_{\widehat{N}-1})\right) \longrightarrow 0,
	\]
	which proves the stated limits for the number of components and the residual mass.
\end{proof}

\subsection{Boundary condensation}

Throughout this subsection, we continue to assume~\eqref{eq:wasymptotic} and, further, that the tilted outer weight sequences has index $-2,$ i.e.
\begin{align}\label{eq:vboundary}
	v_k \sim L_v(k)k^{-2} W(\rho_w)^{-k}.
\end{align}
Recall $H$ from~\eqref{eq:defHboundary}. By~\eqref{eq:vboundary},
\begin{equation}\label{eq:kvkWasymp}
	k v_k W(\rho_w)^k \sim k^{-1} L_v(k).
\end{equation}
Hence,
\[
	H(\infty) = \infty \iff \sum_{k \ge 1} k^{-1}L_v(k) = \infty.
\]

Suppose then that $H(\infty) = \infty$. We assume $L_v$ to be locally bounded on $[1,\infty[$ (we may modify it on any bounded interval, if necessary). The uniform convergence theorem~\ref{thm:uct} gives
\[
    \int_{k}^{k+1}t^{-1} L_v(t) \diff t \sim k^{-1} L_v(k), \quad k \rightarrow \infty.
\]
Since $H(\infty) = \infty$, summing the above comparison and using~\eqref{eq:kvkWasymp} gives
\[
    H(x) \sim \int_{1}^{x} t^{-1}L_v(t)\,\diff t.
\]
Proposition~1.5.9a of~\cite{bgt1989} shows that $H$ is slowly varying and $L_v(x) = o(H(x))$. If $H(\infty) < \infty$, slow variation follows directly from convergence to a finite, positive limit. Hence, in both cases,
\begin{equation}\label{eq:Hratio}
	\frac{H(tx)}{H(x)} \longrightarrow 1, \qquad t > 0.
\end{equation}

For all sufficiently large $n$, let $\widehat{N}_n$ be independent of $(X_i)_{i \ge 1}$, and with distribution
\begin{align}\label{eq:defNhatboundary}
	\Prb{\widehat{N}_n=k} = \frac{k v_k W(\rho_w)^k}{H(m_n)}, \qquad 1 \le k \le m_n.
\end{align}
Conditionally on $\widehat{N}_n$, let $J_n$ be uniformly distributed on $[\widehat{N}_n]$ and independent from $(X_i)_{i \ge 1}$. We set
\begin{align}\label{eq:defXigel}
	\boldsymbol{\Xi}^{\mathrm{g}}_n := (X_1, \ldots, X_{J_n-1}, n-S_{\widehat{N}_n-1}, X_{J_n},\ldots, X_{\widehat{N}_n-1}).
\end{align}
If $S_{\widehat{N}_n-1} > n$, set $\boldsymbol{\Xi}_n^{\mathrm{g}} = (n,0,\ldots,0)$ with exactly $\widehat{N}_n$ coordinates, using the convention in~\eqref{eq:Xi}.

The factor $k$ in~\eqref{eq:defNhatboundary} counts the possible positions for the large component. If $H(\infty) < \infty$, then $\widehat{N}_n$ converges in distribution to $\widehat{N}$. If $H(\infty) = \infty$, then, for every fixed positive integer $M$, and for $n$ large enough,
\[
 \Prb{\widehat{N}_n \le M} = \frac{H(M)}{H(m_n)} \rightarrow 0.
\]
Hence, $\widehat{N}_n$ tends to infinity in probability. In either case, for every fixed $0 < \varepsilon < 1$, slow variation~\eqref{eq:Hratio} gives
\[
 \Prb{\widehat{N}_n > \varepsilon m_n} = 1 - \frac{H(\varepsilon m_n)}{H(m_n)} \rightarrow 0.
\]
It follows that $\widehat{N}_n/m_n \convp 0$.

To obtain condensation without conditioning on $\cG_n$, we assume that the contribution from numbers of components $k = o(m_n)$ dominates the contribution from $k = m_n+O(c_n)$. The precise assumption, sometimes referred to as the single-big-jump dominance, is
\begin{align}\label{eq:geldominance}
	\frac{L_v(m_n)\mu(a_{m_n})}{H(m_n)L(n)} \longrightarrow 0.
\end{align}

For $k \ge 1$, let the sequence of order statistics of $X_1, \ldots, X_k$ be denoted by
\[
	X_{(1):k} \ge X_{(2):k} \ge \cdots \ge X_{(k):k}.
\]
The sequence is empty when $k = 0$, and we set $S_0 := 0$. With the deletion map $\tau$ defined in~\eqref{eq:deletion-map}, we can now state the following uniform single-big-jump approximation. Equation~\eqref{eq:uniformbigjumplocal} is a specialisation of the arguments of Denisov, Dieker and Shneer~\cite{MR2440928}, and its consequence~\eqref{eq:uniformbigjumptv} is due to Armendáriz and Loulakis~\cite{MR2775110}. We provide a brief justification for the reader's convenience.

\begin{lemma}[Uniform single-big-jump approximation]
	\label{lem:uniformbigjump}
	Suppose that~\eqref{eq:wasymptotic} holds, and let $(t_n)_{n \ge 1}$ be a sequence of positive integers satisfying $t_n = o(m_n)$. Then
	\begin{align}\label{eq:uniformbigjumplocal}
		\sup_{1 \le k \le t_n} \left|\frac{\Prb{S_k = n}}{k \Prb{X = n}}-1\right|\longrightarrow 0
	\end{align}
	and
	\begin{align}\label{eq:uniformbigjumptv}
		\sup_{1 \le k \le t_n}d_{\mathrm{TV}} \left(\tau\big((X_1, \ldots, X_k) \mid S_k = n \big),(X_1, \ldots, X_{k-1})\right) \longrightarrow 0.
	\end{align}
\end{lemma}

\begin{proof}
	Since $\mu$ is nondecreasing, $t_n = o(m_n)$ and $b_{m_n} \sim n$ by~\eqref{eq:bdiff2co},
	\[
		\frac{b_{t_n}}{n} \sim\frac{b_{t_n}}{b_{m_n}} = \frac{t_n \mu(a_{t_n})}{m_n \mu(a_{m_n})} \le \frac{t_n}{m_n} \longrightarrow 0.
	\]
	Moreover, $a_{t_n} \le a_{m_n} = o(n)$ by~\eqref{eq:amo}. Hence, uniformly for $1 \le k \le t_n$, and because $\lfloor b_k \rfloor \le b_{t_n},$
	\[
		\frac{n-\lfloor b_k\rfloor}{a_k} = \frac{n-\lfloor b_k \rfloor}{n} \frac{n}{a_k} = \left(1-\frac{\lfloor b_k \rfloor}{n}\right) \frac{n}{a_k} \longrightarrow \infty, \quad \text{ and } n-\lfloor b_k\rfloor\sim n.
	\]
	Equation~\eqref{eq:beloc} yields, uniformly along every sequence $k=k_n \le t_n$ with $k_n \rightarrow \infty$,
	\[
		\Prb{S_k = n} = \Prb{S_k - \lfloor b_k \rfloor = n-\lfloor b_k \rfloor} \sim kL(n)n^{-2} \sim k \Prb{X = n},
	\]
	where the first equivalence follows from Theorem~\ref{thm:uct} and the second from~\eqref{eq:lawX} and~\eqref{eq:wasymptotic}. For each fixed $k$, the same estimate follows from~\cite[Theorem~3]{MR4038060}, applied with the deterministic stopping index $N \equiv k$ and $\alpha = 2$. It is uniform when $k$ ranges over a fixed, finite set of positive integers. A subsequence argument gives~\eqref{eq:uniformbigjumplocal} uniformly over the full claimed range.

	If $t_n \rightarrow \infty$, the stable convergence in~\eqref{eq:snbnclt} gives
	\[
		\frac{S_{t_n}}{n} = \frac{b_{t_n}}{n}+ \frac{a_{t_n}}{n} \frac{S_{t_n}-b_{t_n}}{a_{t_n}} \convp 0,
	\]
	by using the above asymptotic relations and because $(S_{t_n}-b_{t_n})/a_{t_n}$ is tight. If $(t_n)$ is bounded, the same conclusion is immediate, and the general case follows simply by a subsequence argument. This means that there exists a deterministic sequence $(l_n) \rightarrow 0$ (sufficiently slowly) such that
	\begin{equation}\label{eq:smallremainderuniform}
		\Prb{S_{t_n} > l_n n} \longrightarrow 0.
	\end{equation}
	Define
	\begin{equation}\label{eq:ank}
		A_{n,k} := \left\{\mathbf{x} \in \N_0^{k-1}: \sum_{i=1}^{k-1} x_i \le l_n n, \quad \prod_{i=1}^{k-1} \Prb{X = x_i} > 0\right\}.
	\end{equation}
	For $\mathbf{x} = (x_1, \ldots, x_{k-1}) \in A_{n,k}$, write $s = \sum_{i=1}^{k-1}x_i$. Since $l_n \rightarrow 0$, eventually
	\[
		n-s \ge (1-l_n)n > l_n n \ge \max(0, x_1, \ldots, x_{k-1}),
	\]
	and, consequently, $n-s$ gives the unique maximal coordinate, which occupies any of the $k$ positions, yielding, on the event $A_{n,k},$
	\[
		\frac{\Prb{\tau(X_1, \ldots, X_k) = \mathbf{x} \mid S_k = n}}{\Prb{(X_1, \ldots, X_{k-1}) = \mathbf{x}}} = \frac{k \Prb{X = n-s}}{\Prb{S_k = n}} = \frac{\Prb{X = n-s}}{\Prb{X = n}} \frac{k \Prb{X = n}}{\Prb{S_k = n}}.
	\]
	By~\eqref{eq:uniformbigjumplocal} and the uniform convergence theorem~\ref{thm:uct},
	\begin{equation}\label{eq:deltanrat}
		\delta_n := \sup_{k \le t_n}\sup_{\mathbf{x} \in A_{n,k}} \left| \frac{\Prb{\tau(X_1, \ldots, X_k) = \mathbf{x} \mid S_k = n}}{\Prb{(X_1, \ldots, X_{k-1}) = \mathbf{x}}}-1 \right|\longrightarrow 0.
	\end{equation}
	Moreover, and uniformly for $k \le t_n$, by~\eqref{eq:smallremainderuniform}, it holds that
	\[
		\Prb{(X_1, \ldots, X_{k-1}) \in A_{n,k}^{c}} = \Prb{S_{k-1} > l_n n} \le \Prb{S_{t_n} > l_n n} \longrightarrow 0.
	\]
	Equation~\eqref{eq:deltanrat} also yields
	\[
		\Prb{\tau(X_1, \ldots, X_k) \in A_{n,k}^c \mid S_k = n} \le \Prb{(X_1, \ldots, X_{k-1}) \in A_{n,k}^c} + \delta_n.
	\]
	Therefore, and uniformly for $k \le t_n,$
	\[
		d_{\mathrm{TV}} \left(\left(\tau(X_1, \ldots, X_k) \mid S_k = n \right), (X_1, \ldots, X_{k-1})\right) \le \Prb{(X_1, \ldots, X_{k-1}) \in A_{n,k}^c} + \delta_n,
	\]
	which proves~\eqref{eq:uniformbigjumptv}.
\end{proof}

The next theorem combines Lemma~\ref{lem:uniformbigjump} with an estimate of the distribution of $N_n$. It applies irrespective of the value of $H(\infty)$. 

\begin{theorem}[Boundary single-big-jump approximation]
	\label{thm:boundarycondensation}
	Suppose that~\eqref{eq:wasymptotic},~\eqref{eq:critical},~\eqref{eq:vboundary} and~\eqref{eq:geldominance} hold. Then
	\begin{align}\label{eq:boundarypartition}
		u_n \sim \frac{H(m_n)}{W(\rho_w)} w_n.
	\end{align}
	In addition, it holds that
	\begin{align}\label{eq:boundaryvectortv}
		d_{\mathrm{TV}} \left((K_1,\ldots,K_{N_n}),\boldsymbol{\Xi}^{\mathrm{g}}_n \right) \rightarrow 0.
	\end{align}
	In particular,
	\begin{align}\label{eq:boundarynumbertv}
		d_{\mathrm{TV}}(N_n,\widehat{N}_n) \longrightarrow 0.
	\end{align}
\end{theorem}

\begin{proof}
	Using the slowly varying behaviour of $H$ established in~\eqref{eq:Hratio}, there exists a sequence of positive integers $(t_n)$ satisfying
	\begin{align}\label{eq:boundaryrn}
		t_n \rightarrow \infty, \qquad t_n= o(m_n), \qquad\frac{H(t_n)}{H(m_n)} \rightarrow 1.
	\end{align}
	If $H(\infty) < \infty$, take $t_n = \lfloor \sqrt{m_n} \rfloor$. Suppose now that $H(\infty) = \infty$. For each fixed integer $j \ge 2,$ choose increasing thresholds $n_j$ such that, for all $n \ge n_j, $
	\[
		\frac{H(\lfloor m_n/j\rfloor)}{H(m_n)} \ge1-\frac{1}{j}, \qquad m_n \ge j^2.
	\]
	Taking $t_n = \lfloor m_n / j \rfloor$ whenever $n_j \le n < n_{j+1}$ yields, for $t_n \ge j$ and $t_n/m_n \le 1/j$,
	\[
		1-\frac{1}{j} \le \frac{H(t_n)}{H(m_n)} \le 1.
	\]
	Since $j \rightarrow \infty$ as $n \rightarrow \infty$, this proves all three assertions in~\eqref{eq:boundaryrn}. Now, we split the sum over numbers of components into four distinct ranges:
	\[
		I_1 := \{k \le t_n \}, \, I_2 := \{t_n < k < m_n/2 \}, \, I_3 := \{m_n/2 \le k \le 2m_n \}, \, I_4 := \{k > 2m_n \}.
	\]
	For $I_1,$ Equation~\eqref{eq:uniformbigjumplocal} in Lemma~\ref{lem:uniformbigjump}, together with~\eqref{eq:lawX} and the definition of $H$ give
	\begin{align}\label{eq:boundarysmallrange}
		\sum_{k \in I_1}  v_k W(\rho_w)^k \Prb{S_k = n} \sim \Prb{X = n} H(t_n)\sim \frac{w_n \rho_w^n}{W(\rho_w)}H(m_n).
	\end{align}

	For $I_2$, note that Equation~\eqref{eq:i2b} gives, for some constants $C, C' >0,$ uniformly for $k < m_n/2$, and in combination with~\eqref{eq:lawX} and~\eqref{eq:wasymptotic},
	\[
		\Prb{S_k = n} \le CkL(n)n^{-2} \le C'k \Prb{X = n}.
	\]
	Therefore,~\eqref{eq:Hratio} and~\eqref{eq:boundaryrn} yield
	\begin{align}\label{eq:boundaryintermediate}
		\sum_{k \in I_2} v_k W(\rho_w)^k \Prb{S_k = n}&\le C' \Prb{X = n} \left(H(m_n/2)-H(t_n)\right)\notag\\
		&= o\left(\frac{w_n \rho_w^n}{W(\rho_w)}H(m_n)\right).
	\end{align}

	For $I_3$, use Lemma~\ref{lem:renewalwindow}, and in particular~\eqref{eq:renewalcentral} and~\eqref{eq:renewaltail} therein (with $q \equiv n$), gives
	\[
		\sum_{k \in I_3}\Prb{S_k = n} = O\left(\mu(a_{m_n})^{-1}\right).
	\]
	Moreover,~\eqref{eq:vboundary} together with the uniform convergence theorem~\ref{thm:uct} give, uniformly for $m_n / 2 \le k \le 2m_n,$ $v_k W(\rho_w)^k = O(L_v(m_n) m_n^{-2})$. Combining with the above,
	\begin{align}\label{eq:boundarycentral}
		\sum_{k \in I_3} v_k W(\rho_w)^k \Prb{S_k = n} = O\left(L_v(m_n)m_n^{-2}\mu(a_{m_n})^{-1}\right).
	\end{align}

	Finally, for $I_4$, Lemma~\ref{lem:large-count-exponential} and the Potter bounds argument in~\eqref{eq:outer-weight-tail}, now with $b = 2$ therein, give
	\begin{equation}\label{eq:boundarylarge}
		\sum_{k \in I_4}v_k W(\rho_w)^k \Prb{S_k = n} = O\left(L_v(m_n)m_n^{-2}\mu(a_{m_n})^{-1} \exp\left(-2\kappa \frac{m_n}{c_n}\right)\right).
	\end{equation}
	By~\eqref{eq:lawX}, \eqref{eq:wasymptotic} and~\eqref{eq:bdiff2co}, for $n$ large enough,
	\[
		\Prb{X = n}H(m_n) \sim\frac{H(m_n)L(n)}{m_n^2\mu(a_{m_n})^2}.
	\]
	By the single-big-jump dominance assumption~\eqref{eq:geldominance},
	\[
		\frac{L_v(m_n)m_n^{-2}/\mu(a_{m_n})}{\Prb{X = n}H(m_n)} \sim \frac{L_v(m_n)\mu(a_{m_n})}{H(m_n)L(n)} \longrightarrow 0.
	\]
	Hence~\eqref{eq:boundarycentral} and~\eqref{eq:boundarylarge} are negligible when compared to $\Prb{X = n}H(m_n)$.
	Together with~\eqref{eq:boundarysmallrange} and \eqref{eq:boundaryintermediate}, this controls all four ranges
	and gives
	\begin{align}\label{eq:boundaryZ}
		u_n \rho_w^n \sim \Prb{X = n}H(m_n) \sim \frac{w_n \rho_w^n}{W(\rho_w)} H(m_n) \implies u_n \sim \frac{H(m_n)}{W(\rho_w)} w_n,
	\end{align}
    which is precisely~\eqref{eq:boundarypartition}. 
    
    We next compare the distributions of $N_n$ and $\widehat{N}_n$. For $n$ large enough, $t_n \le m_n$ and $\Prb{S_k = n} > 0$ for every $1 \le k \le t_n$, by~\eqref{eq:boundaryrn} and~\eqref{eq:uniformbigjumplocal}. This gives, for this range,
    \[
        \Prb{N_n = k} > 0 \iff v_k > 0 \iff \Prb{\widehat{N}_n = k} > 0.
    \]
    Equations~\eqref{eq:components},~\eqref{eq:defNhatboundary},~\eqref{eq:uniformbigjumplocal} and~\eqref{eq:boundaryZ} give, uniformly for $1 \le k \le t_n$ with $v_k > 0$,
    \begin{equation}\label{eq:boundarymixingratio}
        \frac{\Prb{N_n = k}}{\Prb{\widehat{N}_n = k}} = \frac{\Prb{S_k = n}}{k \Prb{X = n}} \frac{\Prb{X = n}H(m_n)}{u_n \rho_w^n} \longrightarrow 1.
    \end{equation}
    For $1 \le k \le t_n$ with $v_k=0$, both probabilities vanish. Further, by~\eqref{eq:defNhatboundary} and~\eqref{eq:boundaryrn} ($H$ being slowly varying),
	\[
		\Prb{\widehat{N}_n \le t_n} = \frac{H(t_n)}{H(m_n)} \longrightarrow 1,
	\]
	whereas~\eqref{eq:components},~\eqref{eq:boundarysmallrange} and~\eqref{eq:boundaryZ} yield
	\[
		\Prb{N_n \le t_n} = \frac{\sum_{k \le t_n} v_k W(\rho_w)^k \Prb{S_k = n}}{u_n \rho_w^n} \sim\frac{H(t_n)}{H(m_n)} \longrightarrow 1\longrightarrow 1,
	\]
	where the last convergence follows from~\eqref{eq:boundaryrn}. Hence~\eqref{eq:boundarymixingratio} proves~\eqref{eq:boundarynumbertv}.

	Let $I_n$ be the position of the first maximal component. Keep $A_{n,k}$ from~\eqref{eq:ank} and the sequences $l_n$ and $\delta_n$ from the proof of Lemma~\ref{lem:uniformbigjump}. Fix $1 \le k \le t_n$ with $v_k > 0$. For $\mathbf{x} \in A_{n,k}$, the coordinate deleted by $\tau$ from~\eqref{eq:deletion-map} is the unique maximum for all sufficiently large $n$. The enumeration of its $k$ possible insertion positions used in that proof gives
	\[
		\Prb{\tau(K_1, \ldots, K_{N_n}) = \mathbf{x},\ I_n= j \mid N_n = k} = \frac{1}{k} \Prb{\tau(X_1,\ldots,X_k) = \mathbf{x}\mid S_k = n}, \qquad j \in [k].
	\]
	Let $U_k$ be uniformly distributed on $[k]$, independently of $(X_i)_{i \ge 1}$. By the above equality and~\eqref{eq:deltanrat}, for every $\mathbf{x} \in A_{n,k}$ and $j \in [k]$,
	\begin{multline*}
		\left| \Prb{\tau(K_1,\ldots,K_{N_n}) = \mathbf{x},\ I_n= j \mid N_n = k} - \frac{1}{k} \Prb{(X_1,\ldots,X_{k-1}) = \mathbf{x}} \right|  \\ \le \frac{\delta_n}{k} \Prb{(X_1,\ldots,X_{k-1}) = \mathbf{x}}.
	\end{multline*}
	The second term on the left-hand side is precisely that of $((X_1, \ldots, X_{k-1}),U_k)$ at $(\mathbf{x}, j)$. Appealing to the total-variation argument from the proof of Lemma~\ref{lem:uniformbigjump} yields, uniformly for $1 \le k \le t_n$ with $v_k > 0$,
	\begin{equation}\label{eq:boundarymarkedfixedk}
		d_{\mathrm{TV}}\Big(\big((\tau(K_1,\ldots,K_{N_n}),I_n)\mid N_n = k\big),\big((X_1,\ldots,X_{k-1}),U_k\big)\Big)\le \delta_n+\Prb{S_{t_n}>l_n n} =: \varepsilon_n \longrightarrow 0.
	\end{equation}
	Summing these bounds with weights $\Prb{N_n = k}$ gives
	\begin{align}\label{eq:boundarymarkedtv}
		d_{\mathrm{TV}}\Big((N_n,\,&\tau(K_1,\ldots,K_{N_n}),I_n), \left(\widehat{N}_n, (X_1,\ldots,X_{\widehat{N}_n-1}),J_n \right) \Big)\\
		&\le \sum_{\substack{k \ge 1\\ \Prb{N_n = k}>0}}\Prb{N_n = k}\, d_{\mathrm{TV}}\Big(\left((\tau(K_1,\ldots,K_k),I_n)\mid N_n = k\right), \left((X_1,\ldots,X_{k-1}),U_k\right) \Big)\nonumber\\
		&\quad +\frac{1}{2}\sum_{k \ge 1} \left|\Prb{N_n = k}-\Prb{\widehat{N}_n=k}\right|\nonumber\\
		&\le \varepsilon_n \Prb{N_n \le t_n} +\Prb{N_n > t_n} +d_{\mathrm{TV}}(N_n,\widehat{N}_n) \longrightarrow 0.\nonumber
	\end{align}
	The convergence follows from~\eqref{eq:boundarynumbertv} and $\Prb{N_n \le t_n} \rightarrow 1$. All remaining assertions follow from \eqref{eq:tv-kernel-contract}, applied to the following measurable maps: inserting $n$ minus the sum of the remaining coordinates at the recorded position, with the convention in~\eqref{eq:defXigel}, reconstructs $(K_1, \ldots, K_{N_n})$, respectively $\boldsymbol{\Xi}^{\mathrm{g}}_n$, proving~\eqref{eq:boundaryvectortv}.
\end{proof}

When $H(\infty) < \infty$, the auxiliary number of components $\widehat{N}_n$ converges in total variation to $\widehat{N}$, defined before Proposition~\ref{thm:condensation}. The boundary theorem recovers convergent condensation, where both $N_n$ and the residual mass $n-K_{(1)}$ remain bounded in probability.

\begin{corollary}[Convergent boundary condensation]
	\label{pro:boundarycondensation}
	Suppose that the assumptions of Theorem~\ref{thm:boundarycondensation} hold and that $H(\infty) < \infty$. Let $\widehat{N}$ have the size-biased law introduced before Proposition~\ref{thm:condensation}, independently of $(X_i)_{i \ge 1}$. Then,
	\[
		u_n \sim V'(W(\rho_w)) w_n
	\]
	and
	\begin{align}\label{eq:boundaryfinitejoint}
		d_{\mathrm{TV}} \left((K_1, \ldots, K_{N_n}), \boldsymbol{\Xi}_n \right) \rightarrow 0.
	\end{align}
\end{corollary}

\begin{proof}
	By~\eqref{eq:defHboundary}, $H(m_n) \rightarrow H(\infty) = W(\rho_w)V'(W(\rho_w))$. The coefficient asymptotic therefore follows from~\eqref{eq:boundarypartition}. Equation~\eqref{eq:defNhatboundary} gives the distributional identity
	\[
		\widehat{N}_n \, \eqdist \, \left(\widehat{N} \mid \widehat{N} \le m_n \right).
	\]
	Hence, because the (common) state space on which $\hat{N}_n$ and $\hat{N}$ are defined is countable, and $H(\infty) < \infty$ by assumption,
	\begin{equation}\label{eq:distNnNhat}
		d_{\mathrm{TV}}(\widehat{N}_n,\widehat{N}) = \Prb{\widehat{N}>m_n} \longrightarrow 0.
	\end{equation}
	The vectors $\boldsymbol{\Xi}_n$ and $\boldsymbol{\Xi}^{\mathrm{g}}_n$ use the same construction (conditional on their lengths). Given $\widehat{N} = k$, respectively $\widehat{N}_n = k$, sample $k-1$ independent copies of $X$ and insert $n-S_{k-1}$ at an independent uniform position in $[k]$. Now, Equation~\eqref{eq:distNnNhat} gives
	\[
		d_{\mathrm{TV}}(\boldsymbol{\Xi}^{\mathrm{g}}_n, \boldsymbol{\Xi}_n) \le d_{\mathrm{TV}}(\widehat{N}_n, \widehat{N})  \longrightarrow 0.
	\]
	Using Theorem~\ref{thm:boundarycondensation},
	\[
		d_{\mathrm{TV}} \left((K_1, \ldots, K_{N_n}),\boldsymbol{\Xi}_n \right) \le d_{\mathrm{TV}} \left((K_1, \ldots, K_{N_n}),\boldsymbol{\Xi}^{\mathrm{g}}_n \right) + \Prb{\widehat{N} > m_n} \longrightarrow 0.
	\]
\end{proof}

\subsection{Fluctuations and random scales in divergent condensation}\label{sec:fluct}

Under the assumptions of Theorem~\ref{thm:boundarycondensation}, when $H(\infty) = \infty$ the number of components tends to infinity in probability while $N_n/m_n \convp 0$. We study the secondary components, meaning all components other than the giant. Their ranked sizes are $K_{(2)}, K_{(3)}, \ldots$, whereas $\tau(P_n)$ lists the same remaining sizes in their original order. Their sizes are normalised by a random $a_{N_n}$, and their total mass is centered by a random $b_{N_n}$. Related random centerings and scales occur in the study of trimmed Lévy processes by Ipsen, Maller and Resnick~\cite{imr20}.

For i.i.d. samples, the order-statistic approach to stable limits goes back to LePage, Woodroofe and Zinn~\cite{lepageWoodroofeZinn1981}. We use the point-process approach of~\cite[Theorems~6.3 and~7.1]{resnick07} to obtain the joint limit of the centered sum and the rescaled observations. The
total-variation approximation transfers this limit to the secondary components of the Gibbs partition at the random scale $a_{N_n}$, jointly with the transformed component count.

First, set
\[
	U_n := \frac{H(N_n \wedge m_n)}{H(m_n)}, \qquad \cN_n^{\mathrm{res}} := \sum_{\substack{2 \le i \le N_n \\K_{(i)} > 0}} \delta_{K_{(i)}/a_{N_n}}.
\]
We use the Poisson point process $\cN = \sum_{j \ge 1} \delta_{J_j}$  introduced in Corollary~\ref{cor:largest}, with intensity measure given by $\nu(\diff x) = x^{-2} \diff x$. The theorem below shows that this same point-process law describes the secondary component sizes after rescaling by the random quantity $a_{N_n}$.

For $0 < \varepsilon < 1$, let
\begin{equation}\label{eq:Yeps}
	Y_\varepsilon := \sum_{J_j>\varepsilon}J_j-\log(1/\varepsilon).
\end{equation}
The compensator is the mean contribution of the atoms in $]\varepsilon,1]$; see~\cite[Section~12.1]{lastPenrose2017}, that is
\begin{equation*}
	\E{\sum_{\varepsilon < J_j \le 1} J_j} = \int_{\varepsilon}^{1} x \nu(\diff x) = \log(1/\varepsilon).
\end{equation*}

The variables $Y_\varepsilon$ are the approximations at time $t=1$ in the Lévy--Itô representation with Lévy measure $\nu$, with threshold at~$1$, as also uncovered in~\cite[Equation~(5.36)]{resnick07}. For the reader's convenience, we give a proof of the next lemma, which records the error due to small jumps in $Y_\varepsilon$ and identifies the limit with the law in~\eqref{eq:chara}.

\begin{lemma}[Poisson representation of the Cauchy limit]
	\label{lem:meso-poisson-cauchy}
	For $Y_{\varepsilon}$ defined in~\eqref{eq:Yeps}, there exists a random variable $Y$, entirely determined by the points of $\cN$, and such that
	\begin{equation}\label{eq:meso-compensated-sum}
		Y = \lim_{\varepsilon \rightarrow 0}Y_\varepsilon \quad\hbox{in probability}.
	\end{equation}
	Moreover,
	\begin{equation}\label{eq:meso-poisson-small-jumps}
		\E{Y-Y_\varepsilon} =0, \qquad \E{(Y-Y_\varepsilon)^2} = \varepsilon,
	\end{equation}
	and, for $t\ne0$,
	\begin{equation}\label{eq:meso-poisson-characteristic}
		\E{e^{itY}} = \exp\left\{-\frac{\pi}{2}|t| + it\left(1-\gamma_E-\log|t|\right)\right\}.
	\end{equation}
	Hence, $Y$ has exactly the asymmetric Cauchy law specified in~\eqref{eq:chara}.
\end{lemma}

\begin{proof}
	Because $\cN(]1,\infty]) < \infty$, and  $\nu(\{\infty\}) = 0$, it holds that
	\[
		\sum_{J_j > 1} J_j < \infty \text{ a.s. }
	\]
	Hence, we can decompose $Y_{\varepsilon}$ as
	\[
		Y_{\varepsilon} = \sum_{J_j > 1} J_j + C_\varepsilon := \sum_{J_j > 1} J_j+ \left( \sum_{\varepsilon<J_j\le 1} J_j -\int_\varepsilon^1x\,\nu(\diff x)\right).
	\]
    For $0 < \delta < \varepsilon < 1$, the mean and variance formulas for Poisson integrals~\cite[Section~12.1]{lastPenrose2017} give
    \[
        \E{C_\delta-C_\varepsilon}=0, \qquad \E{(C_\delta-C_\varepsilon)^2} = \int_{\delta}^{\varepsilon}  x^2\,\nu(\diff x) = \varepsilon-\delta.
    \]
    It follows that $(C_\varepsilon)$ is Cauchy in $L^2$ as $\varepsilon \rightarrow 0$. Let $C$ be its $L^2$ limit. Since every $C_\varepsilon$ is $\sigma(\cN)$-measurable, $C$ has a $\sigma(\cN)$-measurable version. Letting $\delta \rightarrow 0$ in the above identities yields
    \[
        \E{C-C_\varepsilon}=0,\qquad \E{(C-C_\varepsilon)^2}=\varepsilon.
    \]
    Now, the random variable
    \[
        Y := \sum_{J_j > 1} J_j + C
    \]
    satisfies both~\eqref{eq:meso-compensated-sum} and~\eqref{eq:meso-poisson-small-jumps}, and is fully determined by $\cN$. Note that the $L^2$ assertion only applies to the compensated small atoms, and the limit $Y$ has no finite second moment. It remains to identify its distribution.

	By the finite intensity Poisson representation, the sum $\sum_{J_j > \varepsilon} J_j$ has a compound Poisson distribution. Its characteristic function satisfies (see e.g.~\cite[Theorem~3.4]{schilling16})
	\begin{align*}
		\E{e^{itY_\varepsilon}} &= e^{-it \log(1/\varepsilon)} \E{e^{it \sum_{J_j > \varepsilon} J_j}} \\
		&= \exp\left\{-it \log(1/\varepsilon) + \int_{\varepsilon}^{\infty} (e^{itx}-1) \nu(\diff x) \right\} \\
		&= \exp\left\{\int_\varepsilon^\infty \left(e^{itx}-1-itx\mathbf{1}_{\{x \le 1\}}\right) x^{-2}\,\diff x\right\}.
	\end{align*}
	For fixed $t$, the compensated integrand is bounded near zero, since $|e^{itx}-1-itx| \le t^2x^2/2$, and is $O(x^{-2})$ at infinity. Dominated convergence hence allows us to write
	\begin{align*}
		\Psi(t) := \log\E{e^{itY}} &= \int_0^\infty \left(e^{itx}-1-itx\mathbf{1}_{\{x \le 1\}}\right)x^{-2}\,\diff x \\ &= \int_0^{\infty} \frac{\cos(tx)-1}{x^2} \diff x + i \int_0^{\infty} \frac{\sin(tx) - tx \I{x \le 1}}{x^2} \diff x. \
	\end{align*}
	Now, for $t > 0$, a calculation similar to Janson~\cite[Example~7.7]{janson2022stabledistributions} gives
    \[
        \int_{0}^{\infty} \frac{\cos(tx)-1}{x^2} \diff x = -\frac{\pi t}{2},
        \qquad \int_{0}^{\infty} \frac{\sin(tx)-tx\mathbf{1}_{\{x\le1\}}}{x^2}\diff x = t(1-\gamma_E-\log t).
    \]
    Consequently,
	\[
		\Psi(t) = -\frac{\pi t}{2} + it(1-\gamma_E - \log t),
	\]
	and, using complex conjugation for $t < 0$ to obtain a similar expression, yields overall
	\[
		\E{e^{itY}} = \exp\left\{-\frac{\pi}{2} |t| + it(1-\gamma_E - \log |t|)\right\},
	\]
	which proves \eqref{eq:meso-poisson-characteristic}. Comparison with~\eqref{eq:chara} identifies the law of $Y$.
\end{proof}

\begin{remark}
For fixed $k$, the stationary law of a zero-range process conditioned on total mass $n$ has the product form~\eqref{eq:KasX}, with sites corresponding to components and occupations to their
sizes. Armend\'ariz, Grosskinsky and Loulakis~\cite{armendarizGrosskinskyLoulakis2013} study the exceedances above the typical total mass in the finite-mean case. At sufficiently high densities, Xu~\cite{xu2024condensation} obtains a product approximation after removing the largest occupation, including for weights proportional to $j^{-2}$ for $j\ge 1$.

A conditioned product law with a random number of positive components also appears in renewal models. Godrèche~\cite[Sections~4.1 and~8.1]{godreche2021} studies geometric number of components' weights and gives corresponding convergent-condensation asymptotics. In our setting, the tilted, regularly varying outer weights dominate the competition between counts $k = m_n + O(c_n)$ and $k=o(m_n)$.
\end{remark}

The following theorem identifies the joint limit of the residual-mass fluctuations and the point process of secondary component sizes as the pair $(Y,\cN)$, constructed in Lemma~\ref{lem:meso-poisson-cauchy}.

\begin{theorem}[Divergent condensation: the giant and the secondary point process]
	\label{thm:divergent condensationtv}
	Suppose that the assumptions of Theorem~\ref{thm:boundarycondensation} hold. That is, we assume~\eqref{eq:wasymptotic},~\eqref{eq:critical},~\eqref{eq:vboundary} and~\eqref{eq:geldominance}. Suppose that $H(\infty) = \infty$. Then, the number of components satisfies
	\begin{equation}\label{eq:mesonumbergrowth}
		N_n \convp \infty, \qquad \frac{N_n}{m_n} \convp 0.
	\end{equation}
	The largest component is unique with probability tending to one, and
	\begin{equation}\label{eq:mesolargestfirstorder}
		\frac{K_{(1)}}{n}\convp 1, \qquad \frac{n-K_{(1)}}{b_{N_n}} \convp 1.
	\end{equation}
	Moreover, the following joint limit holds
	\begin{equation}\label{eq:meso-joint}
		\left(U_n,\frac{n-K_{(1)}-b_{N_n}}{a_{N_n}}, \cN_n^{\mathrm{res}} \right) \convdis (U,Y,\cN)
	\end{equation}
	in $[0,1]\times \ndR \times \cM_p(]0,\infty])$. Here $U$ is uniform on $[0,1]$, independent of $\cN$, and $Y$ is constructed from $\cN$ in Lemma~\ref{lem:meso-poisson-cauchy}.
\end{theorem}
\begin{proof}
	For fixed $M$ and $\varepsilon>0$, the auxiliary number of components satisfies
	\[
		\Prb{\widehat{N}_n \le M} = \frac{H(M)}{H(m_n)}\to0,\qquad \Prb{\widehat{N}_n>\varepsilon m_n} =1-\frac{H(\varepsilon m_n)}{H(m_n)} \rightarrow 0.
	\]
	Equation~\eqref{eq:boundarynumbertv} transfers these conclusions to $N_n$, proving~\eqref{eq:mesonumbergrowth}. Since $k v_k W(\rho_w)^k \sim L_v(k)/k \rightarrow 0$, this sequence is bounded, and
	\begin{equation}\label{eq:widehatNmax}
		\max_{k \le m_n} \Prb{\widehat{N}_n=k} \longrightarrow 0.
	\end{equation}
	Let $F_n(k) := H(k)/H(m_n)$ for $1 \le k \le m_n$, and $\widehat{U}_n := F_n(\widehat{N}_n)$. The discrepancy between the distribution function of $F_n(\widehat{N}_n)$ and the uniform distribution function is bounded by the largest jump of $F_n$. Then, Equation~\eqref{eq:widehatNmax} gives $\widehat{U}_n \convdis U$, where $U$ is uniform on $[0,1]$.

	We next prove the joint limit for an independent sample. Let
	\[
		\cN^{(k)} := \sum_{\substack{1 \le i \le k-1\\X_i>0}} \delta_{X_i/a_k}.
	\]
	We claim that
	\begin{equation}\label{eq:meso-iid-joint}
		Z_k := \left(\frac{S_{k-1}-b_k}{a_k},\,\cN^{(k)}\right) \convdis (Y,\,\cN) =: Z, \quad\text{ as } k\rightarrow \infty.
	\end{equation}
	Regular variation and~\eqref{eq:defan2} give $(k-1)\Prb{X>a_kx} \rightarrow x^{-1}$ for every $x > 0$. Resnick's independent-sample point process theorem~\cite[Theorem~6.3]{resnick07} gives
	\begin{equation}\label{eq:meso-iid-point-process}
		\cN^{(k)} \convdis \cN.
	\end{equation}

	Equations~\eqref{eq:common-truncated-power} and~\eqref{eq:common-truncated-first} give, for fixed $0 < \varepsilon < 1$,
	\begin{align}\label{eq:meso-truncated-moments}
		\frac{k}{a_k^2}\E{X^2; X \le \varepsilon a_k} &\longrightarrow \varepsilon,& \frac{k}{a_k}\E{X; \varepsilon a_k < X \le a_k} &\longrightarrow \log(1/\varepsilon).
	\end{align}
	Let $Q_{k,\varepsilon}$ denote the centered and rescaled contribution from the terms with $X_i \le \varepsilon a_k$:
	\[
		Q_{k,\varepsilon} := \frac{1}{a_k}\sum_{i=1}^{k-1} \left(X_i\mathbf{1}_{\{X_i \le \varepsilon a_k\}}-\E{X; X \le \varepsilon a_k}\right).
	\]
	Chebyshev's inequality gives, for $\eta > 0,$
	\begin{equation}\label{eq:meso-scalar-small-jumps}
		\limsup_{k \rightarrow \infty} \Prb{|Q_{k,\varepsilon}|>\eta} \le \limsup_{k \rightarrow \infty} \frac{\Va{Q_{k,\varepsilon}}}{\eta^2} \le \limsup_{k \rightarrow \infty} \frac{k}{\eta^2 a_k^2} \E{X^2; X \le \varepsilon a_k} \le \frac{\varepsilon}{\eta^2}.
	\end{equation}
	   Split the sum according to the large $X_i/a_k > \varepsilon$ and small contributions $X_i/a_k \le \varepsilon$. Splitting its centering at the same threshold gives
	\begin{align*}
		\frac{S_{k-1}-b_k}{a_k} &= \int_{]\varepsilon,\infty]}x\,\cN^{(k)}(\diff x) -\frac{k}{a_k}\E{X;\varepsilon a_k<X \le a_k} + Q_{k,\varepsilon} - \frac{\mu(\varepsilon a_k)}{a_k}.
	\end{align*}
	   For each fixed $\varepsilon$, the final term tends to zero as $k \rightarrow \infty$ by~\eqref{eq:muan}. For the contribution $Q_{k,\varepsilon}$, which is centered, Equation~\eqref{eq:meso-scalar-small-jumps} gives, for every $\eta > 0$,
	\[
		\lim_{\varepsilon \rightarrow 0}\limsup_{k \rightarrow \infty} \Prb{|Q_{k,\varepsilon}| > \eta} =0.
	\]

	For a fixed $\varepsilon$, summing the atoms above $\varepsilon$ is continuous at $\cN$ almost surely, by the atom-matching argument in the proof of Corollary~\ref{cor:largest}. Appeal to the continuous mapping theorem~\cite[Theorem~3.1]{resnick07} to~\eqref{eq:meso-iid-point-process}, keeping the point process coordinate and use~\eqref{eq:meso-truncated-moments}, to obtain
	\[
		\left(\int_{]\varepsilon,\infty]} x \cN^{(k)}(\diff x) -\frac{k}{a_k}\E{X;\varepsilon a_k < X \le a_k},\, \cN^{(k)}\right) \convdis (Y_{\varepsilon},\cN).
	\]
	Now~\eqref{eq:meso-scalar-small-jumps} and \eqref{eq:meso-poisson-small-jumps}, followed by letting $\varepsilon \rightarrow 0$, prove~\eqref{eq:meso-iid-joint} through the second-converging together theorem \cite[Theorem~3.5]{resnick07}. Therefore, Equation~\eqref{eq:meso-iid-joint} holds with $Y$ equal to the compensated sum of the atoms of $\cN$.

	We now evaluate~\eqref{eq:meso-iid-joint} at the independent random index $\widehat{N}_n$, jointly with $\widehat{U}_n$. For a bounded continuous function $g$,~\eqref{eq:meso-iid-joint} gives
	\[
		\E{g\left(Z_k\right)} \rightarrow \E{g(Z)}.
	\]
	Since $\widehat{N}_n$ is independent of $(X_i)_{i \ge 1}$, conditioning on $\widehat{N}_n$ gives, for every bounded continuous function $f$ on $[0,1],$
	\[
		\E{f(\widehat{U}_n)g\left(Z_{\widehat{N}_n}\right)} = \sum_{k = 1}^{m_n}\Prb{\widehat{N}_n=k} f\left(\frac{H(k)}{H(m_n)}\right) \E{g(Z_k)}.
	\]

	Splitting the above sum at $M$, and using $\lvert \E{g(Z_k)} - \E{g(Z)} \rvert \le 2 \lVert g\rVert_{\infty}$, we obtain
	\begin{multline*}
		\left| \E{f(\widehat{U}_n)g(Z_{\widehat{N}_n})} -\E{f(\widehat{U}_n)\E{g(Z)}}\right| \\ \le \|f\|_\infty\left(\sup_{k \ge  M}|\E{g(Z_k)}-\E{g(Z)}|+2\|g\|_\infty\Prb{\widehat{N}_n < M}\right).
	\end{multline*}
	Letting first $n \rightarrow \infty$ and then $M \rightarrow \infty$, and using $\widehat{U}_n \convdis U$ as shown above, this yields
	\[
		\E{f(\widehat{U}_n)g\left(Z_{\widehat{N}_n}\right)} \longrightarrow \E{f(U)}\E{g(Z)}.
	\]
	Taking $f \equiv 1$ gives $Z_{\widehat{N}_n}\convdis Z$. Together with $\widehat{U}_n \convdis U$, this implies tightness of the joint laws. The limit obtained by the product of test functions above identifies every subsequential limit with the product law of $U$ and $Z$. In turn, this gives
	\begin{equation}\label{eq:jointlimitUYN}
		\left(\widehat{U}_n,\,\frac{S_{\widehat{N}_n-1}-b_{\widehat{N}_n}}{a_{\widehat{N}_n}},\, \cN^{(\widehat{N}_n)}\right) \convdis (U,Y,\cN),
	\end{equation}
	where $U$ is independent of $(Y,\cN)$.

     We have proved the joint limit for the independent sample with a random number $\widehat{N}_n-1$ of terms. To obtain the corresponding limit for the partition, we apply the same measurable maps to the two triples compared in~\eqref{eq:boundarymarkedtv}. For $\mathbf{x} \in \ndN_0^{k-1}$, write $\mathbf{x}^{\downarrow}$ for the vector of coordinates rearranged in decreasing order. On the countable space $\bigsqcup_{k \ge 1}(\{k\} \times \ndN_0^{k-1} \times[k])$, consider the map
    \[
        (k,\mathbf{x},j) \mapsto \left(k, \sum_{i=1}^{k-1}x_i, \mathbf{x}^{\downarrow}\right). 
    \]
    The sum and ranked vector are interpreted as $0$ and $\emptyset$ when $k = 1$. Applying this map to~\eqref{eq:boundarymarkedtv}, and then using~\eqref{eq:tv-kernel-contract}, gives
    \begin{multline}\label{eq:boundaryjointtv}
        d_{\mathrm{TV}} \Big(\big(N_n,n-K_{(1)},(K_{(2)},\ldots,K_{(N_n)})\big), \big(\widehat{N}_n,S_{\widehat{N}_n-1},(X_{(1):\widehat{N}_n-1},\ldots,X_{(\widehat{N}_n-1):\widehat{N}_n-1})\big)\Big) \longrightarrow 0.
    \end{multline}
	For $\mathbf{r} =(r_1, \ldots, r_{k-1})$, consider further the measurable map given by
	\[
		T_n(k,s,\mathbf{r}) := \left(\frac{H(k\wedge m_n)}{H(m_n)},\,\frac{s-b_k}{a_k},\, \sum_{\substack{1 \le i \le k-1 \\ r_i>0}} \delta_{r_i/a_k} \right).
	\]
    Applying $T_n$ to the two vectors in~\eqref{eq:boundaryjointtv} gives the vectors in~\eqref{eq:meso-joint} and~\eqref{eq:jointlimitUYN}. Their total-variation distance tends to zero by~\eqref{eq:tv-kernel-contract}. Since the second vector converges in distribution, the first does too, and with the same limit. This proves~\eqref{eq:meso-joint}.
    
	It remains to show the first-order relations~\eqref{eq:mesolargestfirstorder}. We first compare $b_{N_n}$ with $n$. On $\{N_n \le m_n \}$, monotonicity of both $(a_{\cdot})$ and $\mu$, and the definition of $(b_{\cdot})$, give
	\[
		\frac{b_{N_n}}{n} = \frac{N_n}{m_n} \frac{\mu(a_{N_n})}{\mu(a_{m_n})} \frac{b_{m_n}}{n} \le \frac{N_n}{m_n}\frac{b_{m_n}}{n}.
	\]
	By~\eqref{eq:mesonumbergrowth}, $N_n/m_n \convp 0$ and $\Prb{N_n>m_n} \rightarrow 0$, while $b_{m_n}/n \rightarrow 1$ by~\eqref{eq:bdiff2co}. Hence $b_{N_n}/n \convp 0$.

	Next, by~\eqref{eq:defan2} and~\eqref{eq:omun},
	\[
		\frac{a_k}{b_k} = \frac{a_k}{k\mu(a_k)} \sim\frac{L(a_k)}{\mu(a_k)} \longrightarrow 0.
	\]
	Since $N_n \convp \infty$, it follows that $a_{N_n}/b_{N_n} \convp 0$. Moreover, \eqref{eq:meso-joint} implies tightness of
	\[
		\frac{n-K_{(1)}-b_{N_n}}{a_{N_n}}.
	\]
	It follows that
	\[
		\frac{n-K_{(1)}}{b_{N_n}} = 1+\frac{a_{N_n}}{b_{N_n}}\frac{n-K_{(1)}-b_{N_n}}{a_{N_n}} \convp 1.
	\]
	Combining this with $b_{N_n}/n \convp 0$ gives
	\[
		\frac{K_{(1)}}{n} = 1-\frac{b_{N_n}}{n} \frac{n-K_{(1)}}{b_{N_n}} \convp 1.
	\]
	In particular, $\Prb{K_{(1)}>n/2} \rightarrow 1$. Since the component sizes sum up to $n$, this event implies that the largest component is unique, completing the proof.
\end{proof}

Ranking the Poisson points gives the limiting secondary sizes. For fixed $r \ge 0$, removing the giant and the next $r$ largest components leaves the mass $n - \sum_{i=1}^{r+1}K_{(i)}$. After centering by $b_{N_n}$ and dividing by $a_{N_n}$, its limit is obtained by subtracting $J_1 + \cdots + J_r$ from $Y$. This removal of a fixed number of the largest terms is called trimming. Continuity of the corresponding operators and functional limits for trimmed Lévy processes are studied by Buchmann, Ipsen and Maller~\cite{zbMATH06943599}. For related small-time limits of ratios involving a trimmed Lévy process and its largest jumps, see~\cite[Theorem~2.1]{ipsenKeveiMaller2018}.

\begin{corollary}[Ranked secondary sizes and residual mass after trimming]
	\label{cor:meso-ranked}
	Under the assumptions of Theorem~\ref{thm:divergent condensationtv}, for every fixed $r \ge 1$,
	\begin{equation}\label{eq:meso-fixed-ranks}
		\left(\frac{K_{(2)}}{a_{N_n}},\ldots, \frac{K_{(r+1)}}{a_{N_n}}\right) \convdis (J_1,\ldots,J_r) \eqdist(\Gamma_1^{-1},\ldots,\Gamma_r^{-1}),
	\end{equation}
	jointly with~\eqref{eq:meso-joint}. For every fixed $r \ge 0$, also jointly with that limit,
	\begin{equation}\label{eq:meso-trimmed}
		\frac{n-\sum_{i=1}^{r+1}K_{(i)}-b_{N_n}}{a_{N_n}} \convdis Y-\sum_{j=1}^rJ_j.
	\end{equation}
	In particular, for every fixed $j \ge 1$ and $x > 0$,
	\[
		\Prb{\frac{K_{(j+1)}}{a_{N_n}} \le x} \longrightarrow e^{-1/x}\sum_{\ell=0}^{j-1}\frac{x^{-\ell}}{\ell!}.
	\]
\end{corollary}

\begin{proof}
	The continuity argument in Remark~\ref{rem:atom-matching}, applied jointly to~\eqref{eq:meso-joint}, proves~\eqref{eq:meso-fixed-ranks}. Subtracting the first $r$ ranked coordinates from the scalar coordinate gives~\eqref{eq:meso-trimmed}, including $r=0$ with the empty-sum convention. The distribution of $J_j$ was computed in Corollary~\ref{cor:largest}.
\end{proof}

Next, we consider the entire ranked sequence of secondary sizes. Convergence of any fixed number of its coordinates does not, by itself, control the sum over all components. For $q > 1$, the next corollary gives convergence in $\ell^q$ and, consequently, convergence of the sums of $q$th powers. Define
\[
	\mathbf{K}_n := \left(\frac{K_{(j+1)}}{a_{N_n}}\right)_{j \ge 1}.
\]

\begin{corollary}[Ranked sequences in $\ell^q$ and power sums]
	\label{cor:meso-functionals}
	Under the assumptions of Theorem~\ref{thm:divergent condensationtv}, for every $q>1$,
	\begin{equation}\label{eq:meso-ranked-lq}
		\mathbf{K}_n \convdis (J_j)_{j \ge 1}\, \text{ in } \ell^q,
	\end{equation}
	with respect to the norm topology, jointly with~\eqref{eq:meso-joint}. Consequently, jointly with these convergences,
	\begin{equation}\label{eq:meso-power-sums}
		\frac{\sum_{i \ge 2} K_{(i)}^q}{a_{N_n}^q} \convdis \sum_{j \ge 1} J_j^q.
	\end{equation}
	The limit has Laplace transform given, for $s \ge 0$, by
	\begin{equation}\label{eq:meso-power-laplace}
		\E{\exp\left(-s\sum_{j \ge 1}J_j^q\right)} = \exp\left(-\Gamma(1-1/q)s^{1/q}\right).
	\end{equation}
\end{corollary}

\begin{proof}
	Fix $q>1$ and let $\mathbf{J}=(J_j)_{j \ge 1}$. Campbell's formula~\cite[Proposition~2.7]{lastPenrose2017} gives, for every $\varepsilon > 0$,
	\begin{equation}\label{eq:CampbJ}
		\E{\sum_{j :J_j \le \varepsilon} J_j^q} = \int_0^\varepsilon x^{q-2}\,\diff x = \frac{\varepsilon^{q-1}}{q-1}.
	\end{equation}
	Taking $\varepsilon = 1$ in the above, and recalling that $\cN$ has only finitely many atoms above $1$ and none at $\infty$, we obtain $\mathbf{J} \in \ell^q$ almost surely. Moreover, The Laplace functional~\cite[Theorem~3.9]{lastPenrose2017} yields, for $s \ge 0$,
	\begin{align*}
		\E{\exp\left(-s\sum_{j \ge 1}J_j^q \right)} = \exp\left(-\int_0^\infty (1-e^{-sx^q})x^{-2}\diff x\right) = \exp\left(-\Gamma(1-1/q) s^{1/q}\right).
	\end{align*}
	where the last equality follows by the substitution $y=x^q$ and integration by parts. This proves~\eqref{eq:meso-power-laplace}. Hence, $\sum_{j\ge 1} J_j^q$ has a positive stable law of index $1/q$, see~\cite[Section~1.2]{bertoin1999}

	Now, set $\mathbf{X}_k := (X_{(j):k-1}/a_k)_{j \ge 1}$. Let $\mathbf{X}_k^{>\varepsilon}$ and $\mathbf{J}^{>\varepsilon}$ keeping the entries strictly greater than $\varepsilon$, with all other entries replaced by zero. For fixed $\varepsilon > 0$, the vectors $\mathbf{X}_k^{>\varepsilon}$ and $\mathbf{J}^{>\varepsilon}$ are obtained by ranking the atoms above $\varepsilon$ of $\cN^{(k)}$ and $\cN$. The continuity argument in Remark~\ref{rem:atom-matching} shows that this operation is continuous at $\cN$ almost surely.

	The continuous mapping theorem~\cite[Theorem~3.1]{resnick07}, applied jointly to~\eqref{eq:meso-iid-joint}, gives
    \[
        \left(\frac{S_{k-1}-b_k}{a_k},\cN^{(k)},\mathbf{X}_k^{>\varepsilon}\right) \convdis (Y,\cN,\mathbf{J}^{>\varepsilon}),
    \]
    with respect to the norm topology on the third coordinate.

	For the omitted, small entries,~\eqref{eq:common-truncated-power} gives
	\[
		\frac{k}{a_k^q}\E{X^q;X\le\varepsilon a_k} \longrightarrow \frac{\varepsilon^{q-1}}{q-1}.
	\]
	Markov's inequality gives, for every $\eta > 0$,
	\begin{equation}\label{eq:meso-ranked-small-jumps}
		\limsup_{k \rightarrow \infty} \Prb{\lvert\mathbf{X}_k-\mathbf{X}_k^{>\varepsilon}\rvert_q>\eta} \le \frac{\varepsilon^{q-1}}{(q-1)\eta^q}.
	\end{equation}
	The identity obtained from Campbell's formula~\eqref{eq:CampbJ} above gives the corresponding bound for the limiting sequence discrepancy:
	\[
		\Prb{\lvert\mathbf{J}-\mathbf{J}^{>\varepsilon}\rvert_q>\eta} \le \frac{\varepsilon^{q-1}}{(q-1)\eta^q}.
	\]
	Both truncation bounds tend to zero as $\varepsilon \rightarrow 0$. The second converging-together theorem~\cite[Theorem~3.5]{resnick07} yields
    \[
        \left(\frac{S_{k-1}-b_k}{a_k},\cN^{(k)},\mathbf{X}_k\right) \convdis (Y,\cN,\mathbf{J})
    \]
    in $\ndR \times \cM_p(]0,\infty])\times\ell^q$. This proves the joint convergence for the independent sample. Conditioning on $\widehat{N}_n$, using its independence from $(X_i)_{i \ge 1}$ and $\widehat{N}_n \convp \infty$, gives the same joint limit, replacing $k$ by $\widehat{N}_n$. The argument with products of bounded continuous test functions used to prove~\eqref{eq:jointlimitUYN} also gives joint convergence with $\widehat{U}_n$, with its limit $U$ being independent of $(Y,\cN,\mathbf{J})$. To obtain the partition limit, apply the measurable map
    \[
     (k,s,\mathbf{r}) \mapsto \left(T_n(k,s,\mathbf{r}), \left(r_i/a_k\right)_{i \ge 1}\right)
    \]
    to the two vectors in~\eqref{eq:boundaryjointtv}, extending $\mathbf{r}$ with zeros. $T_n$ is the map defined in the proof of Theorem~\ref{thm:divergent condensationtv}. Its coordinates keeps the normalised number of components, the centered and rescaled residual mass, and the size point measure. Inequality~\eqref{eq:tv-kernel-contract} shows that the resulting total-variation distance tends to zero. This proves~\eqref{eq:meso-ranked-lq} jointly with~\eqref{eq:meso-joint}. Finally, continuity of $\mathbf{x} \mapsto \lvert \mathbf{x}\rvert_q^q$ on $\ell^q$ gives~\eqref{eq:meso-power-sums}.
\end{proof}

Note that the restriction $q > 1$ is actually sharp. At the extreme-value scale, $\int_{0}^{1} x^q\nu(\diff x) < \infty$ exactly when $q > 1$. By contrast, $(J_j) \eqdist (\Gamma_j^{-1})$ and $\Gamma_j/j \rightarrow 1$ imply $\sum_j J_j = \infty$ almost surely. The uncentered residual mass has first-order scale $b_{N_n}$, which is larger than $a_{N_n}$, since
\[
	\|\mathbf{K}_n \|_1 = \frac{n-K_{(1)}}{b_{N_n}}\frac{b_{N_n}}{a_{N_n}} \convp \infty.
\]
The main theorem normalises the secondary sizes by the random quantity $a_{N_n}$. The next corollary shows that, for each fixed $i \ge 2$, there is no deterministic sequence $\gamma_n > 0$ such that $K_{(i)}/\gamma_n$ converges in distribution to a positive, finite, nondegenerate random variable.

\begin{corollary}
	\label{cor:meso-scale}
	Under the assumptions of Theorem~\ref{thm:divergent condensationtv}, recall from Section~\ref{sec:background} that $a_k = a(k)$, where $a$ is strictly increasing with inverse function $A$. Extend $A$ by $A(x) =0$ for $0 \le x \le 1$, and set $H(x) =0$ for $x < 1$. For every fixed $r \ge 1$, jointly with~\eqref{eq:meso-joint},
	\[
		\left(\frac{H(A(K_{(2)})\wedge m_n)}{H(m_n)},\ldots,\frac{H(A(K_{(r+1)})\wedge m_n)}{H(m_n)}\right)
		\convdis (U,\ldots,U).
	\]
	For each fixed $i \ge 2$, no deterministic sequence $\gamma_n > 0$ makes $K_{(i)}/\gamma_n$ converge in distribution to a nondegenerate positve random variable. More precisely, for $0 < u < 1$, define
	\[
		k_n(u) := \min \{k \in \ndN : k \le m_n,\ H(k)/H(m_n) \ge u\}.
	\]
	For every fixed $i \ge2$ and $x > 0$,
	\begin{equation}\label{eq:meso-deterministic-quantiles}
		\Prb{K_{(i)} \le x a_{k_n(u)}} \longrightarrow u.
	\end{equation}
	If, additionally, for some $\beta > 0$,
	\[
		\frac{H(x^t)}{H(x)} \longrightarrow t^\beta, \qquad 0 < t \le 1,
	\]
	then, for every fixed $r \ge 1$,
	\[
		\frac{1}{\log n} \left(\log N_n, \left(\log(K_{(j+1)} \vee 1)\right)_{1 \le j \le r} \right) \convdis (U^{1/\beta},U^{1/\beta},\ldots,U^{1/\beta}),
	\]
	jointly with~\eqref{eq:meso-joint}, where $U^{1/\beta}$ has the $\operatorname{Beta}(\beta,1)$ law. In particular, $U^{1/\beta}$ is independent of $(Y,(J_j)_{j \ge 1})$.
\end{corollary}

\begin{proof}
	Fix $i \ge 2$. By Corollary~\ref{cor:meso-ranked}, for every $\eta > 0$ we may choose $0 < c < C < \infty$ such that, for all sufficiently large $n$,
	\[
		\Prb{c\le K_{(i)}/a_{N_n} \le C} \ge1-\eta.
	\]
	Since $A(a_k) = k$, and $A$ is regularly varying with index 1, the uniform convergence theorem~\ref{thm:uct} gives
	\[
		\sup_{x\in [c,C]} \left|\frac{A(xa_k)}{k}-x\right| \longrightarrow 0.
	\]
	In particular, $A(xa_k)/k \in [c/2, 2C]$ for all sufficiently large $k$, uniformly for $x\in [c,C]$. Appealing to Equation~\eqref{eq:common-rv-window} to the slowly varying function $H$ on this explicit interval gives
	\[
		\sup_{x\in [c,C]} \left|\frac{H(A(xa_k))}{H(k)}-1\right| \longrightarrow 0.
	\]
	These (deterministic) asymptotics hold for every $k \ge M$, as soon as $M$ is sufficiently large. We can then substitute $k = N_n$ and $x = K_{(i)}/a_{N_n}$ on $\{N_n \ge M,\ c \le K_{(i)}/a_{N_n} \le C\}$. By~\eqref{eq:mesonumbergrowth}, the complementary event has probability (by a simple union bound) at most $\eta+o(1)$. Since $\eta>0$ is arbitrary, we obtain
	\[
		\frac{A(K_{(i)})}{N_n} = O_p(1), \qquad \frac{H(A(K_{(i)}))}{H(N_n)}\convp 1.
	\]
	Together with $N_n/m_n \convp 0$, the first estimate gives $A(K_{(i)})/m_n \convp 0$. Hence, $N_n \le m_n$ and $A(K_{(i)}) \le m_n$ with probability tending to one. On this event, the difference below equals
	\[
		\frac{H(N_n)}{H(m_n)} \left(\frac{H(A(K_{(i)}))}{H(N_n)}-1\right).
	\]
	The first factor is at most one and the second tends to zero in probability. Recalling the definition of $U_n$, we obtain
	\begin{equation}\label{eq:meso-scale-identification}
		\frac{H(A(K_{(i)})\wedge m_n)}{H(m_n)} - U_n \convp 0.
	\end{equation}
	A finite union bound,~\eqref{eq:meso-scale-identification}, and Slutsky's theorem (see e.g. \cite[Theorem~3.4]{resnick07}) give the joint convergence of these transformed ranks with the same variable $U$.

	We next prove~\eqref{eq:meso-deterministic-quantiles}. Fix $0 < u < 1$. Because $k_n(u)$ is minimal, and using the largest-atom bound~\eqref{eq:widehatNmax}, it follows that
	\[
		0 \le \frac{H(k_n(u))}{H(m_n)}-u \le \frac{H(k_n(u))-H(k_n(u)-1)}{H(m_n)} \le \max_{k \le m_n}\Prb{\widehat{N}_n=k} \longrightarrow 0.
	\]
	Hence $k_n(u) \rightarrow \infty$. Moreover, $k_n(u) = o(m_n)$, since otherwise, along a subsequence with $k_n(u) \ge \delta m_n$ for some $\delta > 0$, monotonicity and slow variation of $H$ would force the above ratio to tend to $1$. For fixed $x > 0$, the uniform convergence theorem~\ref{thm:uct} gives $A(xa_{k_n(u)}) \sim xk_n(u) = o(m_n)$. Therefore, for all sufficiently large $n$,
	\[
		\frac{H(A(xa_{k_n(u)})\wedge m_n)}{H(m_n)} = \frac{H(A(xa_{k_n(u)}))}{H(k_n(u))} \frac{H(k_n(u))}{H(m_n)} \longrightarrow u.
	\]
	Fix some $0 < \delta < \min \{u,1-u\}$. The above limit gives, for all sufficiently large $n$,
	\[
		u-\delta < \frac{H(A(xa_{k_n(u)}) \wedge m_n)}{H(m_n)} < u+\delta.
	\]
	Since the map $y \mapsto H(A(y)\wedge m_n)$ is nondecreasing,
	\begin{align*}
		\Prb{\frac{H(A(K_{(i)})\wedge m_n)}{H(m_n)}<u-\delta} \le \Prb{K_{(i)} \le xa_{k_n(u)}} \le \Prb{\frac{H(A(K_{(i)})\wedge m_n)}{H(m_n)} \le u+\delta}.
	\end{align*}
	By~\eqref{eq:meso-scale-identification} and~\eqref{eq:meso-joint}, the random variable in the outer probabilities converges in distribution to $U \sim \text{Unif}[0,1]$. Those probabilities therefore tend to $u-\delta$, respectively $u+\delta$. Letting $\delta \rightarrow 0$ proves \eqref{eq:meso-deterministic-quantiles}.

	Now fix $0 < u < 1$ and suppose that $K_{(i)}/ \gamma_n \convdis Z$ for deterministic $\gamma_n > 0$, where $0 < Z < \infty$ almost surely. Pass to a subsequence along which $\gamma_n/a_{k_n(u)} \rightarrow c \in [0,\infty]$. If $c = 0$, tightness of $K_{(i)}/\gamma_n$ gives $\Prb{K_{(i)} \le a_{k_n(u)}} \rightarrow 1$, whereas if $c = \infty$, positivity of $Z$ gives convergence of this probability to zero. Both conclusions contradict~\eqref{eq:meso-deterministic-quantiles}. If $0 < c < \infty$, bound $y\gamma_n/ a_{k_n(u)}$ between two fixed positive constants and use~\eqref{eq:meso-deterministic-quantiles}. At every positive continuity point $y$ of $Z$, this gives $\Prb{Z \le y} = u$, which is impossible for a probability distribution on $]0, \infty[$. This proves that a deterministic multiplicative normalisation is not possible if we are to converge to a proper (positive) limit.

	Finally, under the additional hypothesis of the Corollary, monotonicity in $t$ and continuity of $t^\beta$ give
	\[
		\sup_{0 \le t \le 1} \left| \frac{H(m_n^t)}{H(m_n)}-t^\beta \right|\longrightarrow 0.
	\]
	At $t = 0$, convergence follows from $H(1)/H(m_n) \rightarrow 0$. Choose a sufficiently fine and finite partition of $[0,1]$. Monotonicity bounds the ratio at each $t$ between its values at the neighboring partition points, while uniform continuity of $t^\beta$ bounds the difference between those (limiting) values. This proves uniform convergence. On $\{N_n \le m_n\}$, take $t = \log N_n/\log m_n$ and apply the continuous inverse $x\mapsto x^{1/\beta}$ to get
	\[
		\frac{\log N_n}{\log m_n}-U_n^{1/\beta} \convp 0.
	\]
	Potter bounds~\eqref{eq:pb} and Lemma~\ref{lem:seq}(e) give $\log a_k \sim \log k$ and $\log m_n \sim \log n$. Furthermore, Corollary~\ref{cor:meso-ranked} gives
	\[
		\log(K_{(i)} \vee 1)-\log a_{N_n} = O_p(1).
	\]
    	Since $N_n \convp \infty$ and $\Prb{N_n \le m_n} \rightarrow 1$, the relation $\log a_k \sim \log k$ gives
    \[
        \frac{\log a_{N_n}-\log N_n}{\log n} \convp 0.
    \]
    To see this, note that $\log N_n/ \log n$ is bounded in probability. Together with the above $O_p(1)$ estimate, this results in
    \[
        \frac{\log(K_{(i)} \vee 1)-\log N_n}{\log n} = \frac{\log(K_{(i)} \vee 1)-\log a_{N_n}}{\log n} + \frac{\log a_{N_n}-\log N_n}{\log n} \convp 0.
    \]
    Combining this with~\eqref{eq:meso-joint}, $\log m_n \sim \log n$, and the earlier approximation of $\log N_n/ \log m_n$ proves the claimed joint logarithmic limit and its independence from $(Y,(J_j)_{j \ge 1})$. 
\end{proof}

Recalling the definition of $\tau$ in~\eqref{eq:deletion-map}, let
\[
	(\widetilde{K}_1,\ldots,\widetilde{K}_{N_n-1}) := \tau(K_1,\ldots,K_{N_n}).
\]
The residual partial-sum process records these components in their original order. For $0 \le t \le 1$, we define
\[
	Z_n(t) := \frac{\sum_{j=1}^{\lfloor t(N_n-1)\rfloor}\widetilde{K}_j-t b_{N_n}}{a_{N_n}},
\]
and we let
\[
	\Pi_n := \sum_{\substack{1 \le j \le N_n-1\\\widetilde{K}_j>0}} \delta_{(j/(N_n-1),\,\widetilde{K}_j/a_{N_n})}.
\]
In this representation, the $j$th residual component is placed at time $j/(N_n-1)$. The partial sum process includes precisely the residual components, with assigned times at most $t$. When $N_n = 1$, the residual vector and point measure are empty, and the residual sums are zero.

Let $(T_j)_{j \ge 1}$ be independent uniform random variables on $[0,1]$, independent of $(U, \cN)$. $T_j$ is attached as the time of the size point $J_j$. Note that the index $j$ reflects the rank of the size, not its time order. Set
\[
    \Pi := \sum_{j \ge 1} \delta_{(T_j, J_j)}.
\]
By the marking theorem~\cite[Theorem~5.6]{lastPenrose2017}, $\Pi$ is a Poisson point process with intensity $\diff t\,x^{-2}\,\diff x$. We define, for $0 \le t \le 1$,
\begin{equation}\label{eq:meso-compensated-process}
	Y_t := \lim_{\varepsilon \rightarrow 0} \left(\sum_{\substack{T_j\le t\\J_j>\varepsilon}} J_j -t\log(1/\varepsilon)\right).
\end{equation}
The proposition below shows convergence in probability uniformly in $t$, and $Y_1 = Y$ with the Poisson representation in Lemma~\ref{lem:meso-poisson-cauchy}.

\begin{proposition}[Functional residual limit]
	\label{pro:meso-functional}
	Under the assumptions of Theorem~\ref{thm:divergent condensationtv},
	\begin{equation}\label{eq:meso-joint-functional}
		(U_n,Z_n,\Pi_n) \convdis (U,Y_{\cdot},\Pi)
	\end{equation}
	in
	\[
		[0,1]\times\mathbb{D}([0,1],\ndR)
		\times\cM_p([0,1]\times]0,\infty]),
	\]
	with respect to the $J_1$ topology for the second coordinate and the vague topology for the third. Here, $U$ is independent of $\Pi$, and $Y_{\cdot}$ is the Lévy process determined by~\eqref{eq:meso-compensated-process}, while its law at time $t = 1$ is~\eqref{eq:chara}.
	For every $q > 1$, this convergence holds jointly with \eqref{eq:meso-ranked-lq} in $\ell^q$.
\end{proposition}

\begin{proof}
    We use the truncation argument from the proof of \cite[Theorem~7.1]{resnick07}, keeping the times of the atoms. We follow the proof of~\eqref{eq:meso-iid-joint}, replacing the scalar estimates~\eqref{eq:meso-scalar-small-jumps} and~\eqref{eq:meso-poisson-small-jumps} by maximal inequalities.

	For $k = 1$, set $\Pi^{(1)} := 0$ and $Z^{(1)}(t) := -t b_1/a_1$. For $k \ge 2$, let
	\[
		Z^{(k)}(t) := \frac{S_{\lfloor t(k-1)\rfloor}-t b_k}{a_k}, \qquad \Pi^{(k)} := \sum_{\substack{1 \le j \le k-1\\X_j>0}} \delta_{(j/(k-1),X_j/a_k)}.
	\]
	Using the same arguments as in the proof of Theorem~\ref{thm:divergent condensationtv}, the independent-sample point-process theorem \cite[Theorem~6.3]{resnick07} yields $\Pi^{(k)} \convdis \Pi$.

	Fixing $q > 1$, we use the notation $\mathbf{X}_k$, $\mathbf{X}_k^{>\varepsilon},$ $\mathbf{J}^{>\varepsilon}$ and $\mathbf{J}$ from the proof of Corollary~\ref{cor:meso-functionals}. For $0 < \varepsilon < 1$, define the truncated (and its limit) paths
	\[
		\begin{split}
			Z^{(k,\varepsilon)}(t) & := \int_{[0,t]\times]\varepsilon,\infty]}x\,\Pi^{(k)}(\diff s,\diff x) - t\log(1/\varepsilon),\\
			Y_t^{(\varepsilon)} & := \int_{[0,t]\times]\varepsilon,\infty]} x\,\Pi(\diff s,\diff x) -t\log(1/\varepsilon).
		\end{split}
	\]
	With a completely similar reasoning as for $\cN$, the point process $\Pi$ has, almost surely, finitely many atoms above $\varepsilon$, with distinct times in $]0,1[$, and none with size $\varepsilon$ or $\infty$. The truncated summation map is continuous at $\Pi$ for the $J_1$ topology (see \cite[Section~7.2.3]{resnick07}). Combining this with the ranking argument in the proof of Corollary~\ref{cor:largest}, the continuous mapping theorem \cite[Theorem~3.1]{resnick07} gives
	\[
		\left(Z^{(k,\varepsilon)},\Pi^{(k)},\mathbf{X}_k^{>\varepsilon}\right) \convdis \left(Y_{\cdot}^{(\varepsilon)},\Pi,\mathbf{J}^{>\varepsilon}\right).
	\]

	It remains to control the small components. Splitting the sum and its centering similarly as in the proof of~\eqref{eq:meso-iid-joint} gives
	\begin{align*}
		\|Z^{(k)}-Z^{(k, \varepsilon)}\|_\infty
		&\le \max_{0 \le \ell \le k-1} \left| \frac{1}{a_k}\sum_{j=1}^{\ell} \left(X_j\mathbf{1}_{\{X_j \le \varepsilon a_k\}} -\mu(\varepsilon a_k)\right) \right|\\
		&\quad +\frac{2\mu(a_k)}{a_k} +\left| \frac{k}{a_k}\E{X;\varepsilon a_k<X \le a_k} -\log(1/\varepsilon) \right|.
	\end{align*}
	The last two terms tend to zero by~\eqref{eq:muan} and~\eqref{eq:meso-truncated-moments}. Kolmogorov's maximal inequality (see e.g. \cite[Equation~(5.43)]{resnick07}) and the first limit in~\eqref{eq:meso-truncated-moments} consequently yield, for $\eta > 0,$
	\[
		\limsup_{k \rightarrow \infty} \Prb{\|Z^{(k)}-Z^{(k,\varepsilon)}\|_\infty>\eta} \le \frac{\varepsilon}{\eta^2}.
	\]

	For the limit $Y_{\cdot}$, the Poisson construction of a Lévy process \cite[Sections~5.5.1 and~5.5.3]{resnick07} gives a càdlàg process $Y_{\cdot}$ obtained from $Y_{\cdot}^{(\varepsilon)}$ by letting $\varepsilon \rightarrow 0$ uniformly in probability. The difference $Y_{\cdot}-Y_{\cdot}^{(\varepsilon)}$ is a (centered) square-integrable Lévy process. It is obtained by considering, for a fixed $\delta$ such that $0 < \delta < \varepsilon < 1,$
	\[
		Y_t^{(\delta)}-Y_t^{(\varepsilon)} = \sum_{\substack{T_j \le t \\ \delta < J_j \le \varepsilon}} J_j - t \log(\varepsilon/\delta),
	\]
	and then letting $\delta \rightarrow 0.$ The Poisson integral isometry \cite[Lemma~12.2 and Proposition~12.4]{lastPenrose2017} states its variance at time $t = 1$ as
	\[
		\int_0^\varepsilon x^2x^{-2}\,\diff x= \varepsilon.
	\]
	The continuous version of Kolmogorov's maximal inequality (see e.g. \cite[Lemma~5.3]{resnick07}) gives
	\[
		\Prb{\|Y_{\cdot}-Y_{\cdot}^{(\varepsilon)}\|_\infty>\eta} \le \frac{\varepsilon}{\eta^2}.
	\]
	At time $1$, this is the construction of $Y$ in Lemma~\ref{lem:meso-poisson-cauchy}. In particular, $Y_1 = Y$ and its characteristic function is~\eqref{eq:chara}.

	For the ranked coordinates, the required truncation bounds are~\eqref{eq:meso-ranked-small-jumps} and its limiting counterpart in the proof of Corollary~\ref{cor:meso-functionals}. Appealing to the second converging-together theorem \cite[Theorem~3.5]{resnick07} to all three coordinates yields
	\[
		(Z^{(k)}, \Pi^{(k)}, \mathbf{X}_k) \convdis (Y_{\cdot}, \Pi, \mathbf{J})
	\]
	in $\mathbb{D}([0,1],\ndR) \times \cM_p([0,1]\times]0,\infty]) \times \ell^q$.

	Condition on $\widehat{N}_n$ and use its independence from $(X_i)_{i \ge 1}$, as in the proof of Theorem~\ref{thm:divergent condensationtv}. This gives the same joint limit with $k$ replaced by $\widehat{N}_n$, jointly with $\widehat{U}_n$. Once again, the limiting variable $U$ is independent of the point process $\Pi$. The comparison in~\eqref{eq:boundarymarkedtv} keeps the number of components and the remaining sizes in their original order. Applying the same partial-sum and point-measure maps to its two triples, and using~\eqref{eq:tv-kernel-contract}, this gives the corresponding joint limit for the partition. This proves~\eqref{eq:meso-joint-functional}, jointly with~\eqref{eq:meso-ranked-lq}.
\end{proof}

\section{The boundary phase transition}
\label{sec:phase}

We now compare the central and single-big-jump contributions without assuming~\eqref{eq:geldominance}. Assume~\eqref{eq:wasymptotic},~\eqref{eq:critical} and~\eqref{eq:vboundary}, and recall the scales introduced in Section~\ref{sec:intro}:
\[
	\cC_n^{\cD}  := \frac{v_{m_n} W(\rho_w)^{m_n}}{\mu(a_{m_n})}, \qquad \cC_n^{\cG}  := \Prb{X = n}H(m_n), \qquad \Theta_n := \frac{\cC_n^{\cG} }{\cC_n^{\cD} },
\]
and the regimes 
\[
	\cD_n := \{N_n \ge m_n/2\}, \qquad \cG_n := \{N_n<m_n/2\}.
\]
By~\eqref{eq:wasymptotic} and \eqref{eq:vboundary},
\[
	\Theta_n \sim \frac{H(m_n)L(n)m_n^2\mu(a_{m_n})}{n^2L_v(m_n)}.
\]
Since~\eqref{eq:defbn} and~\eqref{eq:bdiff2co} give $n \sim m_n \mu(a_{m_n})$, this simplifies to
\begin{equation}\label{eq:boundary-phase-ratio}
	\Theta_n \sim \frac{H(m_n)L(n)}{L_v(m_n)\mu(a_{m_n})}.
\end{equation}

The dense contribution comes from counts $k= m_n + O(c_n)$, and the condensation contribution from counts $k = o(m_n)$. Conditionally on $\cD_n$, the dense-phase limits hold. Conditionally on $\cG_n$, a giant emerges and, after removing it, the remaining components are approximated in total variation by $(X_1, \ldots, X_{\widehat{N}_n-1})$. The ratio $\Theta_n$ determines the relative probabilities of these two regimes. 

Bloznelis~\cite[Theorem~1(iii)]{MR4038060} gives a two-term local asymptotic for regularly varying summand and count laws when the summand is regularly varying with index below $-2$, the count $N$ has index at most $-2$, and $\E{N}<\infty$. Stufler~\cite[Theorem~3.12]{MR4780503} uses this asymptotic to describe a mixture of dense behaviour and convergent condensation when the two mass indices coincide below $-2$. A truncated first moment of the weights of $N$ also appears
in the finite-variance weighted renewal expansion of \cite[Theorem~3.2(iii) and Remark~3.1]{borovkovBorovkov2014}.

At the Cauchy boundary, Theorem~\ref{thm:b2} below gives, for $N$ as in~\eqref{eq:Nconv},
\[
    \Prb{S_N=n} \sim \frac{\Prb{N=m_n}}{\mu(a_{m_n})} +\Prb{X=n} \E{N;N\le m_n},
\]
with $\E{N;N\le m_n}= H(m_n)/V(W(\rho_w))$. When $H(\infty)=\infty$, the number of components given $\cG_n$ tends to infinity in probability.

\begin{theorem}[Boundary phase transition]
	\label{thm:b2}
	Suppose that~\eqref{eq:wasymptotic},~\eqref{eq:critical} and~\eqref{eq:vboundary} hold. Then
	\[
		u_n \rho_w^n \sim \cC_n^{\cD} +\cC_n^{\cG} ,
	\]
	and
	\[
		\Prb{\cD_n} \sim\frac{1}{1+\Theta_n}, \qquad \Prb{\cG_n} \sim\frac{\Theta_n}{1+\Theta_n}.
	\]
	Conditionally on $\cD_n$,
	\[
		\sup_{k \in \ndZ} \left| c_n \Prb{N_n = k \mid \cD_n} -f\left(-\frac{k-m_n}{c_n}\right) \right|\longrightarrow 0,
	\]
	and the conclusions of Lemma~\ref{lem:tv}, Theorem~\ref{thm:functional}, and Corollary~\ref{cor:largest} remain valid under the law given $\cD_n$. Conditionally on $\cG_n$,
	\begin{equation}\label{eq:b2gelconditional}
		d_{\mathrm{TV}} \left(\left((K_1,\ldots,K_{N_n})\mid \cG_n \right), \boldsymbol{\Xi}^{\mathrm{g}}_n \right) \longrightarrow 0.
	\end{equation}
	If $H(\infty) < \infty$, the distributional and total-variation conclusions of Corollary~\ref{pro:boundarycondensation} hold under the above conditional law. If $H(\infty) = \infty$, the conclusions of Theorem~\ref{thm:divergent condensationtv}, Corollaries~\ref{cor:meso-ranked},~\ref{cor:meso-functionals} and~\ref{cor:meso-scale}, and Proposition~\ref{pro:meso-functional} hold under the conditional law given $\cG_n$.

	If $\Theta_n \to\Theta \in [0,\infty]$, then the limiting behaviour is a $\mathrm{Bernoulli}(1/(1 +\Theta))$ mixture of these two regimes.
\end{theorem}
\begin{proof}
	Split $u_n \rho_w^n$ in two sums, and note that
	\[
		I_n^{\cD} := \sum_{k \ge m_n/2} v_k W(\rho_w)^k \Prb{S_k = n}, \quad I_n^{\cG}  := \sum_{k < m_n/2} v_k W(\rho_w)^k \Prb{S_k = n}, \quad u_n \rho_w^n = I_n^{\cD} + I_n^{\cG} .
	\]
	The estimates~\eqref{eq:almostllt},~\eqref{eq:i3} and \eqref{eq:i4} on $k \ge m_n/2$ remain valid at $b=2$. For fixed $M>0$, the contribution from $|k-m_n| \le M c_n$ (term $I_1$ therein), divided by $\cC_n^{\cD} $, converges to $\int_{-M}^{M}f(-x) \diff x$, while the remaining contribution is at most $CM^{-1+\eta}$ for $0 < \eta < 1$. Letting first $n \rightarrow \infty$ and then $M \rightarrow \infty$, and using that $f$ is a probability density, gives
	\[
		I_n^{\cD} \sim \frac{v_{m_n} W(\rho_w)^{m_n}}{\mu(a_{m_n})} = \cC_n^{\cD} .
	\]
	Choosing $t_n$ as in~\eqref{eq:boundaryrn}, the estimates~\eqref{eq:boundarysmallrange} (term $I_1$ therein) and \eqref{eq:boundaryintermediate} (term $I_2$ therein), whose proofs do not assume the single-big-jump dominance assumption~\eqref{eq:geldominance}, give
	\[
		\sum_{k \le t_n}v_k W(\rho_w)^k \Prb{S_k = n} \sim \cC_n^{\cG} , \qquad \sum_{t_n<k < m_n/2}v_k W(\rho_w)^k \Prb{S_k = n} = o(\cC_n^{\cG} ).
	\]
	Hence, $I_n^{\cG}  \sim \cC_n^{\cG} $. The first assertion of the Theorem then follows immediately, while
	\[
		\Prb{\cD_n} = \sum_{k \ge m_n/2} \Prb{N_n = k} = \frac{\sum_{k \ge m_n/2} v_k W(\rho_w)^k \Prb{S_k = n}}{u_n \rho_w^n} = \frac{I_n^{\cD}}{I_n^{\cD} + I_n^{\cG}}
	\]
	and similar calculations for $\cG_n$ prove the claimed asymptotics for $\Prb{\cD_n}$, respectively $\Prb{\cG_n}.$

	We have proved the partition-function asymptotic and the asymptotics for $\Prb{\cD_n}$ and $\Prb{\cG_n}$. We now show the limits conditional on $\cD_n$. For $k \ge m_n/2$, the conditional distribution of the number of components is
	\[
		\Prb{N_n = k \mid \cD_n} = \frac{v_k W(\rho_w)^k \Prb{S_k = n}}{I_n^{\cD}}.
	\]
	Since $c_n \mu(a_{m_n}) = a_{m_n}$, the identity above, $I_n^{\cD} \sim \cC_n^{\cD} $, the uniform convergence theorem~\ref{thm:uct} and \eqref{eq:outer-weight-tail} give, for $k \ge m_n/2$,
	\[
		c_n \Prb{N_n = k \mid \cD_n} \le C a_{m_n}\Prb{S_k = n}.
	\]
	For $0 < \eta <1$, the large-deviation bound~\eqref{eq:bedev}, the Potter bounds~\eqref{eq:pb} and~\eqref{eq:defan2} control the range $m_n/2 \le k \le 2m_n$. Lemma~\ref{lem:large-count-exponential} controls $k > 2m_n$. Consequently, for some constants $C_\eta,M_0>0$ and all sufficiently large $n$,
	\[
		c_n \Prb{N_n = k \mid \cD_n} \le
		\begin{cases}
			C_\eta\left(\dfrac{|k-m_n|}{c_n}\right)^{-2+\eta},
			& \begin{gathered}
				m_n/2 \le k \le 2m_n,\\
				|k-m_n| > M_0c_n,
			\end{gathered}\\[2ex]
			C\exp(-\kappa k/c_n), & k > 2m_n.
		\end{cases}
	\]
	On each fixed window $|k-m_n| \le Mc_n$, \eqref{eq:almostllt} and $I_n^{\cD} \sim \cC_n^{\cD} $ give the uniform local limit. The above tail bounds and~\eqref{eq:ftail} extend this convergence to all $k \ge m_n/2$, by first letting $n \rightarrow \infty$ and then $M \rightarrow \infty$. It follows that
	\begin{align}
		\label{eq:km2llt}
		\sup_{k \ge m_n/2} \left|c_n \Prb{N_n = k \mid \cD_n} -f\left(-\frac{k-m_n}{c_n}\right) \right|\longrightarrow 0.
	\end{align}
	Uniformly for $k<m_n/2$ we have $\Prb{N_n = k \mid \cD_n} =0$ and $\left| \frac{k-m_n}{c_n} \right| \ge \left| \frac{m_n}{2c_n} \right|\to \infty$, so that $f\left(-\frac{k-m_n}{c_n}\right)\to 0$. Hence~\eqref{eq:km2llt} also holds when $k$ ranges over all integers.
	
	A Riemann summation argument for the conditional local limit gives, for every fixed $x < y$,
	\[
		\Prb{x<\frac{N_n-m_n}{c_n} \le y \middle|\, \cD_n} \longrightarrow \int_x^y f(-t)\,\diff t.
	\]
	Since $f(-\cdot)$ is a probability density, this implies
	\[
		\left(\frac{N_n-m_n}{c_n}\,\middle|\,\cD_n \right) \convdis -Y.
	\]
    
	To establish the conditional version of Lemma~\ref{lem:tv}, fix any sequence of positive integers satisfying $c_n = o(s_n)$ and $s_n = o(m_n)$, and take
	\[
		r_n := m_n-s_n, \qquad h_n := \lfloor\sqrt{c_n s_n}\rfloor, \qquad E_n := \{|N_n-m_n| \le h_n \}.
	\]
	Since $h_n = o(s_n)$ and $s_n = o(m_n)$, we have, for $n$ large enough, $E_n \subseteq \cD_n \cap \{N_n \ge r_n\}$. The conditional limit for the number of components implies that $(N_n-m_n)/c_n$ is bounded in probability given $\cD_n$. Since $h_n/c_n \rightarrow \infty$, it follows that $\Prb{E_n \mid \cD_n} \rightarrow 1$. Furthermore, $I_n^{\cD} \sim v_{m_n} W(\rho_w)^{m_n}/\mu(a_{m_n})$. Therefore, the calculation in~\eqref{eq:likrat} applies given $\cD_n$, with the summation restricted by $|k-m_n| \le h_n$. For every fixed $M > 0$, uniformly over vectors $\mathbf{x}$ of positive probability which satisfy $|\sum_i x_i-b_{r_n}| \le M a_{r_n}$, this yields
	\[
		\frac{\Prb{(K_1,\ldots,K_{r_n}) = \mathbf{x},E_n \mid \cD_n}}{\Prb{(X_1, \ldots, X_{r_n}) = \mathbf{x}}}
		\sim \mu(a_{m_n}) \sum_{j: |j-s_n| \le h_n} \Prb{S_j= n-\sum_i x_i} \longrightarrow 1,
	\]
	where the change of index is $j = k-r_n$ and the final convergence follows from~\eqref{eq:renewaluniform}. The estimates used to prove~\eqref{eq:renewaluniform} depend only on the inner law~\eqref{eq:wasymptotic} and the assumptions on $s_n$, so they apply unchanged here. This is the equivalent to the ratio~\eqref{eq:likrat2}, now under a conditional law. The final argument of Lemma~\ref{lem:tv} gives in this case
	\[
		d_{\mathrm{TV}} \left(\big((K_1, \ldots, K_{\min(N_n,r_n)})\mid \cD_n \big), (X_1, \ldots, X_{r_n})\right) \longrightarrow 0.
	\]
	Because $\cD_n$ is determined by $N_n$, conditioning on $\cD_n$ preserves exchangeability of the component sizes. Therefore, the proofs of Theorem~\ref{thm:functional} and Corollary~\ref{cor:largest} apply with the conditional total-variation approximation and the conditional limits for the number of components proved above. This establishes the claims given $\cD_n$.

	Now, we prove the limits conditional on $\cG_n$. Choose $t_n \rightarrow \infty$ with $t_n = o(m_n)$ and $H(t_n)/H(m_n) \rightarrow 1$. Equation~\eqref{eq:uniformbigjumplocal} and $I_n^{\cG} \sim \cC_n^{\cG} $ give, uniformly for $k \le t_n$,
	\[
		\Prb{N_n = k \mid \cG_n} = \frac{v_k W(\rho_w)^k \Prb{S_k = n}}{I_n^{\cG} } = (1+o(1))\frac{k v_k W(\rho_w)^k}{H(m_n)}.
	\]
	Since $H(t_n)/H(m_n) \rightarrow 1$, summing over $1 \le k \le t_n$, and using $\widehat{N}_n$ in~\eqref{eq:defNhatboundary}, yields
	\[
		d_{\mathrm{TV}}\big((N_n \mid \cG_n),\widehat{N}_n \big) \longrightarrow 0, \qquad \Prb{N_n > t_n \mid \cG_n} \longrightarrow 0.
	\]
	Let $I_n$ be the position of the first maximal component. The bound~\eqref{eq:boundarymarkedfixedk} holds uniformly for $1 \le k \le t_n$ with $v_k > 0$, and its proof does not use the condition~\eqref{eq:geldominance}. For these values of $k$, and for $n$ large enough, the event $\{N_n = k\}$ is contained in $\cG_n$. If we condition additionally on $\cG_n$, this leaves the law given $N_n = k$ unchanged. Summing the bounds with weights $\Prb{N_n = k \mid \cG_n}$, over values of $k$ with positive conditional probability, gives
	\begin{multline}\label{eq:b2gelmarkedconditional}
		d_{\mathrm{TV}}\Big(\big((N_n,\tau(K_1,\ldots,K_{N_n}),I_n)\mid \cG_n \big),\big(\widehat{N}_n,(X_1,\ldots,X_{\widehat{N}_n-1}),J_n \big)\Big)\\
		\le d_{\mathrm{TV}}\big((N_n \mid \cG_n),\widehat{N}_n \big) +\Prb{N_n > t_n \mid \cG_n} + o(1) \longrightarrow 0.
	\end{multline}
	Applying the insertion map to~\eqref{eq:b2gelmarkedconditional} gives~\eqref{eq:b2gelconditional}. Using then the map $(k,\mathbf{x},j) \mapsto (k,\sum_i x_i,\mathbf{x}^{\downarrow})$ gives the conditional version of~\eqref{eq:boundaryjointtv}. The proofs of Corollary~\ref{pro:boundarycondensation}, Theorem~\ref{thm:divergent condensationtv}, and Corollaries~\ref{cor:meso-ranked},~\ref{cor:meso-functionals} and~\ref{cor:meso-scale} therefore apply unchanged under the conditional law. The proof of Proposition~\ref{pro:meso-functional} also applies given $\cG_n$: in its final total-variation comparison, replace~\eqref{eq:boundarymarkedtv} by~\eqref{eq:b2gelmarkedconditional}. This approximation preserves the order of the residual components, as required in order to construct $Z_n$ and $\Pi_n$. 

	If $\Theta_n \rightarrow \Theta \in [0,\infty]$, set $p = 1/(1+\Theta)$, with $p = 0$ when $\Theta = \infty$. Then $\Prb{\cD_n} \rightarrow p$ and $\Prb{\cG_n} \rightarrow 1-p$. Together with the conditional conclusions proved above, this gives the claimed mixture.
\end{proof}

The ratio $\Theta_n$ determines the probabilities of the events $\cD_n$ and $\cG_n$. Conditionally on $\cG_n$, the number of components converges to the finite random variable $\widehat{N}$ when $H(\infty) < \infty$, and tends to infinity in probability when $H(\infty) = \infty$. In particular, when $H(\infty) < \infty$,
\[
	u_n \sim V'(W(\rho_w)) w_n \quad \Longleftrightarrow \quad \Theta_n \rightarrow \infty \quad \Longleftrightarrow \quad d_{\mathrm{TV}}\big((K_1, \ldots, K_{N_n}),\boldsymbol{\Xi}_n \big) \rightarrow 0.
\]
Since $H(m_n) \rightarrow H(\infty) = W(\rho_w) V'(W(\rho_w))$, the coefficient asymptotic gives
\[
	\frac{u_n}{V'(W(\rho_w)) w_n} \sim 1+\Theta_n^{-1},
\]
which proves the first equivalence. For the second, convergence in total variation of the complete vectors implies that $N_n$ is bounded in probability, since the lengths of the comparison vectors are bounded in probability. As $m_n \rightarrow \infty$, this gives $\Prb{\cD_n} \rightarrow 0$, and hence $\Theta_n \rightarrow \infty$ by Theorem~\ref{thm:b2}. Conversely, if $\Theta_n \rightarrow \infty$, then $\Prb{\cD_n} \rightarrow 0$. The approximation in that latter theorem, this time conditional on $\cG_n$, gives, when $H(\infty) < \infty$,
\[
	d_{\mathrm{TV}}\left((K_1,\ldots,K_{N_n}),\boldsymbol{\Xi}_n \right) \le \Prb{\cD_n} +d_{\mathrm{TV}}\left(((K_1,\ldots,K_{N_n})\mid \cG_n),\boldsymbol{\Xi}_n \right) \longrightarrow 0.
\]

The residual-mass fluctuation limit satisfies a local refinement. Appealing to Theorem~\ref{thm:local-big-jump} through the Kolchin representation~\eqref{eq:KasX} gives the following local limit conditional on the number of components. 

\begin{corollary}[Conditional local limit in divergent condensation]
	\label{cor:meso-conditional-llt}
	Under the assumptions of Theorem~\ref{thm:divergent condensationtv},
	\begin{equation}\label{eq:meso-conditional-llt}
		\sup_{j \in \ndZ} \left|a_{N_n}\Prb{n-K_{(1)} = j \mid N_n} - f\left(\frac{j-b_{N_n}}{a_{N_n}}\right)\right| \convp 0.
	\end{equation}
	Under the assumptions of Theorem~\ref{thm:b2} with $H(\infty) = \infty$, the same statement holds given $\cG_n$.
\end{corollary}

\begin{proof}
	We first prove the unconditional assertion. By a diagonal choice using~\eqref{eq:mesonumbergrowth}, there exist deterministic integer sequences $r_n \rightarrow \infty$ and $t_n = o(m_n)$ such that
	\[
		\Prb{r_n \le N_n \le t_n} \longrightarrow 1.
	\]
	Since eventually $t_n \le m_n$, monotonicity of $(a_{\cdot})$ and $\mu$, which implies $a_{t_n} \le a_{m_n}$ and $\mu(a_{t_n}) \le \mu(a_{m_n})$, together with~\eqref{eq:bdiff2co}, gives
	\[
		\frac{b_{t_n}}{n} = \frac{t_n \mu(a_{t_n})}{n} \le \frac{t_n \mu(a_{m_n})}{n} = \frac{t_n}{m_n}\frac{b_{m_n}}{n} \longrightarrow 0.
	\]
	Recalling that~\eqref{eq:amo} gives $a_{m_n} = o(n)$, and combining with the above bound, for all
	sufficiently large $n$,
	\[
		\inf_{r_n \le k \le t_n}\frac{n-b_k}{a_k} \ge \frac{n-b_{t_n}}{a_{m_n}} \ge \frac{n}{2a_{m_n}}\longrightarrow \infty.
	\]
	Writing once again $M_k := \max_{1 \le i \le k} X_i$, Theorem~\ref{thm:local-big-jump} consequently gives
	\begin{equation}\label{eq:meso-conditional-uniform}
		\sup_{r_n \le k \le t_n} \sup_{j \in \ndN_0} \left|a_k \Prb{M_k= j \mid S_k = n} - f\left(\frac{n-b_k-j}{a_k}\right) \right| \longrightarrow 0.
	\end{equation}
	To justify uniformity, suppose that the above fails. Then, there exist $n_i \rightarrow \infty$ and $k_i \in [r_{n_i},t_{n_i}]$ for which the inner supremum stays bounded away from zero. However, $k_i \rightarrow \infty$ and $(n_i-b_{k_i})/a_{k_i} \rightarrow \infty$. Extract a strictly increasing subsequence of $(k_i)$. But then, Theorem~\ref{thm:local-big-jump} applies along these indices, so this gives a contradiction. For all sufficiently large $n$, $\Prb{N_n = k} > 0$ for every integer $k \in [r_n, t_n]$, by~\eqref{eq:vboundary} and~\eqref{eq:uniformbigjumplocal}.

    For these values of $k$ and every integer $j \le n$, the Kolchin representation~\eqref{eq:KasX} gives
    \[
        \Prb{n-K_{(1)}=j \mid N_n = k} = \Prb{M_k=n-j \mid S_k =n}.
    \]
	Substituting $n-j$ for the maximum in~\eqref{eq:meso-conditional-uniform} gives the approximation uniformly for $j \le n$ and $r_n \le k \le t_n$, since
	\[
		\frac{n-b_k-(n-j)}{a_k} = \frac{j-b_k}{a_k}.
	\]
	On the other hand, when $j > n$, we have $\Prb{n-K_{(1)} = j \mid N_n = k} = 0$, since $n-K_{(1)} \le n$. Furthermore,
	\[
		\inf_{\substack{r_n \le k \le t_n \\ j > n}} \frac{j-b_k}{a_k} \ge \frac{n-b_{t_n}}{a_{m_n}} \longrightarrow \infty.
	\]
    Since $f(x) \rightarrow 0$ as $|x| \rightarrow \infty$, we obtain
	\[
		\sup_{r_n \le k \le t_n} \sup_{j \in \Z} \left\vert a_k \Prb{n-K_{(1)} = j \mid N_n = k} - f\left(\frac{j-b_k}{a_k}\right) \right\vert \longrightarrow 0.
	\]

	Evaluating the above estimate at $k = N_n$ on $\{r_n \le N_n \le t_n \}$ (which has probability tending to one) proves~\eqref{eq:meso-conditional-llt}.

	Under the assumptions of Theorem~\ref{thm:b2} with $H(\infty) = \infty$, the conditional limits for the number of components allow the same choice of sequences, which now satisfy
	\[
		\Prb{r_n \le N_n \le t_n \mid \cG_n} \longrightarrow 1.
	\]
	Because $\cG_n$ is determined by $N_n$, for every $k$ with $\Prb{N_n = k \mid \cG_n} > 0$,
	\[
		\Prb{n-K_{(1)} = j \mid N_n = k, \cG_n} = \Prb{n-K_{(1)} = j \mid N_n = k}.
	\]
	On $\{r_n \le N_n \le t_n\}$, the error in~\eqref{eq:meso-conditional-llt} is bounded by the deterministic supremum above, which tends to zero. The complementary event has conditional probability tending to zero given $\cG_n$. This proves the second assertion.
\end{proof}

We proceed with an example in which both slowly varying factors are logarithmic, in order to illustrate the competition. 

\begin{example}[Logarithmic weights and phase transitions]
	\label{ex:logs}
	Suppose that~\eqref{eq:wasymptotic},~\eqref{eq:critical} and~\eqref{eq:vboundary} hold, with constants $c_w,c_v>0$ and exponents $r,s\in \ndR$ such that
	\[
		L_w(n)\sim c_w\log^r n, \qquad L_v(n)\sim c_v\log^s n.
	\]
	Set $c:=c_w/W(\rho_w)$, so that $L(n)\sim c \log^r n$. For $r < -1$, write $\mu:= \E{X}<\infty$. Also, $H(\infty) < \infty$ if and only if $s<-1$, in which case
	\[
		H(\infty) =W(\rho_w)V'(W(\rho_w)).
	\]
	Example~\ref{ex:addme} gives $\log m_n \sim \log n$ and
	\[
		\mu(a_{m_n}) \sim
		\begin{cases}
			\displaystyle \frac{c}{r+1}\log^{r+1}n, & r>-1,\\
			c \log \log n, & r=-1,\\
			\mu, & r < -1,
		\end{cases}
	\]
	where $\mu:= \E{X}<\infty$ in the last case. Similarly, the definition of $H$ yields
	\[
		H(m_n) \sim
		\begin{cases}
			\displaystyle \frac{c_v}{s+1}\log^{s+1}n, & s>-1,\\
			c_v\log \log n, & s=-1,\\
			H(\infty), & s<-1,
		\end{cases}
	\]
	where $H(\infty) < \infty$ in the last case. Consequently,
	\begin{align*}
		\frac{L_v(m_n)\mu(a_{m_n})}{H(m_n)L(n)} \sim \frac{c_v}{c}\frac{\log^{s-r}n \,\mu(a_{m_n})}{H(m_n)}.
	\end{align*}
	% \begin{figure}[t]
	% 	\centering
	% 	\begin{minipage}{\textwidth}
	% 		\centering
	% 		\includegraphics[width=0.7\linewidth]{regimes2.pdf}
	% 		\caption[]{Asymptotic regimes for log-power slowly varying factors.}
	% 		\label{fi:regimes}
	% 	\end{minipage}
	% \end{figure}

    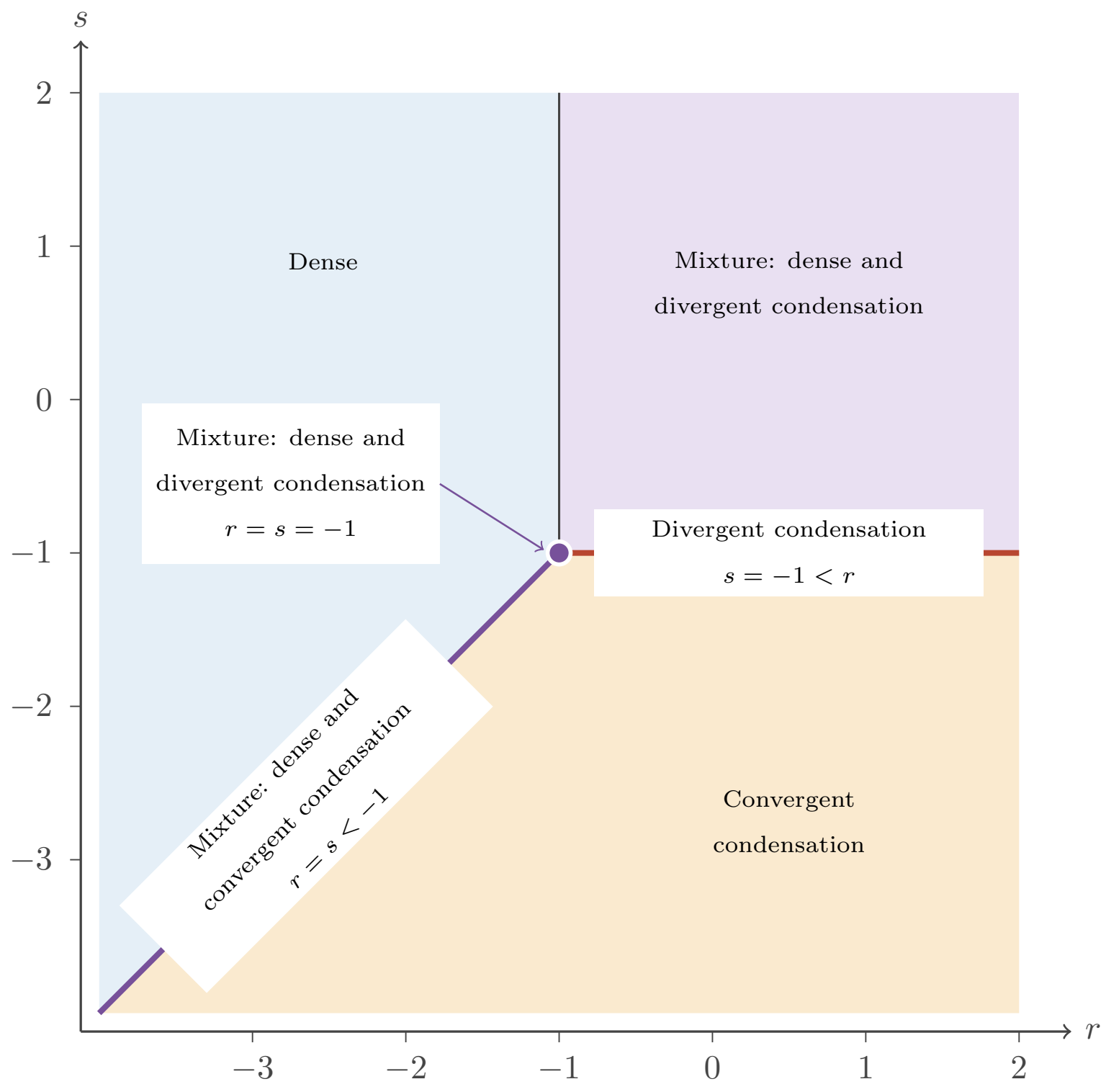
\begin{figure}[t]
    \centering
    \resizebox{0.75\linewidth}{!}{
        \input{regimes_tikz.tex}
    }
    \caption[]{Asymptotic regimes for log-power slowly varying factors.}
    \label{fig:regimes}
    \end{figure}
    
	Substituting the estimates for $m_n$ and $\mu(a_{m_n})$ gives
	\[
		\cC_n^{\cD}  \sim \begin{cases}\dfrac{c c_v}{r+1}n^{-2}\log^{r+s+1}n, & r>-1,\\[1ex]
			c c_v n^{-2}\log^s n \log \log n, & r=-1,\\
			\mu c_v n^{-2}\log^s n, & r < -1.\end{cases}
	\]
	Similarly, summing the outer weights in $H(m_n)$ gives
	\[
		\cC_n^{\cG}  \sim \begin{cases}\dfrac{c c_v}{s+1}n^{-2}\log^{r+s+1}n, & s>-1,\\[1ex]
			c c_v n^{-2}\log^r n \log \log n, & s=-1,\\
			cH(\infty)n^{-2}\log^r n, & s<-1.\end{cases}
	\]
	By Theorem~\ref{thm:b2}, these equivalents determine the tilted partition function asymptotics as
	\[
		u_n \rho_w^n \sim \cC_n^{\cD} +\cC_n^{\cG} .
	\]
	Taking the ratio of the two contributions and using $\Prb{\cG_n}\sim\Theta_n/(1+\Theta_n)$ gives
	\[
		\Prb{\cG_n}\longrightarrow \begin{cases}\dfrac{r+1}{r+s+2}, & r>-1,\ s>-1,\\[1ex]
			\dfrac{1}{2}, & r= s=-1,\\[1ex]
			\dfrac{cH(\infty)}{cH(\infty)+c_v\mu}, & r= s<-1,\\[1ex]
			1, & s<\min \{r,-1\}\text{ or }s=-1<r,\\
			0, & r\le-1,\ s>r.\end{cases}
	\]
	In particular, an asymptotic mixture behaviour occurs throughout $r,s>-1$ and on the diagonal $r= s\le-1$.

	The conditional conclusions of Theorem~\ref{thm:b2} apply even when $\Prb{\cG_n} \rightarrow 0$. Given $\cG_n$, condensation is convergent if $s < -1$ and divergent if $s \ge -1$, according to whether $H(\infty)$ is finite or infinite. See Figure~\ref{fig:regimes}.
    
    The calculations to obtain the above quantities are straightforward and therefore omitted. 
\end{example}

\section{Applications to conditioned BGW trees and forests}
\label{sec:bgw}

We apply the results of the earlier sections to three related questions: the maximal outdegree of a conditioned Bienaymé--Galton--Watson (BGW) tree, the sizes of the subtrees attached to its root, and the component trees of a conditioned forest.

Throughout this section, let $\xi$ be a nonnegative integer-valued random variable with mean
\begin{align}\label{eq:bgw1}
	\mathfrak{m} := \E{\xi} \in ]0,1]
\end{align}
and probabilities
\begin{align}\label{eq:bgw2}
	p_k := \Prb{\xi = k} \sim L_\xi(k)k^{-2}, \qquad k \rightarrow \infty,
\end{align}
where $L_\xi$ is slowly varying. Let $\cT$ be a BGW tree with offspring law $\xi$, and let $\cT_n$ have the law of $\cT$ conditioned on $|\cT| = n$. The tail assumption gives $p_k > 0$ for all sufficiently large $k$, so the conditioning is well defined for all sufficiently large $n$.

\subsection{The maximal degree of a conditioned BGW tree}

Write $\Delta_{\mathrm{max}}(\cT_n)$ for the largest outdegree of $\cT_n$. Construct the scale $\alpha_n$ for $\xi$ in the same way that $a_n$ was constructed for $X$ in Section~\ref{sec:background}, so that
\[
	n \Prb{\xi > \alpha_n} \longrightarrow 1, \qquad \alpha_n \sim nL_\xi(\alpha_n),
\]
and set
\[
	B_n := n \E{\xi;\xi\le\alpha_n}, \qquad d_n := n-1-B_n.
\]
This way, $B_n$ is the counterpart of the centering $b_n$ for the offspring distribution. We write $f$ for the density of the variable $Y$ in~\eqref{eq:chara}. The following theorem gives a local limit counterpart of the distributional results of Kortchemski and Richier~\cite[Theorem~1]{zbMATH07057468} and Kortchemski and Vetter~\cite[Corollary~4.1]{zbMATH08114289}, with the centering and Cauchy normalisation fixed above. Note that, for subcritical BGW trees with offspring tail index $\alpha > 1$, a uniform local limit for the maximal degree is given in \cite[Lemma~2.2]{MR4144886}. 

\begin{theorem}[Local limit theorem for the maximal degree]\label{thm:bgwmax}
	Suppose that~\eqref{eq:bgw1} and~\eqref{eq:bgw2} hold. Then
	\begin{equation}\label{eq:bgwmax}
		\sup_{k \in \ndZ} \left|\alpha_n \Prb{\Delta_{\mathrm{max}}(\cT_n) =k}-f\left(\frac{d_n-k}{\alpha_n}\right)\right| \longrightarrow 0.
	\end{equation}
\end{theorem}
\begin{proof}
	Let $\xi_1,\xi_2,\ldots$ be independent copies of $\xi$. The tree $\cT_n$ admits a well-known coupling with a balls-in-boxes model, so that the outdegrees of $\cT_n$ in depth-first-search order correspond to a cyclic shift of the conditioned vector $\left( (\xi_1, \ldots, \xi_n) \mid \sum_{i=1}^n (\xi_i -1) = -1 \right)$. See~\cite[Lem. 17.1]{MR2908619}. Hence
	\[
		\Delta_{\mathrm{max}}(\cT_n)\eqdist\left(\max_{1 \le i\le n}\xi_i\,\middle|\,\sum_{i=1}^n \xi_i= n-1\right).
	\]
	We are going to argue that
	\begin{align}\label{eq:showme}
		d_n/\alpha_n \rightarrow \infty.
	\end{align}
	If $\mathfrak{m} < 1$, then $d_n \sim (1-\mathfrak{m})n$, whereas $\alpha_n = o(n)$ by~\eqref{eq:finitemean}. This verifies~\eqref{eq:showme} in this case. If $\mathfrak{m} = 1$, define
	\[
		\ell_\xi^*(x):= \int_x^\infty\frac{L_\xi(t)}{t}\,\diff t.
	\]
	In this case, $d_n= n \E{\xi;\xi>\alpha_n}-1$. Equations~\eqref{eq:ellstarasymp} and~\eqref{eq:ellstarl}, together with $\alpha_n \sim nL_\xi(\alpha_n)$, yield
	\begin{equation}\label{eq:bgwscale}
		d_n \sim n \ell_\xi^*(\alpha_n), \qquad \frac{d_n}{\alpha_n} \sim \frac{\ell_\xi^*(\alpha_n)}{L_\xi(\alpha_n)} \longrightarrow \infty.
	\end{equation}
	Hence~\eqref{eq:showme} holds in both cases. This way, we can rely on Theorem~\ref{thm:local-big-jump}, yielding
	\[
		\sup_{j \in \ndZ} \left|\alpha_n \Prb{\Delta_{\mathrm{max}}(\cT_n) = j}-f\left(\frac{d_n-j}{\alpha_n}\right)\right|\longrightarrow 0.
	\]
\end{proof}

\subsection{The partition at the root}

Set $q_j := \Prb{|\cT| = j}$ for $j \ge 1$, and $q_0 := 0$. For $n \ge 2$, let $D_n$ be the root degree of $\cT_n$, and write $P_n = (K_1,\ldots,K_{D_n})$ for the ordered sizes of the subtrees attached to the root. These sizes sum to $n-1$. The vector $P_n$ is the Gibbs partition of $n-1$ with weights $v_k = p_k$ and $w_j = q_j$.

We reuse our framework to this specific Gibbs partition. The following corollary complements the condensation results of Kortchemski and Richier~\cite{zbMATH07057468} and Kortchemski and Vetter~\cite{zbMATH08114289} with a conditional local limit for the root degree. Note that, for every fixed $k\ge 1$, the root degree satisfies $\Prb{D_n=k}\to kp_k$ by~\cite[Theorem~7.10]{MR2908619}. The limiting law on $\ndN_0\cup\{\infty\}$ has mass $1-\mathfrak{m}$ at infinity, see also~\cite[Section~20]{MR2908619} for the relation between the root degree and the maximal degree
in the subcritical case. 

\begin{corollary}[Partition at the root]\label{cor:bgwroot}
	Under the assumptions of Theorem~\ref{thm:bgwmax},
	\begin{equation}\label{eq:bgwsize}
		\Prb{|\cT| = n} \sim \Prb{\xi = \lfloor d_n \rfloor},
	\end{equation}
	and $(\Prb{|\cT| = n})$ is regularly varying with index $-2$. 

	Let $m_n$ and $c_n$ be the centering and fluctuation sequences defined in~\eqref{eq:defmn} and~\eqref{eq:defcn}, respectively, constructed for the inner distribution $(q_j)$. Define
	\[
		\cG_n := \{D_n<m_{n-1}/2\}, \qquad \cD_n := \{D_n \ge m_{n-1}/2\}.
	\]
	Then
	\begin{equation}\label{eq:rootprob}
		\Prb{\cG_n}\longrightarrow \mathfrak{m}, \qquad \Prb{\cD_n}\longrightarrow 1\longrightarrow 1-\mathfrak{m}.
	\end{equation}

	Conditionally on $\cD_n$, the root degree satisfies the local limit theorem
	\begin{equation}\label{eq:rootllt}
		\sup_{k \in \ndZ} \left|c_{n-1}\Prb{D_n=k\mid\cD_n}-f\left(-\frac{k-m_{n-1}}{c_{n-1}}\right)\right| \longrightarrow 0,
	\end{equation}
	where $f$ is the density of the random variable $Y$ characterised in~\eqref{eq:chara}. This conditional limit also holds when $\mathfrak{m} = 1$, although $\Prb{\cD_n} \rightarrow 0$ in this case.

	Conditionally on $\cG_n$, the approximation~\eqref{eq:boundaryfinitejoint} holds with total mass $n-1$, inner distribution $(q_j)$ and limiting number of components $\widehat{\xi}$ given by
	\[
		\Prb{\widehat{\xi} =k} = \frac{kp_k}{\mathfrak{m}}, \qquad k \ge 1.
	\]
	Explicitly, the approximating vector is obtained by taking $\widehat{\xi}-1$ independent sizes with distribution $(q_j)$ and inserting the remaining mass ($n-1$ minus the sum of the sampled sizes) at an independent uniform position. On the event that the remaining mass is negative, use the convention in~\eqref{eq:Xi}.
\end{corollary}
\begin{proof}
	By Dwass' identity~\cite{MR0253433}, recalled in~\cite[Theorem~15.5]{MR2908619}, this gives
	\[
		q_n = \frac{1}{n}\Prb{\xi_1+\cdots+\xi_n= n-1}.
	\]
	By the proof of Theorem~\ref{thm:bgwmax}, the mass value $n-1$ exceeds the centering $B_n$ by $n-1-B_n = d_n$, and $d_n/\alpha_n \rightarrow \infty$. The local large-deviation estimate~\eqref{eq:beloc} therefore proves~\eqref{eq:bgwsize}. In the critical case $\mathfrak{m}=1$, this recovers the local total-progeny estimate of Kortchemski and Richier~\cite[Lemma~17(i)]{zbMATH07057468}.

	The sequence $(d_n)$ is regularly varying with index $1$: this follows from $d_n \sim(1-\mathfrak{m})n$ in the subcritical case and from~\eqref{eq:bgwscale} in the critical case. Since $(p_k)$ varies regularly with index $-2$, the uniform convergence theorem~\ref{thm:uct} gives, for every $t > 0$,
	\[
		\frac{q_{\lfloor tn \rfloor}}{q_n} \sim \frac{p_{\lfloor d_{\lfloor tn \rfloor}\rfloor}}{p_{\lfloor d_n \rfloor}}\longrightarrow t^{-2}.
	\]
	This proves the claimed regular variation of $(q_n)$.

	The branching property of BGW trees gives, for $k \ge 1$, $n_i \ge 1$ and $n_1+\cdots+n_k= n-1$,
	\[
		\Prb{P_n=(n_1,\ldots,n_k)} = \frac{p_k}{q_n}\prod_{i=1}^kq_{n_i}.
	\]
	This is the Gibbs law with $v_k = p_k$, $w_j = q_j$ and partition function $u_{n-1} = q_n$. The generating series of these weight sequences have radius of convergence $1$ and satisfy
	\[
		W(1) = \sum_{j \ge 1} q_j = 1, \qquad H(\infty) = V'(1) = \sum_{k \ge 1}kp_k= \mathfrak{m}.
	\]
	The assumptions of Theorem~\ref{thm:b2} are therefore satisfied.

	At total mass $n-1$, the contribution from $D_n < m_{n-1}/2$ is asymptotically equivalent to
	\[
		q_{n-1}H(m_{n-1}) \sim \mathfrak{m}q_n.
	\]
	Dividing by the partition function $q_n$ gives $\Prb{\cG_n} \rightarrow \mathfrak{m}$ and $\Prb{\cD_n} \rightarrow 1-\mathfrak{m}$, proving~\eqref{eq:rootprob}. The conditional local limit theorem~\eqref{eq:rootllt} follows from the assertion of Theorem~\ref{thm:b2} conditional on~$\cD_n$.
\end{proof}

\subsection{BGW forests}

We provide an example of divergent condensation at the root of a BGW forest. Suppose that the offspring law $\xi$ satisfies $\Ex{\xi} = \mathfrak{m} = 1$ and
\[
	p_k\sim\frac{c}{k^2(\log k)^2}, \qquad c>0.
\]
Equation~\eqref{eq:example-an} and direct integration give
\[
	\alpha_n \sim\frac{cn}{(\log n)^2}, \qquad \ell_\xi^*(x)\sim\frac{c}{\log x}.
\]
Since $\log\alpha_n \sim \log n$, \eqref{eq:bgwscale} gives $d_n\sim cn/\log n$, as in~\cite[Example~19]{zbMATH07057468}. Hence $\log d_n \sim \log n$, and the local total-progeny estimate~\eqref{eq:bgwsize} gives
\begin{equation}\label{eq:foresttail}
	q_n \sim p_{\lfloor d_n \rfloor} \sim \frac{c}{d_n^2(\log d_n)^2} \sim \frac{1}{cn^2}.
\end{equation}
Let $\cT^{(1)},\cT^{(2)},\ldots$ be independent copies of $\cT$. Let $R$ denote an independent random positive integer such that for some $c_R>0$ and $r \in \ndR$
\[
	\Prb{R=k} \sim c_Rk^{-2}(\log k)^r, \qquad k \rightarrow \infty.
\]
Condition the ordered forest $(\cT^{(1)},\ldots,\cT^{(R)})$ to have total mass $n$, and denote the resulting forest by $\cF_n$. This way, the sizes of the trees in $\cF_n$ is a Gibbs partition of $n$ with weights
\[
	v_k = \Prb{R=k}, \qquad w_j = \Prb{\lvert \cT \rvert =j}.
\]
Write $N_n$ for its number of trees. From now on, the sequences $a_k,b_k,m_n,c_n$ refer to the inner law $(q_j)$, and as before let
\[
	\cG_n := \{N_n<m_n/2\}, \qquad \cD_n := \{N_n \ge m_n/2\}.
\]
By~\eqref{eq:foresttail}, Example~\ref{ex:logs} applies with inner logarithmic exponent $0$ and outer exponent $r$. It gives
\[
    \Prb{\cG_n}\longrightarrow
    \begin{cases}
        \dfrac{1}{r+2}, & r>-1,\\[1ex]
        1, & r\le -1.
    \end{cases}
\]
For $r < -1$, $H(\infty) < \infty$ and the Gibbs partition is in the convergent condensation phase. For $r \ge -1$, we have $H(\infty) = \E{R} = \infty$. At $r = -1$, condensation is divergent and has probability tending to one. For $r > -1$, the events $\{N_n \ge m_n/2\}$ and $\{N_n < m_n/2\}$ have positive limiting probabilities, corresponding to dense and divergent condensation. Theorem~\ref{thm:b2} describes the asymptotic behaviour of component number and sizes conditionally on each of these events.

\section{Acknowledgement}

This research was funded in part by the Austrian Science Fund (FWF) 10.55776/PAT6732623 and 10.55776/PAT1382025. For open access purposes, the authors have applied a CC BY public copyright license to any author-accepted manuscript version arising from this submission. 

\bibliographystyle{abbrv}
\bibliography{gibo}

\end{document}

%% file: regimes_tikz.tex
    \begingroup
    \definecolor{regimeDense}{HTML}{E5EFF7}
    \definecolor{regimeConvergent}{HTML}{FAEACF}
    \definecolor{regimeMixture}{HTML}{E9E0F2}
    \definecolor{regimePurple}{HTML}{77519A}
    \definecolor{regimeRed}{HTML}{B84630}
    \definecolor{regimeAxis}{HTML}{444444}
    \begin{tikzpicture}[
      x=1.6cm,
      y=1.6cm,
      >=to,
      line cap=butt,
      line join=miter,
      regime label/.style={font=\fontsize{7}{11.25}\selectfont, align=center, inner sep=0pt, text=black},
      regime box/.style={regime label, fill=white, outer sep=0pt},
      regime axis/.style={draw=regimeAxis, line width=0.7bp},
      regime tick/.style={draw=regimeAxis, line width=0.4bp},
      regime tick label/.style={font=\fontsize{10}{12}\selectfont, text=regimeAxis, inner sep=0pt}
    ]
      % Keep the white margins and aspect ratio close to the supplied figure.
      \path[use as bounding box] (-4.93,-4.82) rectangle (2.87,2.86);
    
      % Three open regions, displayed in the window [-4,2] x [-4,2].
      \fill[regimeDense] (-4,-4) -- (-1,-1) -- (-1,2) -- (-4,2) -- cycle;
      \fill[regimeConvergent] (-4,-4) -- (2,-4) -- (2,-1) -- (-1,-1) -- cycle;
      \fill[regimeMixture] (-1,-1) rectangle (2,2);
    
      % Boundaries r=-1, r=s<-1 and s=-1<r.
      \draw[draw=regimeAxis, line width=0.65bp] (-1,-1) -- (-1,2);
      \draw[draw=regimePurple, line width=1.7bp] (-4,-4) -- (-1,-1);
      \draw[draw=regimeRed, line width=1.7bp] (-1,-1) -- (2,-1);
    
      % Labels of the open regions.
      \node[regime label] at (-2.54,0.90) {Dense};
      \node[regime label] at (0.50,0.75) {Mixture: dense and\\divergent condensation};
      \node[regime label] at (0.50,-2.75) {Convergent\\condensation};
    
      % Labels on the boundaries.
      \node[regime box, rotate=45, minimum width=4.225cm, minimum height=1.29cm] at (-2.65,-2.65)
        {Mixture: dense and\\convergent condensation\\$r=s<-1$};
      \node[regime box, minimum width=4.065cm, minimum height=0.902cm] at (0.50,-1)
        {Divergent condensation\\$s=-1<r$};
    
      % Callout for point r=s=-1.
      \node[regime box, minimum width=3.11cm, minimum height=1.675cm] (critical-label) at (-2.75,-0.55)
        {Mixture: dense and\\divergent condensation\\$r=s=-1$};
      \draw[draw=regimePurple, line width=0.55bp, ->] (critical-label.east) -- (-1.10,-0.98);
      \fill[white] (-1,-1) circle[radius=4bp];
      \fill[regimePurple] (-1,-1) circle[radius=2.8bp];
    
      % Axes sit just outside the shaded plotting window.
      \draw[regime axis, ->] (-4.12,-4.12) -- (2.34,-4.12);
      \draw[regime axis, ->] (-4.12,-4.12) -- (-4.12,2.34);
      \node[regime tick label, right=4pt] at (2.34,-4.12) {$r$};
      \node[regime tick label, above=4pt] at (-4.12,2.34) {$s$};
    
      \foreach \v in {-3,-2,-1,0,1,2} {
        \draw[regime tick] (\v,-4.12) -- (\v,-4.19);
        \node[regime tick label, below=4pt] at (\v,-4.19) {$\v$};
        \draw[regime tick] (-4.12,\v) -- (-4.19,\v);
        \node[regime tick label, left=5pt] at (-4.19,\v) {$\v$};
      }
\end{tikzpicture}
\endgroup